\documentclass[a4paper,10pt]{article}
\usepackage{amsmath,amsthm,amssymb,vmargin}
\title{\vspace{25ex}
\vbox{\textbf{Joint distribution of the process\\ and its sojourn time
for pseudo-processes\\ governed by high-order heat equation}}}
\author{Valentina \textsc{Cammarota}\footnote{
Dipartimento di Statistica, Probabilit\`a e Statistiche Applicate,
\textsc{University of Rome `La Sapienza'}, P.le A. Moro 5, 00185 Rome,
\textsc{Italy}. \textit{E-mail address}: \texttt{valentina.cammarota@uniroma1.it}}$\,$
and Aim\'e \textsc{Lachal}\footnote{
P\^ole de Math\'ematiques/Institut Camille Jordan/CNRS UMR5208,
B\^at. L. de Vinci,
\textsc{Institut National des Sciences Appliqu\'ees de Lyon},
20 av. A. Einstein, 69621 Villeurbanne Cedex, \textsc{France}.
\textit{E-mail address:} \texttt{aime.lachal@insa-lyon.fr},
\textit{Web page:} \texttt{http://maths.insa-lyon.fr/$\mbox{}^{\sim}$lachal}}
}
\date{}

\setmarginsrb{25mm}{15mm}{25mm}{20mm}{0mm}{10mm}{0mm}{10mm}

\newtheorem{teo}{Theorem}[section]
\newtheorem{prop}{Proposition}[section]
\newtheorem{defi}{Definition}[section]
\newtheorem{lem}{Lemma}[section]
\theoremstyle{definition}
\newtheorem{rem}{Remark}[section]
\newtheorem{ex}{Example}[section]
\numberwithin{equation}{section}

\newcommand{\Dim}{\noindent\textsc{Proof\\}}
\newcommand{\EndDim}{ $\blacksquare$ \\}
\newcommand{\ind}{1\hspace{-.27em}\mbox{\rm l}}
\newcommand{\inde}{1\hspace{-.23em}\mathrm{l}}
\newcommand{\lqn}[1]{\noalign{\noindent $\displaystyle{#1}$}}

\begin{document}
\maketitle

\begin{abstract}
Consider the high-order heat-type equation $\partial u/\partial
t=\pm\,\partial^N u/\partial x^N$ for an integer $N>2$ and introduce
the related Markov pseudo-process $(X(t))_{t\ge 0}$. In this paper,
we study the sojourn time $T(t)$ in the interval $[0,+\infty)$ up to
a fixed time $t$ for this pseudo-process. We provide explicit expressions
for the joint distribution of the couple $(T(t),X(t))$.
\end{abstract}

\begin{small}
\noindent \textbf{Keywords:}
pseudo-process, joint distribution of the process and its sojourn time,
Spitzer's identity.

\noindent\textbf{AMS 2000 Subject Classification:}
Primary 60G20; Secondary 60J25, 60K35, 60J05.
\end{small}

\newpage
\section{Introduction}

Let $N$ be an integer equal or greater than $2$ and $\kappa_{\!_{ N}}=(-1)^{1+N/2}$
if $N$ is even, $\kappa_{\!_{ N}}=\pm1$ if $N$ is odd. Consider the heat-type
equation of order $N$:
\begin{equation}
\frac{\partial u}{\partial t}=\kappa_{\!_{ N}}\,\frac{\partial^N u}{\partial x^N}.
\label{EDP}
\end{equation}
For $N=2$, this equation is the classical normalized heat equation and its
relationship with linear Brownian motion is of the most well-known. For $N>2$,
it is known that no ordinary stochastic process can be associated with this
equation.
Nevertheless a Markov ``pseudo-process'' can be constructed by imitating the
case $N=2$. This pseudo-process, $X=(X(t))_{t\ge 0}$ say, is driven by a signed
measure as follows. Let $p(t;x)$ denote the elementary solution of
Eq.~(\ref{EDP}), that is, $p$ solves~(\ref{EDP}) with the initial
condition $p(0;x)=\delta(x)$.
This solution is characterized by its Fourier transform (see, e.g., \cite{2007})
$$
\int_{-\infty}^{+\infty} e^{i\mu x}\,p(t;x)\,\mathrm{d}x
=e^{\kappa_{\!_{ N}}t(-i\mu)^N}.
$$
The function $p$ is real, not always positive and its total mass is equal to one:
$$
\int_{-\infty}^{+\infty} p(t;x)\,\mathrm{d}x=1.
$$
Moreover, its total absolute value mass $\rho$ exceeds one:
$$
\rho=\int_{-\infty}^{+\infty} |p(t;x)|\,\mathrm{d}x>1.
$$
In fact, if $N$ is even, $p$ is symmetric and $\rho<+\infty$,
and if $N$ is odd, $\rho=+\infty$.
The signed function $p$ is interpreted as the pseudo-probability for
$X$ to lie at a certain location at a certain time. More precisely,
for any time $t>0$ and any locations $x,y\in\mathbb{R}$, one defines
$$
\mathbb{P}\{X(t)\in \mathrm{d}y|X(0)=x\}/\mathrm{d}y=p(t;x-y).
$$
Roughly speaking, the distribution of the pseudo-process $X$ is defined through its
finite-dimensional distributions according to the Markov rule:
for any $n>1$, any times $t_1,\dots,t_n$ such that $0<t_1<\dots<t_n$
and any locations $x,y_1,\dots,y_n\in\mathbb{R}$,
$$
\mathbb{P}\{X(t_1)\in \mathrm{d}y_1,\dots X(t_n)\in \mathrm{d}y_n|X(0)=x\}
/\mathrm{d}y_1\dots\mathrm{d}y_n=\prod_{i=1}^n p(t_i-t_{i-1};y_{i-1}-y_i)
$$
where $t_0=0$ and $y_0=x$.

This pseudo-process has been studied by several authors: see the references
\cite{bho} to \cite{bor} and the references \cite{hoch} to \cite{ors}.

Now, we consider the sojourn time of $X$ in the interval $[0,+\infty)$ up to
a fixed time $t$:
$$
T(t)=\int_0^t \ind_{[0,+\infty)} (X(s))\,\mathrm{d}s.
$$
The computation of the pseudo-distribution of $T(t)$ has been done by
Beghin, Hochberg, Nikitin, Orsingher and Ragozina in some particular cases
(see \cite{bho,bor,ho,no,ors}), and by Krylov and the second author
in more general cases (see \cite{kry,2003}).

The method adopted therein is the use of the Feynman-Kac functional which leads
to certain differential equations. We point out that the pseudo-distribution
of $T(t)$ is actually a genuine probability distribution and in the case where
$N$ is even, $T(t)$ obeys the famous Paul L\'evy's arcsine law, that is
$$
\mathbb{P}\{T(t)\in\mathrm{ds}\}/\mathrm{ds}=\frac{\ind_{(0,t)}(s)}{\pi\sqrt{s(t-s)}}.
$$
We also mention
that the sojourn time of $X$ in a small interval $(-\varepsilon,\varepsilon)$
is used in~\cite{bo} to define a local time for $X$ at 0.
The evaluation of the pseudo-distribution of the sojourn time $T(t)$ together
with the up-to-date value of the pseudo-process, $X(t)$, has been tackled only in
the particular cases $N=3$ and $N=4$ by Beghin, Hochberg, Orsingher and Ragozina
(see \cite{bho,bor}).
Their results have been obtained by solving certain differential equations
leading to some linear systems.
In~\cite{bho,bor,2003}, the Laplace transform of the sojourn time
serves as an intermediate tool for computing the distribution of
the up-to-date maximum of $X$.

In this paper, our aim is to derive the joint pseudo-distribution of the
couple $(T(t),X(t))$ for any integer $N$. Since the Feynman-Kac approach
used in~\cite{bho,bor} leads to very cumbersome calculations, we employ an
alternative method based on Spitzer's identity. Since the pseudo-process $X$ is
properly defined only in the case where $N$ is an even integer, the results we
obtain are valid in this case. Throughout the paper, we shall then
assume that $N$ is even. Nevertheless, we formally perform all computations
also in the case where $N$ is odd, even if they are not justified.

The paper is organized as follows.
\begin{itemize}
\item
In Section~\ref{section-settings}, we write down the settings
that will be used. Actually, the pseudo-process $X$ is not well defined on the
whole half-line $[0,+\infty)$. It is properly defined on dyadic times
$k/2^n$, $k,n\in\mathbb{N}$. So, we introduce ad-hoc definitions for $X(t)$ and
$T(t)$ as well as for some related pseudo-expectations.
For instance, we shall give a meaning to the quantity
$$
E(\lambda,\mu,\nu)=\mathbb{E}\!\left[\int_0^{\infty}
e^{-\lambda t+i \mu X(t)-\nu T(t)}\,\mathrm{d}t\right]
$$
which is interpreted as the $3$-parameters Laplace-Fourier transform of
$(T(t),X(t))$. We also recall in this part some algebraic known results.

\item
In Section~\ref{section-evaluation}, we explicitly compute $E(\lambda,\mu,\nu)$
with the help of Spitzer's identity. This is Theorem~\ref{theorem}.

\item
Sections~\ref{section-inverting-mu}, \ref{section-inverting-nu}
and~\ref{section-inverting-lambda} are devoted to successively inverting the
Laplace-Fourier transform with respect to $\mu$, $\nu$ and $\lambda$ respectively.
More precisely, in Section~\ref{section-inverting-mu}, we perform the inversion
with respect to $\mu$; this yields Theorem~\ref{theorem-wrt.mu}. Next, we perform the inversion
with respect to $\nu$ which gives Theorems~\ref{theorem-wrt.nu+1} and~\ref{theorem-wrt.nu+2}.
Finally, we carry out the inversion with respect to $\lambda$ and the main
results of this paper are Theorems~\ref{theorem-wrt.lambda+1} and~\ref{theorem-wrt.lambda+2}.
In each section, we examine the particular cases $N=2$ (case of
rescaled Brownian motion), $N=3$ (case of an asymmetric pseudo-process)
and $N=4$ (case of the biharmonic pseudo-process). Moreover, our results
recover several known formulas concerning the marginal distribution of $T(t)$.

\item
The final appendix (Section~\ref{section-appendix}) contains a discussion
on Spitzer's identity as well as some technical computations.
\end{itemize}

\section{Settings}\label{section-settings}

\subsection{A first list of settings}

In this part, we introduce for each integer $n$ a step-process $X_n$ coinciding
with the pseudo-process $X$ on the times $k/2^n$, $k\in\mathbb{N}$.
Fix $n\in\mathbb{N}$. Set, for any $k\in\mathbb{N}$, $X_{k,n}=X(k/2^n)$ and
for any $t\in [k/2^n,(k+1)/2^n),$ $X(t)=X_{k,n}$. We can write globally
$$
X_n(t)=\sum_{k=0}^{\infty} X_{k,n} \ind_{[k/2^n,(k+1)/2^n)}(t).
$$
Now, we recall from~\cite{2007} the definitions of tame functions, functions of discrete
observations, and admissible functions associated with the pseudo-process $X$.
They were introduced by Nishioka~\cite{nish2} in the case $N=4$.
%
\begin{defi}
Fix $n\in\mathbb{N}$. A tame function for $X$ is a function of a finite number $k$
of observations of the pseudo-process $X$ at times $j/2^n$, $1\le j\le k$, that is
a quantity of the form $F_{k,n}=F(X(1/2^n),\dots,X(k/2^n))$
for a certain $k$ and a certain bounded Borel function
$F:\mathbb{R}^k\longrightarrow \mathbb{C}$.
The ``expectation'' of $F_{k,n}$ is defined as
$$
\mathbb{E}(F_{k,n})=\int\dots\int_{\mathbb{R}^k} F(x_1,\dots,x_k)
\,p(1/2^n;x-x_1)\dots p(1/2^n;x_{k-1}-x_k)\,\mathrm{d}x_1\dots \mathrm{d}x_k.
$$
\end{defi}
%
\begin{defi}\label{def2}
Fix $n\in\mathbb{N}$. A function of the discrete observations of $X$ at times
$k/2^n$, $k\ge 1$, is a convergent series of tame functions:
$F_{X_n}=\sum_{k=1}^{\infty} F_{k,n}$ where $F_{k,n}$ is a tame function for all
$k\ge 1$. Assuming the series $\sum_{k=1}^{\infty} |\mathbb{E}(F_{k,n})|$
convergent, the ``expectation'' of $F_{X_n}$ is defined as
$$
\mathbb{E}(F_{X_n})=\sum_{k=1}^{\infty} \mathbb{E}(F_{k,n}).
$$
\end{defi}
%
\begin{defi}\label{def3}
An admissible function is a functional $F_X$ of the pseudo-process $X$ which is
the limit of a sequence $(F_{X_n})_{n\in\mathbb{N}}$ of functions of discrete
observations of $X$: $F_X=\lim_{n\to\infty} F_{X_n},$ such that the sequence
$(\mathbb{E}(F_{X_n}))_{n\in\mathbb{N}}$ is convergent.
The ``expectation'' of $F_X$ is defined as
$$
\mathbb{E}(F_X)=\lim_{n\to\infty}\mathbb{E}(F_{X_n}).
$$
\end{defi}
%

In this paper, we are concerned with the sojourn time of $X$ in $[0,+\infty)$:
$$
T(t)=\int_0^t \ind_{[0,+\infty)}(X(s))\,\mathrm{d}s.
$$
In order to give a proper meaning to this quantity, we introduce the similar
object related to $X_n$:
$$
T_n(t)=\int_0^t \ind_{[0,+\infty)}(X_n(s))\,\mathrm{d}s.
$$
For determining the distribution of $T_n(t)$, we compute its $3$-parameters
Laplace-Fourier transform:
$$
E_n(\lambda,\mu,\nu)=\mathbb{E}\!\left[\int_0^{\infty}
e^{-\lambda t+i\mu X_n(t)-\nu T_n(t)}\,\mathrm{d}t\right]\!.
$$
In Section~\ref{section-evaluation}, we prove that the sequence
$(E_n(\lambda,\mu,\nu))_{n\in\mathbb{N}}$ is convergent and we compute its limit:
$$
\lim_{n\to\infty} E_n(\lambda,\mu,\nu)=E(\lambda,\mu,\nu).
$$
Formally, $E(\lambda,\mu,\nu)$ is interpreted as
$$
E(\lambda,\mu,\nu)=\mathbb{E}\!\left[\int_0^{\infty}
e^{-\lambda t+i\mu X(t)-\nu T(t)}\,\mathrm{d}t\right]
$$
where the quantity $\int_0^{\infty}e^{-\lambda t+i\mu X(t)-\nu T(t)}\,\mathrm{d}t$
is an admissible function of $X$.
This computation is performed with the aid of Spitzer's identity.
This latter concerns the classical random walk. Nevertheless, since it hinges
on combinatorial arguments, it can be applied to the context of pseudo-processes.
We clarify this point in Section~\ref{section-evaluation}.

\subsection{A second list of settings}

We introduce some algebraic settings.
Let $\theta_i$, $1 \le i \le N$, be the $N^{\mathrm{th}}$ roots of $\kappa_{\!_{ N}}$ and
$$
J=\{i \in \{1,\dots, N\}:\,\Re{\theta_i}>0\},\qquad
K=\{i \in \{1,\dots, N\}:\,\Re{\theta_i}<0\}.
$$
Of course, the cardinalities of $J$ and $K$ sum to $N$: $\#J+\#K=N$.
We state several results related to the $\theta_i$'s which are proved
in~\cite{2003,2007}. We have the elementary equalities
\begin{equation}
\sum_{j \in J} \theta_j+\sum_{k \in K} \theta_k=\sum_{i=1}^N \theta_i=0,
\qquad \left(\vphantom{\prod_{k}}\right.\!\prod_{j\in J} \theta_j
\!\left.\vphantom{\prod_{k}}\right) \!\! \left(\prod_{k\in K} \theta_k\right)
=\prod_{i=1}^N \theta_i=(-1)^{N-1}\kappa_{\!_{ N}}
\label{set4}
\end{equation}
and
\begin{equation}
\prod_{i=1}^{N} (x-\theta_i)=\prod_{i=1}^{N} (x-\bar{\theta}_i)=x^N-\kappa_{\!_{ N}}.
\label{set6}
\end{equation}
Moreover, from formula (5.10) in \cite{2007},
\begin{equation}
\prod_{k\in K}(x-\theta_k)=\sum_{\ell=0}^{\# K} (-1)^{\ell}
\sigma_{\ell} \,x^{\# K-\ell},
\label{set10}
\end{equation}
where $\sigma_{\ell}=\sum_{\substack{k_1<\dots<k_\ell\\ k_1,\dots, k_\ell\in K}}
\theta_{k_1} \dots \theta_{k_\ell}.$
We have by Lemma~11 in \cite{2003}
\begin{equation}
\sum_{j \in J} \theta_j \prod_{i \in J\setminus \{j\}}
\frac{\theta_i x -\theta_j}{\theta_i-\theta_j}=\sum_{j \in J} \theta_j
=-\sum_{k \in K} \theta_k=
\begin{cases}
\displaystyle{\frac{1}{\sin \frac{\pi}{N}}}& \mbox{if } N \mbox{\,is\, even},
\\[2ex]
\displaystyle{\frac{1}{2 \sin \frac{\pi}{2N}}=\frac{\cos \frac{\pi}{2 N}}{\sin \frac{\pi}{N}}}
& \mbox{if } N \mbox{\,is\, odd}.
\end{cases}
\label{set7}
\end{equation}

Set $A_j=\prod_{i \in J\setminus  \{j\}} \frac{\theta_i}{\theta_i-\theta_j}$
for $j \in J$, and $B_k=\prod_{i \in K\setminus  \{k\}}
\frac{\theta_i}{\theta_i-\theta_k}$ for $k \in K$. The $A_j$'s and $B_k$'s
solve a Vandermonde system: we have
\begin{gather}
\sum_{j \in J} A_j=\sum_{k \in K} B_k=1
\nonumber\\[-3ex]
\label{set13}
\\[-1ex]
\sum_{j \in J} A_j\theta_j^m=0 \mbox{ for } 1\le m\le \# J-1,\qquad
\sum_{k \in K}  B_k\theta_k^m=0 \mbox{ for } 1\le m\le \# K-1.
\nonumber
\end{gather}
Observing that $1/\theta_j=\bar{\theta}_j$ for $j\in J$, that
$\{\theta_j,j\in J\}=\{\bar{\theta}_j,j\in J\}$ and similarly for the $\theta_k$'s,
$k\in K$, formula (2.11) in \cite{2007} gives
\begin{align}
\sum_{j \in J} \frac{A_j\theta_j}{\theta_j-x}
&=
\sum_{j \in J}\frac{A_j}{1-\bar{\theta}_j x}= \frac{1}{\prod_{j \in J}(1-\theta_j x)}
=-\kappa_{\!_{ N}}\,\frac{\prod_{k \in K} (1-\theta_k x)}{x^N-\kappa_{\!_{ N}}}
=-\kappa_{\!_{ N}}\,\frac{\prod_{k \in K} (1-\bar{\theta}_k x)}{x^N-\kappa_{\!_{ N}}},
\nonumber\\
\label{set5}
\\
\sum_{k \in K} \frac{B_k\theta_k}{\theta_k-x}
&=
\sum_{k \in K} \frac{B_k}{1-\bar{\theta}_k x}= \frac{1}{\prod_{k \in K}(1-\theta_k x)}
=-\kappa_{\!_{ N}}\,\frac{\prod_{j \in J} (1-\theta_j x)}{x^N-\kappa_{\!_{ N}}}
=-\kappa_{\!_{ N}}\,\frac{\prod_{j \in J} (1-\bar{\theta}_j x)}{x^N-\kappa_{\!_{ N}}}.
\nonumber
\end{align}
In particular,
\begin{equation}
\sum_{j \in J}\frac{A_j\theta_j}{\theta_j-\theta_k}=\frac{1}{N B_k},
\hspace{1cm} \sum_{k \in K}\frac{B_k\theta_k}{\theta_k-\theta_j}
=\frac{1}{N A_j}.
\label{set11}
\end{equation}

Set, for any $m\in\mathbb{Z}$, $\alpha_{m}=\sum_{j \in J} A_j\theta_j^m$ and
$\beta_m=\sum_{k \in K} B_k\theta_k^m$. We have, by formula (2.11)  of~\cite{2007},
$\beta_{\#K}=(-1)^{\#K -1}\prod_{k \in K} \theta_k$. Moreover,
$\beta_{\#K+1}=(-1)^{\#K -1}\left(\prod_{k \in K} \theta_k\right)\!\left(\sum_{k \in K} \theta_k\right)$.
The proof of this claim is postponed to Lemma~\ref{lemma-vdm} in the appendix. We sum up
this information and~(\ref{set13}) into
\begin{equation}
\beta_m=\begin{cases}
1 & \mbox{if } m=0,
\\
0 & \mbox{if } 1 \le m \le \#K-1,
\\
(-1)^{\#K-1} \prod_{k \in K} \theta_k & \mbox{if } m=\# K ,
\\
(-1)^{\#K-1} \left(\prod_{k \in K} \theta_k\right)\!
\left(\sum_{k \in K} \theta_k\right) &  \mbox{if } m=\# K+1,
\\
\kappa_{\!_{ N}} & \mbox{if } m=N.
\end{cases}
\label{set26}
\end{equation}
We also have
$$
\alpha_{-m}=\sum_{j \in J} \frac{A_j}{\theta_j^m}
=\kappa_{\!_{ N}} \sum_{j \in J}   A_j\theta_j^{N-m}=\kappa_{\!_{ N}} \alpha_{N-m}
$$
and then
\begin{equation}
\alpha_{-m}=\begin{cases}
1 & \mbox{if } m=0,
\\
\kappa_{\!_{ N}} (-1)^{\# J -1} \left(\vphantom{\prod_{\in}}\right.\!\prod_{j \in J}
\theta_j\left.\vphantom{\prod_{\in}}\!\right)\!
\left(\vphantom{\prod_{\in}}\right.\!\sum_{j \in J} \theta_j
\left.\vphantom{\prod_{\in}}\!\right) & \mbox{if } m=\# K-1,
\\
\kappa_{\!_{ N}} (-1)^{\#J -1}  \prod_{j \in J} \theta_j & \mbox{if } m=\#K,
\\
0 & \mbox{if } \#K+1 \le m \le N-1,
\\
\kappa_{\!_{ N}} & \mbox{if } m=N.
\end{cases}
\label{set14}
\end{equation}
In particular, by~(\ref{set4}),
\begin{equation}
\alpha_0\beta_0=\alpha_{-N}\beta_N=1,\quad
\alpha_{-\#K}\beta_{\#K}=-1,\quad
\alpha_{-\#K}\beta_{\#K+1}=\sum_{j \in J} \theta_j,\quad
\alpha_{1-\#K}\beta_{\#K}=\sum_{k \in K} \theta_k.
\label{set15}
\end{equation}
With $\sigma_0=1$, $\sigma_{\#K-1}=\left(\prod_{k \in K} \theta_k\right)\!
\left(\sum_{k \in K} \bar{\theta}_k\right)$ and $\sigma_{\#K}=\prod_{k \in K} \theta_k$,
we also have
\begin{equation}
\bar{\sigma}_0\beta_0=1,\quad
\bar{\sigma}_{\#K-1}\beta_{\#K}=\bar{\sigma}_{\#K}\beta_{\#K+1}
=(-1)^{\#K-1}\sum_{k \in K} \theta_k,\quad
\bar{\sigma}_{\#K}\beta_{\#K}=(-1)^{\#K-1}.
\label{set16}
\end{equation}

Concerning the kernel $p$, we have from Proposition~1 in \cite{2003}
\begin{equation}
p(t;0)=\begin{cases}
\displaystyle{\frac{\Gamma\!\left(\vphantom{\frac aN}\right.\!\!\frac{1}{N}
\!\!\left.\vphantom{\frac aN}\right)}{N \pi t^{1/N}}}
& \mbox{if $N$ is even},
\\[2ex]
\displaystyle{\frac{\Gamma\!\left(\vphantom{\frac aN}\right.\!\!\frac{1}{N}
\!\!\left.\vphantom{\frac aN}\right) \cos\!\left(\vphantom{\frac aN}\right.\!\!
\frac{\pi}{2 N}\!\!\left.\vphantom{\frac aN}\right)}{N \pi t^{1/N}}}
& \mbox{if $N$ is odd}.
\end{cases}
\label{p-at-zero}
\end{equation}
Proposition~3 in \cite{2003} states
\begin{equation}
\mathbb{P}\{X(t)\ge 0\}=\int_0^{\infty} p(t;-\xi)  \,\mathrm{d} \xi=\frac{\# J}{N},
\qquad \mathbb{P}\{X(t)\le 0\}=\int_{-\infty}^0 p(t;-\xi)  \,\mathrm{d} \xi=\frac{\# K}{N}
\label{set3}
\end{equation}
and formulas (4.7) and (4.8) in \cite{2007} yield, for $\lambda>0$ and $\mu\in\mathbb{R}$,
\begin{align}
\int_0^{\infty} \frac{e^{-\lambda t}}{t} \,\mathrm{d}t
\int_{-\infty}^0\left(e^{i\mu\xi}-1\right) p(t;-\xi) \,\mathrm{d}\xi
&=
\log \!\left(\,\prod_{k \in K} \frac{\!\sqrt[N]{\lambda}}{\!\sqrt[N]{\lambda}
-i\mu\theta_k}\right)\!,
\nonumber\\[-1ex]
\label{set1}\\[-1ex]
\int_0^{\infty} \frac{e^{-\lambda t}}{t} \,\mathrm{d}t
\int_0^{\infty}\left(e^{i\mu\xi}-1\right) p(t;-\xi) \,\mathrm{d}\xi
&=
\log \!\left(\vphantom{\prod_{\in}}\right.\!\prod_{j \in J} \frac{\!\sqrt[N]{\lambda}}{\!\sqrt[N]{\lambda}
-i\mu\theta_j}\!\left.\vphantom{\prod_{\in}}\right)\!.
\nonumber
\end{align}

Let us introduce, for $m \le N-1$ and $x \ge 0$,
\begin{equation}
I_{j,m}(\tau;x)=\frac{N i}{2 \pi} \left(e^{-i \frac{m}{N}\pi}\!
\int_0^{\infty} \xi^{N-m-1} e^{-\tau\xi^N-\theta_j e^{i \frac{\pi}{N}} x \xi}
\,\mathrm{d}\xi
-e^{i \frac{m}{N}\pi} \!\int_0^{\infty} \xi^{N-m-1} e^{-\tau\xi^N-\theta_j
e^{-i \frac{\pi}{N}} x \xi} \,\mathrm{d}\xi\right)\!.
\label{set18}
\end{equation}
Formula (5.13) in \cite{2007} gives, for $0 \le m \le N-1$ and $x\ge 0$,
\begin{equation}
\int_0^{\infty} e^{-\lambda \tau} I_{j,m}(\tau;x)\,\mathrm{d}\tau
=\lambda^{-\frac{m}{N}}e^{-\theta_j\!\!\sqrt[N]{\lambda}\,x}.
\label{set24}
\end{equation}

%
\begin{ex}\label{example1}
Case $N=2$: we have $\kappa_2=+1$. This is the case of rescaled Brownian motion.
The square roots of $\kappa_2$ are $\theta_1=1$, $\theta_2=-1$
and then $J=\{1\}$, K=\{2\}, $A_1=B_2=1$, $\alpha_0=\alpha_{-1}=1$,
$\beta_0=1$, $\beta_{-1}=-1$. Moreover,
$$
I_{1,0}(\tau;x)=\frac{i}{\pi} \left(\int_0^{\infty}\xi \,e^{-\tau \xi^2-ix \xi}
\,\mathrm{d}\xi-\int_0^{\infty}\xi \,e^{-\tau \xi^2+ix \xi}\,\mathrm{d}\xi\right)\!.
$$
The function $I_{1,0}$ can be simplified. In fact, we have
\begin{align*}
I_{1,0}(\tau;x)
&=
\frac{i}{\pi} \int_{-\infty}^{\infty} \xi \,e^{-\tau \xi^2-ix \xi} \,\mathrm{d}\xi
=\frac{i}{\pi} \,e^{-\frac{x^2}{4 \tau}} \int_{-\infty}^{\infty} \xi\,
e^{-\tau(\xi+\frac{ix}{2 \tau})^2}\,\mathrm{d}\xi
\\
&=
\frac{i}{\pi} \,e^{-\frac{x^2}{4 \tau}} \int_{-\infty}^{\infty}
\left(\xi-\frac{i x}{2 \tau}\right) e^{-\tau\xi^2}\,\mathrm{d}\xi
=\frac{x \,e^{-\frac{x^2}{4 \tau}}}{2 \pi \tau} \int_{-\infty}^{\infty}e^{-\tau \xi^2} \,\mathrm{d}\xi
=\frac{x \,e^{-\frac{x^2}{4 \tau}}}{2 \pi \tau} \int_0^{\infty}\frac{e^{-\tau \xi}}{\sqrt{\xi}} \,\mathrm{d}\xi.
\end{align*}
Finally,
\begin{equation}
I_{1,0}(\tau;x)=\frac{x \,e^{-\frac{x^2}{4 \tau}}}{2 \sqrt \pi \,\tau^{3/2}}.
\label{case2I}
\end{equation}
\end{ex}
%
\begin{ex}\label{example2}
Case $N=3$.

\noindent $\bullet$ For $\kappa_3=+1$, the third roots of $\kappa_3$ are
$\theta_1=1$, $\theta_2=e^{i \frac{2\pi}{3}}$, $\theta_3=e^{-i \frac{2\pi}{3}}$,
and the settings read $J=\{1\}$, $K=\{2,3\}$, $A_1=1$,
$B_2=\frac{e^{-i\frac{\pi}{6}}}{\sqrt{3}}$, $B_3=\frac{e^{i\frac{\pi}{6}}}{\sqrt{3}}$,
$\alpha_0=\alpha_{-1}=\alpha_{-2}=1$, $\beta_0=1$, $\beta_{-1}=-1$. Moreover,
$$
I_{1,0}(\tau;x)=\frac{3i}{2\pi} \left(
\int_0^{\infty} \xi^2\,e^{-\tau\xi^3-e^{i\frac{\pi}{3}}x\xi}\,\mathrm{d}\xi
-\int_0^{\infty} \xi^2\,e^{-\tau\xi^3-e^{-i\frac{\pi}{3}}x\xi}\,\mathrm{d}\xi\right)\!.
$$

\noindent $\bullet$ For $\kappa_3=-1$, the third roots of $\kappa_3$ are
$\theta_1=e^{i \frac{\pi}{3}}$, $\theta_2=e^{-i \frac{\pi}{3}}$, $\theta_3=-1$.
The settings read $J=\{1,2\}$, $K=\{3\}$,
$A_1=\frac{e^{i \frac{\pi}{6}}}{\sqrt{3}}$, $A_2=\frac{e^{-i\frac{\pi}{6}}}{\sqrt{3}}$,
$B_3=1$, $\alpha_0=\alpha_{-1}=1$, $\beta_0=\beta_{-2}=1$, $\beta_{-1}=-1$. Moreover,
\begin{align*}
I_{1,1}(\tau;x)
&=
\frac{3i}{2\pi} \left(e^{-i\frac{\pi}{3}}
\int_0^{\infty} \xi\,e^{-\tau\xi^3-e^{i\frac{2\pi}{3}}x\xi}\,\mathrm{d}\xi
-e^{i\frac{\pi}{3}} \int_0^{\infty} \xi\,e^{-\tau\xi^3-x\xi}\,\mathrm{d}\xi\right)\!,
\\
I_{2,1}(\tau;x)
&=
\frac{3i}{2\pi} \left(e^{-i\frac{\pi}{3}}
\int_0^{\infty} \xi\,e^{-\tau\xi^3-x\xi}\,\mathrm{d}\xi
-e^{i\frac{\pi}{3}} \int_0^{\infty} \xi\,e^{-\tau\xi^3-e^{-i\frac{2\pi}{3}}x\xi}
\,\mathrm{d}\xi\right)\!.
\end{align*}
Actually, the three functions $I_{1,0}$, $I_{1,1}$ and $I_{2,1}$ can be
expressed by mean of the Airy function $\mathrm{Hi}$ defined as
$\mathrm{Hi}(z)=\frac{1}{\pi}\int_0^{\infty} e^{-\frac{\xi^3}{3}+z\xi}\,\mathrm{d}\xi$
(see, e.g., \cite[Chap.~10.4]{as}).
Indeed, we easily have by a change of variables, differentiation and integration by parts,
for $\tau>0$ and $z\in\mathbb{C}$,
\begin{align*}
\int_0^{\infty} e^{-\tau\xi^3+z\xi}\,\mathrm{d}\xi
&=
\frac{\pi}{(3\tau)^{4/3}}\,\mathrm{Hi}\!\left(\frac{z}{\sqrt[3]{3\tau}}\right)\!,
\\
\int_0^{\infty} \xi\,e^{-\tau\xi^3+z\xi}\,\mathrm{d}\xi
&=
\frac{\pi}{(3\tau)^{2/3}}\,\mathrm{Hi}'\!\left(\frac{z}{\sqrt[3]{3\tau}}\right)\!,
\\
\int_0^{\infty} \xi^2 \,e^{-\tau\xi^3+z\xi}\,\mathrm{d}\xi
&=
\frac{\pi z}{(3\tau)^{4/3}}\,\mathrm{Hi}\!\left(\frac{z}{\sqrt[3]{3\tau}}\right)+\frac{1}{3\tau}.
\end{align*}
Therefore,
\begin{align}
I_{1,0}(\tau;x)
&=
\frac{x}{2\sqrt[3]{3} \,\tau^{4/3}} \left[\vphantom{\frac{z}{\sqrt{3}}}\right.\!
e^{i\frac{\pi}{6}} \mathrm{Hi}\!\left(\vphantom{\frac{z}{\sqrt{3}}}\right.\!\!
-\frac{e^{-i\frac{\pi}{3}}x}{\sqrt[3]{3\tau}} \!\left.\vphantom{\frac{z}{\sqrt{3}}} \right)
+e^{-i\frac{\pi}{6}} \mathrm{Hi}\!\left(\vphantom{\frac{z}{\sqrt{3}}}\right.\!\!
-\frac{e^{i\frac{\pi}{3}}x}{\sqrt[3]{3\tau}} \!\left.\vphantom{\frac{z}{\sqrt{3}}} \right)\!
\!\left.\vphantom{\frac{z}{\sqrt{3}}} \right]\!,
\label{case3I10}\\
I_{1,1}(\tau;x)
&=
\frac{\sqrt[3]3}{2\tau^{2/3}} \left[\vphantom{\frac{z}{\sqrt{3}}}\right.\!
e^{i\frac{\pi}{6}} \mathrm{Hi}'\!\left(\vphantom{\frac{z}{\sqrt{3}}}\right.\!\!
-\frac{e^{i\frac{2\pi}{3}}x}{\sqrt[3]{3\tau}} \!\left.\vphantom{\frac{z}{\sqrt{3}}} \right)
+e^{-i\frac{\pi}{6}} \mathrm{Hi}'\!\left(-\frac{x}{\sqrt[3]{3\tau}}\right)
\!\!\left.\vphantom{\frac{z}{\sqrt{3}}}\right]\!,
\label{case3I11}\\
I_{2,1}(\tau;x)
&=
\frac{\sqrt[3]3}{2\tau^{2/3}} \left[\vphantom{\frac{z}{\sqrt{3}}}\right.\!
e^{i\frac{\pi}{6}} \mathrm{Hi}'\!\left(-\frac{x}{\sqrt[3]{3\tau}}\right)
+e^{-i\frac{\pi}{6}} \mathrm{Hi}'\!\left(\vphantom{\frac{z}{\sqrt{3}}}\right.\!\!
-\frac{e^{-i\frac{2\pi}{3}}x}{\sqrt[3]{3\tau}}\!\left.\vphantom{\frac{z}{\sqrt{3}}}\right)
\!\!\left.\vphantom{\frac{z}{\sqrt{3}}} \right]\!.
\label{case3I21}
\end{align}

\end{ex}
%
\begin{ex}\label{example3}
Case $N=4$: we have $\kappa_4=-1$. This is the case of the biharmonic pseudo-process.
The fourth roots of $\kappa_4$ are $\theta_1=e^{-i \frac{\pi}{4}}$, $\theta_2=e^{i \frac{\pi}{4}}$,
$\theta_3=e^{i\frac{3\pi}{4}}$, $\theta_4=e^{-i\frac{3\pi}{4}}$ and the
notations read in this case $J=\{1,2\}$, $K=\{3,4\}$,
$A_1=B_3=\frac{e^{-i \frac{\pi}{4}}}{\sqrt{2}}$,
$A_2=B_4=\frac{e^{i\frac{\pi}{4}}}{\sqrt{2}}$, $\alpha_0=\alpha_{-2}=1$,
$\alpha_{-1}=\sqrt 2$, $\beta_0=\beta_{-2}=1$, $\beta_{-1}=-\sqrt 2$. Moreover,
\begin{align}
I_{1,1}(\tau;x)
&=
\frac{2}{\pi} \left(e^{i\frac{\pi}{4}}
\int_0^{\infty} \xi^2\,e^{-\tau\xi^4-x\xi}\,\mathrm{d}\xi
+e^{-i\frac{\pi}{4}} \int_0^{\infty} \xi^2\,e^{-\tau\xi^4+ix\xi}\,\mathrm{d}\xi\right)\!,
\nonumber\\[-1ex]
\label{case4I}\\[-1ex]
I_{2,1}(\tau;x)
&=
\frac{2}{\pi} \left(e^{i\frac{\pi}{4}}
\int_0^{\infty} \xi^2\,e^{-\tau\xi^4-ix\xi}\,\mathrm{d}\xi
+e^{-i\frac{\pi}{4}} \int_0^{\infty} \xi^2\,e^{-\tau\xi^4-x\xi}
\,\mathrm{d}\xi\right)\!.
\nonumber\end{align}
\end{ex}

\section{Evaluation of $E(\lambda,\mu,\nu)$}\label{section-evaluation}

The goal of this section is to evaluate the limit
$E(\lambda,\mu,\nu)=\lim_{n \to \infty} E_n(\lambda,\mu,\nu)$.
We write $E_n(\lambda,\mu,\nu)=\mathbb{E}[F_n(\lambda,\mu,\nu)]$ with
$$
F_n(\lambda,\mu,\nu)=\int_0^{\infty} e^{-\lambda t+i \mu X_n(t)-\nu T_n(t)} \,\mathrm{d}t.
$$
Let us rewrite the sojourn time $T_n(t)$ as follows:
\begin{align*}
T_n(t)
&=
\sum_{j=0}^{[2^n t]}\int_{j/2^n}^{(j+1)/2^n} \ind_{[0,+\infty)}(X_n(s)) \,\mathrm{d}s
-\int_t^{([2^n t]+1)/2^n} \ind_{[0,+\infty)}(X_n(s)) \,\mathrm{d}s
\\
&=
\sum_{j=0}^{[2^n t]}\int_{j/2^n}^{(j+1)/2^n} \ind_{[0,+\infty)}(X_{j,n}) \,\mathrm{d}s
-\int_t^{([2^n t]+1)/2^n} \ind_{[0,+\infty)}(X_{[2^n t],n}) \,\mathrm{d}s
\\
&=
\frac{1}{2^n} \sum_{j=0}^{[2^n t]} \ind_{[0,+\infty)}(X_{j,n})
+\left(t-\frac{[2^n t]+1}{2^n}\right)\ind_{[0,+\infty)}(X_{[2^n t],n}).
\end{align*}
Set $T_{0,n}=0$ and, for $k\ge 1$,
$$
T_{k,n}=\frac{1}{2^n} \sum_{j=1}^k \ind_{[0,+\infty)}(X_{j,n}).
$$
For $k\ge 0$ and $t\in[k/2^n,(k+1)/2^n)$, we see that
$$
T_n(t)=T_{k,n}+\left(t-\frac{k+1}{2^n}\right) \ind_{[0,+\infty)}(X_{k,n})+\frac{1}{2^n}.
$$
With this decomposition at hand, we can begin to compute $F_n(\lambda,\mu,\nu)$:
\begin{align*}
F_n(\lambda,\mu,\nu)
&=
\int_0^{\infty} e^{-\lambda t+i\mu X_n(t)-\nu T_n(t)}\,\mathrm{d}t
\\
&=
\sum_{k=0}^{\infty} \int_{k/2^n}^{(k+1)/2^n}
e^{-\lambda t+i\mu X_{k,n}-\nu T_{k,n}-\frac{\nu}{2^n}
+\nu(\frac{k+1}{2^n}-t)\inde_{[0,+\infty)}(X_{k,n})} \,\mathrm{d}t
\\
&=
e^{-\nu/2^n} \left(\sum_{k=0}^{\infty} \int_{k/2^n}^{(k+1)/2^n}
e^{-\lambda t+\nu(\frac{k+1}{2^n}-t)\inde_{[0,+\infty)}(X_{k,n})} \,\mathrm{d}t\right)
e^{i\mu X_{k,n}-\nu T_{k,n}}.
\end{align*}
The value of the above integral is
$$
\int_{k/2^n}^{(k+1)/2^n} e^{-\lambda t+\nu(\frac{k+1}{2^n}-t)
\inde_{[0,+\infty)}(X_{k,n})} \,\mathrm{d}t
=e^{-\lambda(k+1)/2^n}\,\frac{e^{[\lambda
+\nu\inde_{[0,+\infty)}(X_{k,n})]/2^n}-1}{\lambda+\nu\ind_{[0,+\infty)}(X_{k,n})}.
$$
Therefore,
\begin{align*}
F_n(\lambda,\mu,\nu)
&=
\frac{1-e^{-(\lambda+\nu)/2^n}}{\lambda+\nu} \sum_{k=0}^{\infty}
e^{-\lambda k/2^n+i\mu X_{k,n}-\nu T_{k,n}} \ind_{[0,+\infty)}(X_{k,n})
\\
&\hphantom{=\,}
+ e^{-\nu/2^n}\,\frac{1-e^{-\lambda/2^n}}{\lambda} \sum_{k=0}^{\infty}
e^{-\lambda k/2^n+i\mu X_{k,n}-\nu T_{k,n}} \ind_{(-\infty,0)}(X_{k,n}).
\end{align*}
Before applying the expectation to this last expression, we have to
check that it defines a function of discrete observations of the pseudo-process $X$
which satisfies the conditions of Definition~\ref{def2}.
This fact is stated in the proposition below.
%
\begin{prop}\label{proposition}
Suppose $N$ even and fix an integer $n$. For any complex $\lambda$ such
that $\Re(\lambda)>0$ and any $\nu>0$, the series
$\sum_{k=0}^{\infty} e^{-\lambda k/2^n} \mathbb{E}\!\left[e^{i \mu X_{k,n}
-\nu T_{k,n}}\ind_{[0,+\infty)}(X_{k,n})\right]$ and
$\sum_{k=0}^{\infty} e^{-\lambda k/2^n} \mathbb{E}\!\left[e^{i \mu X_{k,n}
-\nu T_{k,n}}\right.$ $\left.\ind_{(-\infty,0)}(X_{k,n})\right]$ are absolutely convergent
and their sums are given by
\begin{align*}
\sum_{k=0}^{\infty} e^{-\lambda k/2^n} \mathbb{E}\!\left[e^{i \mu X_{k,n}
-\nu T_{k,n}}\ind_{[0,+\infty)}(X_{k,n})\right]
&=
\frac{e^{\nu/2^n}-S_n^+(\lambda,\mu,\nu)}{e^{\nu/2^n}-1},
\\
\sum_{k=0}^{\infty} e^{-\lambda k/2^n} \mathbb{E}\!\left[e^{i \mu X_{k,n}
-\nu T_{k,n}}\ind_{(-\infty,0)}(X_{k,n})\right]
&=
\frac{e^{\nu/2^n}[S_n^-(\lambda,\mu,\nu)-1]}{e^{\nu/2^n}-1},
\end{align*}
where
\begin{align*}
S_n^+(\lambda,\mu,\nu)
&=
\exp\!\left(-\sum_{k=1}^{\infty} \left(1-e^{-\nu k/2^n}\right)
\frac{e^{-\lambda k/2^n}}{k}\,
\mathbb{E}\!\left[e^{i \mu X_{k,n}}\ind_{[0,+\infty)}(X_{k,n})\right]\right)\!,
\\[1ex]
S_n^-(\lambda,\mu,\nu)
&=
\exp\!\left(\,\sum_{k=1}^{\infty} \left(1-e^{-\nu k/2^n}\right)
\frac{e^{-\lambda k/2^n}}{k}\,
\mathbb{E}\!\left[e^{i \mu X_{k,n}}\ind_{(-\infty,0)}(X_{k,n})\right]\right)\!.
\end{align*}
\end{prop}
%
\Dim
\noindent \textbf{$\bullet$ Step 1.}
First, notice that for any $k\ge 1$, we have
\begin{align*}
\lqn{\left|\mathbb{E}\!\left[e^{i \mu X_{k,n}-\nu T_{k,n}}
\ind_{[0,+\infty)}(X_{k,n})\right]\right|}
&=
\left|\,\int \right.\!\!\dots\! \int_{\mathbb{R}^{k-1}\times [0,+\infty)}
e^{i \mu x_k-\frac{\nu}{2^n} \sum_{j=1}^k \inde_{[0,+\infty)}(x_j)}\,
\mathbb{P}\{X_{1,n}\in \mathrm{d}x_1,\dots,X_{k,n}\in \mathrm{d}x_k\}
\!\left.\vphantom{\int}\right|
\\
&=
\left|\,\int\right.\!\!\dots\! \int_{\mathbb{R}^{k-1}\times [0,+\infty)}
e^{i \mu x_k-\frac{\nu}{2^n}\sum_{j=1}^k \inde_{[0,+\infty)}(x_j)}\,
p\!\left(\frac{1}{2^n}; x_1\right) \prod_{j=1}^{k-1}
p\!\left(\frac{1}{2^n}; x_j-x_{j+1}\right) \mathrm{d}x_1 \dots \mathrm{d}x_k
\!\left.\vphantom{\int}\right|
\\
&\le
\int\!\dots\!\int_{\mathbb{R}^k} \left|\vphantom{\int}\right.\!
p\!\left(\frac{1}{2^n}; x_1\right)
\prod_{j=1}^{k-1} p\!\left(\frac{1}{2^n}; x_j-x_{j+1}\right)
\!\!\left.\vphantom{\int}\right| \mathrm{d}x_1 \dots \mathrm{d}x_k
\\
&=
\int\!\dots\!\int_{\mathbb{R}^k} \prod_{j=1}^k \left|\, p\!\left(\frac{1}{2^n};
y_j\right)\right| \mathrm{d}y_1 \dots \mathrm{d}y_k
=\prod_{j=1}^k \int_{-\infty}^{+\infty} \left|\, p\!\left(\frac{1}{2^n}; y_j\right)\right|\mathrm{d}y_j=\rho^k.
\end{align*}
Hence, we derive the following inequality:
$$
\sum_{k=1}^{\infty}\left| e^{-\lambda k/2^n} \mathbb{E}\!\left[e^{i \mu X_{k,n} -\nu T_{k,n}}
\ind_{[0,+\infty)}(X_{k,n})\right]\right| \le \sum_{k=1}^{\infty} \rho^k \left| e^{-\lambda k/2^n}\right|
=\frac{1}{1-\rho e^{-\Re(\lambda)/2^n}}.
$$
We can easily see that this bound holds true also when the factor
$\ind_{[0,+\infty)}(X_{k,n})$ is replaced by $\ind_{(-\infty,0)}(X_{k,n})$.
This shows that the two series of Proposition~\ref{proposition} are finite for
$\lambda \in \mathbb{C}$ such that $ \rho e^{-\Re(\lambda)/2^n} <1$,
that is $\Re(\lambda)>2^n \log \rho$.
\\

\noindent \textbf{$\bullet$ Step 2.}
For $\lambda \in \mathbb{C}$ such that $\Re(\lambda)>2^n \log \rho$,
the Spitzer's identity~(\ref{spitzer-identitybis}) (see Lemma~\ref{lemma-spitzer}
in the appendix) gives for the first series of Proposition~\ref{proposition}
\begin{align}
\lqn{\sum_{k=0}^{\infty} e^{-\lambda k/2^n} \mathbb{E}\!
\left[e^{i \mu X_{k,n} -\nu T_{k,n}}\ind_{[0,+\infty)}(X_{k,n})\right]}
&=
\frac{1}{e^{\nu/2^n}-1}\left[e^{\nu/2^n}-
\exp\!\left(-\sum_{k=1}^{\infty} \left(1-e^{-\nu k/2^n}\right)\frac{e^{-\lambda k/2^n}}{k}\,
\mathbb{E}\!\left[e^{i \mu X_{k,n}}\ind_{[0,+\infty)}(X_{k,n})\right]\right)\right]\!.
\label{series}
\end{align}
The right-hand side of~(\ref{series}) is an analytic continuation of the
Dirichlet series lying in the left-hand side of~(\ref{series}),
which is defined on the half-plane $\{\lambda \in \mathbb{C}: \Re(\lambda)>0\}$.
Moreover, for any $\varepsilon>0$, this continuation is bounded over
the half-plane $\{\lambda \in \mathbb{C}: \Re(\lambda)\ge \varepsilon\}$.
Indeed, we have
$$
\left| \mathbb{E}\!\left[e^{i \mu X_{k,n}}\ind_{[0,+\infty)}(X_{k,n})\right]\right|
= \left|\, \int_0^{+\infty}  e^{i \mu \xi} p\!\left(\frac{k}{2^n};-\xi\right)  \mathrm{d}\xi\right|
\le \int_0^{+\infty} \left|\, p\!\left(\frac{k}{2^n}; -\xi\right)\right| \mathrm{d}\xi <\rho
$$
and then
\begin{align*}
\lqn{\left|\,\exp\!\left(-\sum_{k=1}^{\infty} \left(1-e^{-\nu k/2^n}\right)\frac{e^{-\lambda k/2^n}}{k}\,
\mathbb{E}\!\left[e^{i \mu X_{k,n}} \ind_{[0,+\infty)}(X_{k,n})\right]\right)\right|}
&\le
\exp\!\left(\rho \sum_{k=1}^{\infty} \frac{e^{-\Re(\lambda) k/2^n}}{k}\right)
=\exp\!\left(-\rho \log (1-e^{-\Re(\lambda)/2^n})\right)
=\frac{1}{(1-e^{-\Re(\lambda)/2^n})^{\rho}}.
\end{align*}
Therefore, if $\Re(\lambda)\ge\varepsilon$,
\begin{align*}
\left|\,\exp\!\left(-\sum_{k=1}^{\infty} \left(1-e^{-\nu k/2^n}\right)\frac{e^{-\lambda k/2^n}}{k}\,
\mathbb{E}\!\left[e^{i \mu X_{k,n}} \ind_{[0,+\infty)}(X_{k,n})\right]\right)\right|
&\le
\frac{1}{(1-e^{-\varepsilon/2^n})^{\rho}}.
\end{align*}
This proves that the left-hand side of this last inequality is bounded
for $\Re(\lambda) \ge \varepsilon$.
By a lemma of Bohr (\cite{bohr}), we deduce that the abscissas of convergence,
absolute convergence and boundedness
of the Dirichlet series $\sum_{k=0}^{\infty} e^{-\lambda k/2^n}
\mathbb{E}\!\left[e^{i\mu X_{k,n}-\nu T_{k,n}}\ind_{[0,+\infty)}(X_{k,n})\right]$ are identical.
So, this series converges absolutely on the half-plane $\{\lambda \in \mathbb{C}: \Re(\lambda)>0\}$
and~(\ref{series}) holds on this half-plane. A similar conclusion holds for the second series
of Proposition~\ref{proposition}. The proof is finished.
\EndDim

Thanks to Proposition~\ref{proposition}, we see that the functional
$F_n(\lambda,\mu,\nu)$ is a function of the discrete observations of $X$
and, by Definition~\ref{def2}, its expectation can be computed as follows:
\begin{align}
E_n(\lambda,\mu,\nu)
&=
\frac{1-e^{-(\lambda+\nu)/2^n}}{\lambda+\nu}\,\frac{e^{\nu/2^n}-S_n^+(\lambda,\mu,\nu)}{e^{\nu/2^n}-1}
+\frac{1-e^{-\lambda/2^n}}{\lambda}\,\frac{S_n^-(\lambda,\mu,\nu)-1}{e^{\nu/2^n}-1}
\nonumber\\
&=
\left(\frac{e^{\nu/2^n}(1-e^{-(\lambda+\nu)/2^n})}{(\lambda+\nu)(e^{\nu/2^n}-1)}
-\frac{1-e^{-\lambda/2^n}}{\lambda(e^{\nu/2^n}-1)}\right)
\nonumber\\
&\hphantom{=\,}
+\frac{1-e^{-\lambda/2^n}}{\lambda (e^{\nu/2^n}-1)}\,S_n^-(\lambda,\mu,\nu)
-\frac{1-e^{-(\lambda+\nu)/2^n}}{(\lambda+\nu) (e^{\nu/2^n}-1)}\,S_n^+(\lambda,\mu,\nu).
\label{En}
\end{align}
Now, we have to evaluate the limit $E(\lambda,\mu,\nu)$ of $E_n(\lambda,\mu,\nu)$ as
$n$ goes toward infinity. It is easy to see that this limit exists; see the
proof of Theorem~\ref{theorem} below. Formally, we write
$E(\lambda,\mu,\nu)=\mathbb{E}[F(\lambda,\mu,\nu)]$
with
$$
F(\lambda,\mu,\nu)=\int_0^{\infty} e^{-\lambda t+i \mu X(t)-\nu T(t)} \,\mathrm{d}t.
$$
Then, we can say that the functional $F(\lambda,\mu,\nu)$ is an admissible
function of $X$ in the sense of Definition~\ref{def3}.
The value of its expectation $E(\lambda,\mu,\nu)$ is given in the following theorem.
%
\begin{teo}\label{theorem}
The $3$-parameters Laplace-Fourier transform of the couple $(T(t),X(t))$ is given by
\begin{equation}
E(\lambda,\mu,\nu)=\frac{1}{\prod_{j \in J}(\!\sqrt[N]{\lambda+\nu}-i\mu \theta_j )
\prod_{k \in K} (\!\sqrt[N]{\lambda}-i \mu \theta_k )}.
\label{expressionE}
\end{equation}
\end{teo}
%
\Dim
It is plain that the term lying within the biggest parentheses in the last equality
of~(\ref{En}) tends to zero as $n$ goes towards infinity and that the coefficients
lying before $S_n^+(\lambda,\mu,\nu)$ and $S_n^-(\lambda,\mu,\nu)$ tend to $1/\nu$.
As a byproduct, we derive at the limit when $n\to \infty$,
\begin{equation}
E(\lambda,\mu,\nu)=\frac{1}{\nu}\left[S^-(\lambda,\mu,\nu)-S^+(\lambda,\mu,\nu)\right]
\label{expressionEinter}
\end{equation}
where we set
\begin{align*}
S^+(\lambda,\mu,\nu)
&=
\lim_{n\to\infty}S_n^+(\lambda,\mu,\nu)
=\exp\!\left(-\int_0^{\infty} \mathbb{E}\!\left[\vphantom{e^X}\right.\!\! e^{i \mu X(t)}
\ind_{[0,+\infty)}(X(t))\!\!\left.\vphantom{e^X}\right]
\!\!\left(1-e^{-\nu t}\right)\frac{e^{-\lambda t}}{t}\,\mathrm{d}t\right)\!,
\\
S^-(\lambda,\mu,\nu)
&=
\lim_{n\to\infty}S_n^-(\lambda,\mu,\nu)
=\exp\!\left(\,\int_0^{\infty} \mathbb{E}\!\left[\vphantom{e^X}\right.\!\! e^{i \mu X(t)}
\ind_{(-\infty,0)}(X(t))\!\!\left.\vphantom{e^X}\right]\!\!
\left(1-e^{-\nu t}\right)\frac{e^{-\lambda t}}{t}\,\mathrm{d}t\right)\!.
\end{align*}
We have
\begin{align*}
\lqn{\int_0^{\infty} \mathbb{E}\!\left[\vphantom{e^X}\right.\!\! e^{i \mu X(t)}
\ind_{[0,+\infty)}(X(t))\!\!\left.\vphantom{e^X}\right]
\!\!\left(1-e^{-\nu t}\right)\frac{e^{-\lambda t}}{t}\,\mathrm{d}t}
&=
\int_0^{\infty} \mathbb{E}\!\left[\left(\!\!\vphantom{e^X}\right.\right.
e^{i \mu X(t)}-1\left.\left.\vphantom{e^X}\!\!\right)\!
\ind_{[0,+\infty)}(X(t))\right]\frac{e^{-\lambda t}}{t}\,\mathrm{d}t
-\int_0^{\infty} \mathbb{E}\!\left[\left(\!\!\vphantom{e^X}\right.\right.
e^{i \mu X(t)}-1\left.\left.\vphantom{e^X}\!\!\right)\!
\ind_{[0,+\infty)}(X(t))\right]\frac{e^{-(\lambda+\nu) t}}{t}\,\mathrm{d}t
\\
&\hphantom{=\,}
+\int_0^{\infty} \mathbb{P}\{X(t)\ge 0\} \,\frac{e^{-\lambda t}-e^{-(\lambda+\nu) t}}{t}\,\mathrm{d}t
\\
&=
\int_0^{\infty} \frac{e^{-\lambda t}}{t} \,\mathrm{d}t
\int_0^{\infty}\left(e^{i \mu \xi}-1\right) p(t;-\xi)\,\mathrm{d}\xi
-\int_0^{\infty} \frac{e^{-(\lambda+\nu)t}}{t} \,\mathrm{d}t
\int_0^{\infty}\left(e^{i \mu \xi}-1\right) p(t;-\xi)\,\mathrm{d}\xi
\\
&\hphantom{=\,}
+\mathbb{P}\{X(1)\ge 0\} \int_0^{\infty} \frac{e^{-\lambda t}-e^{-(\lambda+\nu) t}}{t}\,\mathrm{d}t.
\end{align*}
In view of~(\ref{set3}) and~(\ref{set1}) and using the elementary equality
$\int_0^{\infty} \frac{e^{-\lambda t}-e^{-(\lambda+\nu) t}}{t} \,\mathrm{d}t
=\log \left(\frac{\lambda+\nu}{\lambda}\right)$, we have
\begin{align*}
\lqn{\int_0^{\infty} \mathbb{E}\!\left[\vphantom{e^X}\right.\!\! e^{i \mu X(t)}
\ind_{[0,+\infty)}(X(t))\!\!\left.\vphantom{e^X}\right]\!\!
\left(1-e^{-\nu t}\right)\frac{e^{-\lambda t}}{t}\,\mathrm{d}t}
&=
\log\!\left(\!\vphantom{\prod_{\in}}\right.
\prod_{j \in J} \frac{\sqrt[N]{\lambda}}{\!\sqrt[N]{\lambda}-i\mu \theta_j}\left.\vphantom{\prod_{\in}}\!\right)\!
-\log\!\left(\!\vphantom{\prod_{\in}}\right.
\prod_{j \in J} \frac{\!\sqrt[N]{\lambda+\nu}}{\!\sqrt[N]{\lambda+\nu}-i\mu \theta_j}
\left.\vphantom{\prod_{\in}}\!\right)\!
+\frac{\#J}{N} \log\!\left(\frac{\lambda+\nu}{\lambda}\right)\!
=\log\!\left(\vphantom{\prod_{k}}\right.\!
\prod_{j \in J} \frac{\!\sqrt[N]{\lambda+\nu}-i\mu \theta_j}{\!\sqrt[N]{\lambda}-i\mu \theta_j}
\!\left.\vphantom{\prod_{k}}\right)\!.
\end{align*}
We then deduce the value of $S^+(\lambda,\mu,\nu)$. By~(\ref{set6}),
\begin{align}
S^+(\lambda,\mu,\nu)
&=
\prod_{j \in J} \frac{\!\sqrt[N]{\lambda}-i\mu \theta_j}{\!\sqrt[N]{\lambda+\nu}-i\mu \theta_j}
=\frac{\prod_{\ell=1}^N (\!\sqrt[N]{\lambda}-i\mu \theta_{\ell})}
{\prod_{j \in J} (\!\sqrt[N]{\lambda+\nu}-i\mu \theta_j)
\prod_{k \in K} (\!\sqrt[N]{\lambda}-i\mu \theta_k)}
\nonumber\\
&=
\frac{\lambda-\kappa_{\!_{ N}}(i\mu)^N}{\prod_{j \in J} (\!\sqrt[N]{\lambda+\nu}-i\mu \theta_j)
\prod_{k \in K} (\!\sqrt[N]{\lambda}-i\mu \theta_k)}.
\label{S+}
\end{align}
Similarly, the value of $S^-(\lambda,\mu,\nu)$ is given by
\begin{equation}
S^-(\lambda,\mu,\nu)=\prod_{k \in K} \frac{\!\sqrt[N]{\lambda+\nu}-i\mu \theta_k}{\!\sqrt[N]{\lambda}-i\mu \theta_k}
=\frac{\lambda+\nu-\kappa_{\!_{ N}}(i\mu)^N}{\prod_{j \in J} (\!\sqrt[N]{\lambda+\nu}-i\mu \theta_j)
\prod_{k \in K} (\!\sqrt[N]{\lambda}-i\mu \theta_k)}.
\label{S-}
\end{equation}
Finally, putting~(\ref{S+}) and~(\ref{S-}) into~(\ref{expressionEinter})
immediately leads to~(\ref{expressionE}).
\EndDim
%
\begin{rem}\label{remark}
In the particular case $\mu=0$, we get the very simple result:
$$
E(\lambda, 0,\nu)=\int_0^{\infty} e^{-\lambda t} \,\mathbb{E}
\!\left[\!\vphantom{e^t}\right. e^{-\nu T(t)}\!\left.\vphantom{e^t}\right] \mathrm{d}t
=\frac{1}{\lambda^{\frac{\# K}{N}} (\lambda+\nu)^{\frac{\# J}{N}}} .
$$
\pagebreak\noindent
This is formula (20) of \cite{2003}. On the other hand,
we can rewrite~(\ref{expressionE}) as
\begin{equation}
E(\lambda,\mu,\nu)=\frac{1}{\lambda^{\frac{\# K}{N}} (\lambda+\nu)^{\frac{\# J}{N}}}
\prod_{j \in J} \frac{\!\sqrt[N]{\lambda+\nu}}{\!\sqrt[N]{\lambda+\nu}-i\mu \theta_j}
\prod_{k \in K} \frac{\!\sqrt[N]{\lambda}}{\!\sqrt[N]{\lambda}-i \mu \theta_k}.
\label{expressionEbis}
\end{equation}
Actually, this form is more suitable for the inversion of the Laplace-Fourier transform.
\end{rem}
%

In the three next sections, we progressively invert the $3$-parameters
Laplace-Fourier transform $E(\lambda,\mu,\nu)$.

\section{Inverting with respect to $\mu$}\label{section-inverting-mu}

In this part, we invert $E(\lambda,\mu,\nu)$ given by~(\ref{expressionEbis})
with respect to $\mu$.
%
\begin{teo}\label{theorem-wrt.mu}
We have, for $\lambda,\nu>0$,
\begin{align}
\lqn{\int_0^{\infty} e^{-\lambda t} \left[\mathbb{E}\!\left(\!\!\vphantom{e^t}\right.\right.
e^{-\nu T(t)},\, X(t) \in \mathrm{d}x \!\!\left.\left.\vphantom{e^t}\right)\!/\mathrm{d}x\right]\mathrm{d}t}
&=
\begin{cases}
\displaystyle{\frac{1}{\lambda^{\frac{\# K-1}{N}}(\lambda+\nu)^{\frac{\# J-1}{N}}}
\sum_{j \in J}  A_j\theta_j \left(\,\sum_{k \in K}
\frac{B_k\theta_k}{\theta_k\!\!\sqrt[N]{\lambda}-\theta_j\!\!\sqrt[N]{\lambda+\nu}}\right)
e^{-\theta_j\!\!\sqrt[N]{\lambda+\nu}\,x}} & \mbox{if } x \ge 0,
\\[3ex]
\displaystyle{\frac{1}{\lambda^{\frac{\# K-1}{N}}(\lambda+\nu)^{\frac{\# J-1}{N}}}
\sum_{k \in K}B_k\theta_k \left(\vphantom{\prod_{k}}\right.\!\sum_{j \in J}
\frac{A_j\theta_j}{\theta_k\!\!\sqrt[N]{\lambda}-\theta_j\!\!\sqrt[N]{\lambda+\nu}}
\!\left.\vphantom{\prod_{k}}\right)
e^{-\theta_k\!\!\sqrt[N]{\lambda}\,x}} &  \mbox{if } x \le 0.
\end{cases}
\label{wrt.mu}
\end{align}
\end{teo}
%
\Dim
By~(\ref{set5}) applied to $x=i\mu/\!\sqrt[N]{\lambda+\nu}$ and
$x=i\mu/\!\sqrt[N]{\lambda}$, we have

\begin{align*}
\prod_{j \in J}  \frac{\!\sqrt[N]{\lambda+\nu}}{ \sqrt[N]{\lambda+\nu}-i\mu \theta_j}
\prod_{k \in K} \frac{\!\sqrt[N]{\lambda}}{ \sqrt[N]{\lambda}  -i \mu \theta_k}
&=
\prod_{j \in J}  \frac{1}{1-\frac{i\mu}{\!\sqrt[N]{\lambda+\nu}}\,\theta_j}
\prod_{k \in K}  \frac{1}{1-\frac{i\mu}{\!\sqrt[N]{\lambda}}\,\theta_k}
\\
&=
\sum_{j \in J}  \frac{A_j\theta_j}{\theta_j-\frac{i\mu}{\!\sqrt[N]{\lambda+\nu}}}
\sum_{k \in K}  \frac{B_k\theta_k}{\theta_k-\frac{i\mu}{\!\sqrt[N]{\lambda}}}
\\
&=
\sqrt[N]{\lambda (\lambda+\nu)} \sum_{\substack{j \in J \\k \in K}}
\frac{A_j B_k\theta_j \theta_k}{(\theta_j\!\!\sqrt[N]{\lambda+\nu}-i\mu ) (\theta_k\!\!\sqrt[N]{\lambda}-i\mu )}.
\end{align*}
Let us write that
\begin{align*}
\frac{1}{(\theta_j\!\!\sqrt[N]{\lambda+\nu}-i\mu ) (\theta_k\!\!\sqrt[N]{\lambda}-i\mu )}
&=
\frac{1}{\theta_k\!\!\sqrt[N]{\lambda}-\theta_j\!\!\sqrt[N]{\lambda+\mu}}
\left(\frac{1}{\theta_j\!\!\sqrt[N]{\lambda+\nu}-i \mu}
-\frac{1}{\theta_k\!\!\sqrt[N]{\lambda}  -i \mu}\right)
\\
&=
\frac{1}{\theta_k\!\!\sqrt[N]{\lambda}-\theta_j\!\!\sqrt[N]{\lambda+\mu}}
\left(\int_0^{\infty} e^{(i \mu -\theta_j\!\!\sqrt[N]{\lambda+\mu}) x} \,\mathrm{d}x
+\int_{-\infty}^0 e^{(i \mu -\theta_k\!\!\sqrt[N]{\lambda}) x} \,\mathrm{d}x\right)\!.
\end{align*}
Therefore, we can rewrite $E(\lambda,\mu,\nu)$ as
\begin{align*}
E(\lambda,\mu,\nu)
&=
\frac{1}{\lambda^{\frac{\# K-1}{N}}(\lambda+\nu)^{\frac{\# J-1}{N}}}
\\
&\hphantom{=\,}
\times\sum_{\substack{j \in J \\k \in K}}
\frac{A_j B_k\theta_j \theta_k}{\theta_k\!\!\sqrt[N]{\lambda}-\theta_j\!\!\sqrt[N]{\lambda+\nu}}
\int_{-\infty}^{\infty} e^{i \mu x} \left(e^{-\theta_k\!\!\sqrt[N]{\lambda}\,x}
\ind_{(-\infty,0]}(x)+ e^{-\theta_j\!\!\sqrt[N]{\lambda+\nu}\,x} \ind_{[0,\infty)}(x)\right) \mathrm{d}x
\end{align*}
which is nothing but the Fourier transform with respect to $\mu$ of the
right-hand side of~(\ref{wrt.mu}).
\EndDim

%
\begin{rem}
$\bullet$
By integrating~(\ref{wrt.mu}) on $(-\infty,0]$, we obtain
$$
\int_0^{\infty} e^{-\lambda t} \,\mathbb{E}\!\left(\!\vphantom{e^t}\right.
e^{-\nu T(t)},\, X(t) \le 0 \left.\!\vphantom{e^t}\right) \mathrm{d}t
=-\frac{1}{\lambda^{\frac{\# K}{N}}(\lambda+\nu)^{\frac{\# J-1}{N}}}
\sum_{k \in K}  B_k \left(\!\vphantom{\prod_{k}}\right.
\sum_{j \in J}\frac{A_j\theta_j}{\theta_k\!\!\sqrt[N]{\lambda}-\theta_j\!\!\sqrt[N]{\lambda+\nu}}
\!\left.\vphantom{\prod_{k}}\right)\!.
$$
Using~(\ref{set5}) applied to $x=\theta_k\!\!\sqrt[N] \lambda/ \!\sqrt[N]{\lambda+\nu}$
and~(\ref{set10}), we see that
\pagebreak
\begin{align}
\sum_{j \in J} \frac{A_j\theta_j}{\theta_j-\theta_k\!\!\sqrt[N]{\frac{\lambda}{\lambda+\nu}}}
&=
\frac{\lambda+\nu}{\nu} \prod_{i\in K} \left(1-\bar{\theta}_i\theta_k
\!\!\sqrt[N]{\frac{\lambda}{\lambda+\nu}}\,\right)
\nonumber\\
&=
\theta_k^{\#K} \frac{\lambda+\nu}{\nu} \left(\frac{\lambda}{\lambda+\nu}\right)^{\!\frac{\#K}{N}}
\prod_{i\in K} \left(\bar{\theta}_k\!\!\sqrt[N]{\frac{\lambda+\nu}{\lambda}}-\bar{\theta}_i\right)
\nonumber\\
&=
\frac{1}{\nu}\,\theta_k^{\#K} \lambda^{\frac{\#K}{N}}(\lambda+\nu)^{\frac{\#J}{N}}\,
\sum_{\ell=0}^{\#K} (-1)^{\ell} \bar{\sigma}_{\ell}\,
\bar{\theta}_k^{\#K-\ell} \left(\frac{\lambda+\nu}{\lambda}\right)^{\!\frac{\#K-\ell}{N}}
\nonumber\\
&=
\frac{1}{\nu}\, \lambda^{\frac{\#K}{N}}(\lambda+\nu)^{\frac{\#J}{N}}\,
\sum_{\ell=0}^{\#K} (-1)^{\ell} \bar{\sigma}_{\ell}\,
\theta_k^{\ell} \left(\frac{\lambda+\nu}{\lambda}\right)^{\!\frac{\#K-\ell}{N}}.
\label{sum-inter}
\end{align}
This entails that
\begin{align*}
\int_0^{\infty} e^{-\lambda t} \,\mathbb{E}\!\left(\!\vphantom{e^t}\right.
e^{-\nu T(t)},\, X(t) \le 0 \left.\!\vphantom{e^t}\right) \mathrm{d}t
&=
\frac{1}{\nu} \,\sum_{k \in K} B_k \sum_{\ell=0}^{\# K} (-1)^{\ell}
\bar{\sigma}_{\ell} \,\theta_k^{\ell}
\left(\frac{\lambda+\nu}{\lambda}\right)^{\!\frac{\#K-\ell}{N}}
\\
&=
\frac{1}{\nu} \,\sum_{\ell=0}^{\# K} (-1)^{\ell} \bar{\sigma}_{\ell}\,
\beta_{\ell} \left(\frac{\lambda+\nu}{\lambda}\right)^{\!\frac{\#K-\ell}{N}}.
\end{align*}
By~(\ref{set26}), we know that all the $\beta_\ell$, $1\le \ell\le \#K-1$ vanish
and it remains, with~(\ref{set16}),
\begin{align}
\int_0^{\infty} e^{-\lambda t} \,\mathbb{E}\!\left(\!\vphantom{e^t}\right.
e^{-\nu T(t)},\, X(t) \le 0 \left.\!\vphantom{e^t}\right) \mathrm{d}t
&=
\frac{1}{\nu} \left[\bar{\sigma}_0\beta_0\left(\frac{\lambda+\nu}{\lambda}
\right)^{\!\frac{\# K}{N}}+(-1)^{\#K}\bar{\sigma}_{\#K}\beta_{\# K}\right]
\nonumber\\
&=
\frac{1}{\nu} \left[\left(\frac{\lambda+\nu}{\lambda}\right)^{\!\frac{\# K}{N}} -1\right]\!.
\label{le2}
\end{align}
We retrieve (30) of \cite{2003}.

$\bullet$
Likewise, we have
\begin{equation}
\int_0^{\infty} e^{-\lambda t} \,\mathbb{E}\!\left(\!\vphantom{e^t}\right.
e^{-\nu T(t)},\, X(t) \le 0 \left.\!\vphantom{e^t}\right) \mathrm{d}t
=\frac{1}{\nu} \left[1-\left(\frac{\lambda}{\lambda+\nu}\right)^{\!\frac{\# J}{N}}\right]\!,
\label{ge2}
\end{equation}
which coincides with (29) of \cite{2003}.

$\bullet$
Adding formulas~(\ref{le2}) and~(\ref{ge2}) we obtain
\begin{align*}
\int_0^{\infty} e^{-\lambda t} \,\mathbb{E}\!\left(\!\vphantom{e^t}\right.
e^{-\nu T(t)} \left.\!\vphantom{e^t}\right) \mathrm{d}t
&=
\frac{1}{\nu} \left[\left(\frac{\lambda+\nu}{\lambda}\right)^{\!\frac{\# K}{N}}
-\left(\frac{\lambda}{\lambda+\nu}\right)^{\!\frac{\# J}{N}}\right]
=\frac{(\lambda+\nu)^{\frac{\#J+\#K}{N}}-\lambda^{\frac{\# J+\#K}{N}}}{\nu
\lambda^{\frac{\# K}{N}}(\lambda+\nu)^{\frac{\# J}{N}}}
\\
&=
\frac{1}{\lambda^{\frac{\#K}{N}}(\lambda+\nu)^{\frac{\# J}{N}}}.
\end{align*}
This is formula (10) of \cite{2003} which has already been pointed out in Remark~\ref{remark}.
Another way of checking this formula consists of integrating~(\ref{wrt.mu})
with respect to $x$ directly on $\mathbb{R}$. Indeed,
\begin{align*}
\lqn{\int_0^{\infty} e^{-\lambda t} \,\mathbb{E}\!\left(\!\vphantom{e^t}\right.
e^{-\nu T(t)} \left.\!\vphantom{e^t}\right) \mathrm{d}t}
\\[-5ex]
&=
\frac{1}{\lambda^{\frac{\#K-1}{N}}(\lambda+\nu)^{\frac{\# J-1}{N}}}
\left(\frac{1}{\!\sqrt[N]{\lambda+\nu}} \sum_{{\substack{j\in J \\ k \in K}}}
\frac{A_j B_k\theta_k}{\theta_k\!\!\sqrt[N]{\lambda}-\theta_j\!\!\sqrt[N]{\lambda+\nu}}
-\frac{1}{\!\sqrt[N]{\lambda}} \sum_{{\substack{j\in J \\ k \in K}}}
\frac{A_j B_k\theta_j}{\theta_k\!\!\sqrt[N]{\lambda}-\theta_j\!\!\sqrt[N]{\lambda+\nu}}\right)
\\
&=
\frac{1}{\lambda^{\frac{\#K-1}{N}}(\lambda+\nu)^{\frac{\# J-1}{N}}}
\sum_{\substack{j\in J \\ k \in K}} \left(\frac{A_j B_k\theta_k}{\!\sqrt[N]{\lambda+\nu}\,
(\theta_k\!\!\sqrt[N]{\lambda}-\theta_j\!\!\sqrt[N]{\lambda+\nu})}
-\frac{A_j B_k\theta_j}{\!\sqrt[N]{\lambda}\,(\theta_k\!\!\sqrt[N]{\lambda}
-\theta_j\!\!\sqrt[N]{\lambda+\nu})}\right)
\\
&=
\frac{1}{\lambda^{\frac{\#K-1}{N}}(\lambda+\nu)^{\frac{\# J-1}{N}}}
\sum_{\substack{j\in J \\ k \in K}} \frac{A_j B_k}{\!\sqrt[N]{\lambda(\lambda+\nu)}}.
\end{align*}
By~(\ref{set13}), we have $\sum_{\substack{j\in J \\ k \in K}} A_j B_k
=\sum_{j\in J} A_j \sum_{k\in K} B_k=1$ and then
$$
\int_0^{\infty} e^{-\lambda t} \,\mathbb{E}\!\left(\!\vphantom{e^t}\right.
e^{-\nu T(t)} \left.\!\vphantom{e^t}\right) \mathrm{d}t
=\frac{1}{\lambda^{\frac{\#K}{N}}(\lambda+\nu)^{\frac{\# J}{N}}}.
$$
\end{rem}
%
\begin{rem}\label{remark-x0}
By replacing $x$ by $0$ into~(\ref{wrt.mu}) and by using~(\ref{sum-inter}), we get
\begin{align*}
\int_0^{\infty} e^{-\lambda t} \,\mathbb{E}\!\left(\!\vphantom{e^t}\right.
e^{-\nu T(t)},\, X(t)\in\mathrm{d}x \left.\!\vphantom{e^t}\right)\!/\mathrm{d}x\,\Big|_{x=0} \,\mathrm{d}t
&=
\frac{1}{\lambda^{\frac{\# K-1}{N}}(\lambda+\nu)^{\frac{\# J-1}{N}}}
\sum_{k \in K} B_k\theta_k \left(\!\vphantom{\prod_{k}}\right.\sum_{j \in J}
\frac{A_j\theta_j}{\theta_k\!\!\sqrt[N]{\lambda}-\theta_j\!\!\sqrt[N]{\lambda+\nu}}
\left.\!\vphantom{\prod_{k}}\right)\!.
\\
&=
-\frac{\!\sqrt[N]{\lambda}}{\nu} \,\sum_{k \in K} B_k\theta_k
\sum_{\ell=0}^{\# K} (-1)^{\ell} \bar{\sigma}_{\ell} \,\theta_k^{\ell}
\left(\frac{\lambda+\nu}{\lambda}\right)^{\!\frac{\#K-\ell}{N}}
\\
&=
-\frac{\!\sqrt[N]{\lambda}}{\nu}\,\sum_{\ell=0}^{\# K} (-1)^{\ell}
\bar{\sigma}_{\ell} \,\beta_{\ell+1}
\left(\frac{\lambda+\nu}{\lambda}\right)^{\!\frac{\#K-\ell}{N}}.
\end{align*}
In this sum, all the $\beta_{\ell+1}$, $0\le\ell\le \#K-2$, vanish
and it remains, with~(\ref{set16}),
\begin{align*}
\lqn{\int_0^{\infty} e^{-\lambda t} \,\mathbb{E}\!\left(\!\vphantom{e^t}\right.
e^{-\nu T(t)},\, X(t)\in\mathrm{d}x \left.\!\vphantom{e^t}\right)\!
/\mathrm{d}x\,\Big|_{x=0} \,\mathrm{d}t}
&=
-\frac{\!\sqrt[N]{\lambda}}{\nu}\left(
(-1)^{\#K-1}\bar{\sigma}_{\#K-1}\beta_{\# K}\!\sqrt[N]{\frac{\lambda+\nu}{\lambda}}
+(-1)^{\#K}\bar{\sigma}_{\#K}\beta_{\# K+1}\right)
\\
&=
-\frac{\!\sqrt[N]{\lambda}}{\nu} \,\sum_{k \in K} \theta_k\!
\left(\!\sqrt[N]{\frac{\lambda+\nu}{\lambda}}-1\right)
=\left(\!\vphantom{\prod_{k}}\right.\sum_{j \in J}\theta_j
\left.\!\!\vphantom{\prod_{k}}\right) \frac{\!\sqrt[N]{\lambda+\nu}-\sqrt[N]\lambda}{\nu}.
\end{align*}
We retrieve formula (26) of~\cite{2003}.
\end{rem}
%
\begin{rem}
For $\nu=0$, formula~(\ref{wrt.mu}) yields with~(\ref{set11}) the $\lambda$-potential of $X$:
\begin{align*}
\int_0^{\infty} e^{-\lambda t} \left[\mathbb{P}\{X(t) \in \mathrm{d}x\}/\mathrm{d}x\right] \mathrm{d}t
&=
\begin{cases}
\displaystyle{\frac{1}{\lambda^{1-\frac{1}{N}}} \sum_{j \in J}  A_j\theta_j \sum_{k \in K}
\frac{B_k\theta_k}{\theta_k-\theta_j} \,e^{-\theta_j\!\!\sqrt[N]{\lambda} \,x}}& \mbox{if } x \ge 0,
\\
\displaystyle{\frac{1}{\lambda^{1-\frac{1}{N}}} \sum_{k \in K} B_k\theta_k \sum_{j \in J}
\frac{A_j\theta_j}{\theta_k-\theta_j} \,e^{-\theta_k\!\!\sqrt[N]{\lambda} \,x}}& \mbox{if } x \le 0,
\end{cases}
\\[2ex]
&=
\begin{cases}
\displaystyle{\frac{1}{N \lambda^{1-\frac{1}{N}}} \sum_{j \in J} \theta_j\,
e^{-\theta_j\!\!\sqrt[N]{\lambda} \,x}}& \mbox{if } x \ge 0,
\\
\displaystyle{-\frac{1}{N \lambda^{1-\frac{1}{N}}}\sum_{k \in K} \theta_k
e^{-\theta_k\!\!\sqrt[N]{\lambda} \,x}} & \mbox{if } x \le 0.
\end{cases}
\end{align*}
We retrieve (12) of \cite{2003}.
\end{rem}
%
\begin{ex}
For $N=2$, formula~(\ref{wrt.mu}) gives, with the numerical values of Example~\ref{example1},
$$
\int_0^{\infty} e^{-\lambda t} \left[\mathbb{E}\!\left(\!\!\vphantom{e^t}\right.\right.
e^{-\nu T(t)},\, X(t) \in \mathrm{d}x \!\!\left.\left.\vphantom{e^t}\right)\!/\mathrm{d}x\right]\mathrm{d}t
=\begin{cases}
\displaystyle{\frac{1}{\sqrt{\lambda}+\sqrt{\lambda+\nu}}
\,e^{-\sqrt{\lambda+\nu}\,x}} & \mbox{if } x \ge 0,
\\[3ex]
\displaystyle{ \frac{1}{\sqrt{\lambda}+\sqrt{\lambda+\nu}}
\,e^{\sqrt{\lambda}\,x}}&  \mbox{if } x \le 0.
\end{cases}
$$
This is formula 1.4.5, p. 129, of~\cite{bs}.
\end{ex}
%
\begin{ex}
For $N=3$, we have two cases to consider.
Formula~(\ref{wrt.mu}) yields, with the numerical values of Example~\ref{example2},
in the case $\kappa_3=1$,
\pagebreak
\begin{align*}
\lqn{\int_0^{\infty} e^{-\lambda t} \left[\mathbb{E}\!\left(\!\!\vphantom{e^t}\right.\right.
e^{-\nu T(t)},\, X(t) \in \mathrm{d}x \!\!\left.\left.\vphantom{e^t}\right)\!/\mathrm{d}x\right]\mathrm{d}t}
=\begin{cases}
\displaystyle\frac{e^{-\sqrt[3]{\lambda +\nu}\,x}}{\lambda^{2/3}
+\sqrt[3]{\lambda(\lambda+\nu)}+ (\lambda+\nu)^{2/3}} & \mbox{if } x\ge 0,
\\[3ex]
\displaystyle
\frac{e^{\frac{\sqrt[3]{\lambda}}{2}\,x}}{\sqrt3\,\sqrt[3] \lambda}
\frac{\sqrt3 \,\sqrt[3] \lambda \, \cos \!\left(\!\vphantom{\frac12}\right.
\frac{\sqrt3\sqrt[3]{\lambda}}{2}\,x \left.\vphantom{\frac12}\!\right)
-(2\sqrt[3]{\lambda+\nu}+\sqrt[3]\lambda\,) \sin\!\left(\!\vphantom{\frac12}\right.
\frac{\sqrt3\sqrt[3]{\lambda}}{2}\,x \left.\vphantom{\frac12}\!\right)}{
\lambda^{2/3} +\sqrt[3]{\lambda(\lambda+\nu)}+ (\lambda+\nu)^{2/3}} & \mbox{if } x\le 0,
\end{cases}
\end{align*}
and in the case $\kappa_3=-1$,
\begin{align*}
\lqn{\int_0^{\infty} e^{-\lambda t} \left[\mathbb{E}\!\left(\!\!\vphantom{e^t}\right.\right.
e^{-\nu T(t)},\, X(t) \in \mathrm{d}x \!\!\left.\left.\vphantom{e^t}\right)\!/\mathrm{d}x\right]\mathrm{d}t}
=\begin{cases}
\displaystyle\frac{e^{-\frac{\sqrt[3]{\lambda+\nu}}{2}\,x}}{ \sqrt 3 \sqrt[3]{ \lambda+\nu}}\,
\frac{\sqrt 3 \,\sqrt[3]{\lambda+\nu} \,\cos\!\left(\!\vphantom{\frac12}\right.
\frac{\sqrt3\sqrt[3]{\lambda+\nu}}{2}\,x \left.\vphantom{\frac12}\!\right)
+(\sqrt[3]{\lambda+\nu}+2 \sqrt[3] \lambda\,) \sin\!\left(\!\vphantom{\frac12}\right.
\frac{\sqrt3\sqrt[3]{\lambda+\nu}}{2}\,x \left.\vphantom{\frac12}\!\right)}{\lambda^{2/3}
+\sqrt[3]{\lambda(\lambda+\nu)}+ (\lambda+\nu)^{2/3}} & \mbox{if } x\ge 0,
\\[2ex]
\displaystyle
\frac{e^{\sqrt[3]{\lambda}\,x}}{\lambda^{2/3} +\sqrt[3]{\lambda(\lambda+\nu)}+ (\lambda+\nu)^{2/3}}
& \mbox{if } x\le 0.
\end{cases}
\end{align*}
\end{ex}
%
\begin{ex}
For $N=4$, formula~(\ref{wrt.mu}) supplies, with the numerical values of Example~\ref{example3},
\begin{align*}
\lqn{\int_0^{\infty} e^{-\lambda t} \left[\mathbb{E}\!\left(\!\!\vphantom{e^t}\right.\right.
e^{-\nu T(t)},\, X(t) \in \mathrm{d}x \!\!\left.\left.\vphantom{e^t}\right)\!/\mathrm{d}x\right]\mathrm{d}t}
=\begin{cases}
\displaystyle\frac{\sqrt2\,e^{-\frac{\sqrt[4]{\lambda+\nu}}{\sqrt{2}}\,x}}{\sqrt[4]{\lambda+\nu}\,(\sqrt \lambda
+\sqrt{\lambda+\nu})(\sqrt[4] \lambda +\sqrt[4]{\lambda+\nu})}
\left[\sqrt[4]{\lambda+\nu} \,\cos \!\left(\frac{\sqrt[4]{\lambda+\nu}}{\sqrt{2}}\,x\right)\!
+ \sqrt[4]{\lambda} \,\sin \!\left(\frac{\sqrt[4]{\lambda+\nu}}{\sqrt{2}}\,x\right)\right]& \mbox{if } x\ge 0,
\\[3ex]
\displaystyle
\frac{\sqrt2\,e^{\frac{\sqrt[4]{\lambda}}{\sqrt{2}}\,x}}{\sqrt[4]{\lambda} \,(\sqrt \lambda
+\sqrt{\lambda+\nu})(\sqrt[4] \lambda +\sqrt[4]{\lambda+\nu})}
\left[\vphantom{\frac{\sqrt[4]{\lambda+\nu}}{\sqrt{2}}}\!\!\right.
\sqrt[4]{\lambda} \,\cos \!\left(\vphantom{\frac{\sqrt[4]{\lambda+\nu}}{\sqrt{2}}}\right.\!\!
\frac{\sqrt[4]{\lambda}}{\sqrt{2}}\,x \left.\vphantom{\frac{\sqrt[4]{\lambda+\nu}}{\sqrt{2}}}\!\!\right)\!
-\sqrt[4]{\lambda+\nu} \,\sin \!\left(\vphantom{\frac{\sqrt[4]{\lambda+\nu}}{\sqrt{2}}}\right.\!\!
\frac{\sqrt[4]{\lambda}}{\sqrt{2}}\,x \left.\vphantom{\frac{\sqrt[4]{\lambda+\nu}}{\sqrt{2}}}\!\!\right)
\left.\vphantom{\frac{\sqrt[4]{\lambda+\nu}}{\sqrt{2}}}\!\!\right]
& \mbox{if } x\le 0.
\end{cases}
\end{align*}
\end{ex}
%
\begin{rem}
Using quite analogous computations to those of Remark~\ref{remark-x0}, we could
derive another expression for formula~(\ref{wrt.mu}).
Actually it will  not be used for the inversion.
\begin{align}
\lqn{\int_0^{\infty} e^{-\lambda t} \left[\mathbb{E}\!\left(\!\!\vphantom{e^t}\right.\right.
e^{-\nu T(t)},\, X(t) \in \mathrm{d}x \!\!\left.\left.\vphantom{e^t}\right)\!/\mathrm{d}x\right]\mathrm{d}t}
&=
\begin{cases}
\displaystyle{\frac{\!\sqrt[N]{\lambda+\nu}-\sqrt[N]{\lambda}}{ \nu} \sum_{j \in J} \theta_j
\left(\prod_{i \in J \setminus \{j\}}  \frac{\theta_i\!\!\sqrt[N]{\frac{\lambda}{\lambda+\nu}}
-\theta_j}{\theta_i-\theta_j}\right) e^{-\theta_j\!\!\sqrt[N]{\lambda+\nu}\,x}} & \mbox{if } x \ge 0,
\\[4ex]
\displaystyle{-\frac{\!\sqrt[N]{\lambda+\nu}-\sqrt[N]{\lambda}}{ \nu} \sum_{k \in K}
\theta_k  \left(\prod_{i \in K \setminus \{k\}}  \frac{\theta_i\!\!\sqrt[N]{\frac{\lambda+\nu}{\lambda}}
-\theta_k}{\theta_i-\theta_k}\right)e^{-\theta_k\!\!\sqrt[N]{\lambda}\,x}} &  \mbox{if } x \le 0.
\end{cases}
\label{wrt.mu2}
\end{align}
Formula~(\ref{wrt.mu2}) looks like formula~(24) in~\cite{2003}.
Nevertheless, (\ref{wrt.mu2}) involves the distribution of $(T(t), X(t))$
when the pseudo-process starts at zero while (24) of~\cite{2003} involves
the density of $(T(t), X(t))$ evaluated at the extremity $X(t)=0$ when the
starting point is $x$. Actually, both formulas are identical by invoking the
duality upon changing $x$ into $-x$, but they were obtained through different approaches.
\end{rem}
%

\section{Inverting with respect to $\nu$}\label{section-inverting-nu}

In this section, we carry out the inversion with respect to the parameter $\nu$.
The cases $x\le 0$ and $x\ge 0$ lead to results which are not quite analogous.
This is due to the asymmetry of our problem. So, we split our analysis into two
subsections related to the cases $x\le 0$ and $x\ge 0$.

\subsection{The case $x\le 0$}
%
\begin{teo}\label{theorem-wrt.nu+1}
The Laplace transform with respect to $t$ of the density of the couple
$(T(t),X(t))$ is given, when $x\le 0$, by
\pagebreak
\begin{align}
\lqn{\int_0^{\infty} e^{-\lambda t}\left[\mathbb{P}\{T(t)\in \mathrm{d}s,
X(t) \in \mathrm{d}x\} /(\mathrm{d}s \,\mathrm{d}x)\right]\mathrm{d}t}
&=
-\frac{e^{-\lambda s}}{\lambda^{\frac{\# K -1}{N}} s^{\frac{\# K}{N}}}
\sum_{m=0}^{\#K} \alpha_{-m}\,(\lambda s)^{\frac mN} E_{1,\frac{m+\# J}{N}}(\lambda s)
\sum_{k \in K} B_k\theta_k^{m+1} e^{-\theta_k\!\!\sqrt[N]{\lambda}\,x}.
\end{align}
\label{wrt.nu}
\end{teo}
%
\Dim
Recall~(\ref{wrt.mu}) in the case $x \le 0$:
\begin{align*}
\lqn{\int_0^{\infty} e^{-\lambda t} \left[\mathbb{E}\!\left(\!\!\vphantom{e^t}\right.\right.
e^{-\nu T(t)},\, X(t) \in \mathrm{d}x \!\!\left.\left.\vphantom{e^t}\right)\!/\mathrm{d}x\right]\mathrm{d}t}
&=
\frac{1}{\lambda^{\frac{\# K-1}{N}}(\lambda+\nu)^{\frac{\# J-1}{N}}}
\sum_{k \in K} B_k\theta_k \left(\!\vphantom{\prod_{k}}\right.
\sum_{j \in J} \frac{A_j\theta_j}{\theta_k\!\!\sqrt[N]{\lambda}
-\theta_j\!\!\sqrt[N]{\lambda+\nu}}\left.\!\vphantom{\prod_{k}}\right)e^{-\theta_k\!\!\sqrt[N]{\lambda}\,x}.
\end{align*}
We have to invert with respect to $\nu$ the quantity
$$
\frac{1}{(\lambda+\nu)^{\frac{\# J-1}{N}}} \sum_{j \in J} \frac{A_j\theta_j}{\theta_k\!\!\sqrt[N]{\lambda}-\theta_j\!\!\sqrt[N]{\lambda+\nu}}
=-\sum_{j \in J}  \frac{A_j}{(\lambda+\nu)^{\frac{\# J-1}{N}} (\!\sqrt[N]{\lambda+\nu}
-\frac{\theta_k}{\theta_j} \sqrt[N]{\lambda})}.
$$
By using the following elementary equality, which is valid for $\alpha>0$,
$$
\frac{1}{(\lambda+\nu)^{\alpha}}=\frac{1}{\Gamma(\alpha)} \int_0^{\infty}
e^{-(\lambda+\nu)s}s^{\alpha-1}\,\mathrm{d}s=\int_0^{\infty} e^{-\nu s}
\left(\frac{s^{\alpha-1} e^{-\lambda s}}{\Gamma(\alpha)}\right) \mathrm{d}s,
$$
we obtain, for $|\beta|<\sqrt[N]{\lambda+\nu}$,
\begin{align*}
\frac{1}{\!\sqrt[N]{\lambda+\nu}-\beta}
&=
\frac{1}{\!\sqrt[N]{\lambda+\nu}} \,\frac{1}{1-\frac{\beta}{\!\sqrt[N]{\lambda+\nu}}}
=\sum_{r=0}^{\infty} \frac{\beta^r}{(\lambda+\nu)^{\frac{r+1}{N}}}
= \sum_{r=0}^{\infty} \frac{\beta^r}{\Gamma\!\left(\vphantom{\frac aN}\right.\!\!
\frac{r+1}{N} \!\!\left.\vphantom{\frac aN}\right)}\int_0^{\infty}e^{-(\lambda+\nu)s}
s^{\frac{r+1}{N}-1} \,\mathrm{d}s
\\
&=
\int_0^{\infty} e^{-\nu s} \left(s^{\frac{1}{N}-1}e^{-\lambda s} \sum_{r=0}^{\infty}
\frac{(\beta \!\sqrt[N]{s}\,)^r}{\Gamma\!\left(\vphantom{\frac aN}\right.\!\!
\frac{r+1}{N}\!\!\left.\vphantom{\frac aN}\right)}\right) \mathrm{d}s.
\end{align*}
The sum lying in the last displayed equality can be expressed by means of
the Mittag-Leffler function (see \cite[Chap.~\textsc{xviii}]{erdelyi}):
$E_{a,b}(\xi)=\sum_{r=0}^{\infty} \frac{\xi^r}{\Gamma(ar+b)}$.
Then,
\begin{equation}
\frac{1}{\!\sqrt[N]{\lambda+\nu}-\beta}=\int_0^{\infty} e^{-\nu s}
\left(s^{\frac{1}{N}-1}e^{-\lambda s} E_{\frac{1}{N},\frac{1}{N}}(\beta \sqrt[N]{s}\,)\right) \mathrm{d}s.
\label{mittag}
\end{equation}
Next, we write
\begin{equation}
\sum_{j \in J}  \frac{A_j}{ \sqrt[N]{\lambda+\nu}-\frac{\theta_k}{\theta_j}
\sqrt[N]{\lambda}}=\int_0^{\infty} e^{-\nu s} \left[\vphantom{\sum_\in}\right.\!
s^{\frac{1}{N}-1}e^{-\lambda s}
\sum_{j \in J} A_j \,E_{\frac{1}{N},\frac{1}{N}} \!\left(\frac{\theta_k}{\theta_j}\!
\sqrt[N]{\lambda s}\right)\!\!\left.\vphantom{\sum_\in}\right] \mathrm{d}s,
\label{a_sum}
\end{equation}
where
$$
\sum_{j \in J} A_j \,E_{\frac{1}{N},\frac{1}{N}} \!\left(\frac{\theta_k}{\theta_j}\!
\sqrt[N]{\lambda s}\right)=\sum_{j \in J} A_j \sum_{r=0}^{\infty}
\left(\frac{\theta_k}{\theta_j}\right)^{\!r} \frac{(\lambda s)^{\frac{r}{N}}}{\Gamma
\!\left(\vphantom{\frac aN}\right.\!\!\frac{r+1}{N}\!\!\left.\vphantom{\frac aN}\right)}
=\sum_{r=0}^{\infty} \left(\!\vphantom{\prod_{k}}\right.
\theta_k^r \,\sum_{j \in J} \frac{A_j}{\theta_j^r}\!\left.\vphantom{\prod_{k}}\right)
\frac{(\lambda s)^{\frac{r}{N}}}{\Gamma\!\left(\vphantom{\frac aN}\right.\!\!
\frac{r+1}{N}\!\!\left.\vphantom{\frac aN}\right)}.
$$
When performing the euclidian division of $r$ by $N$, we can write $r$ as
$r=\ell N+m$ with $\ell\ge 0$ and $0 \le m \le N-1$. With this, we have
$\theta_j^{-r}=(\theta_j^N)^{-\ell} \,\theta_j^{-m}=\kappa_{\!_{ N}}^{\ell} \,\theta_j^{-m}$ and
$\theta_k^r=\kappa_{\!_{ N}}^{\ell} \,\theta_k^m$. Then,
$$
\theta_k^r \,\sum_{j \in J} \frac{A_j}{\theta_j^r}
=\theta_k^m \,\sum_{j \in J} \frac{A_j}{\theta_j^m}=\theta_k^m \alpha_{-m}.
$$
Hence, since by~(\ref{set14}) the $\alpha_{-m}$, $\#K+1\le m\le N$, vanish,
$$
\sum_{j \in J} A_j \,E_{\frac{1}{N},\frac{1}{N}}\!\left(\frac{\theta_k}{\theta_j}\!\sqrt[N]{\lambda s}\right)
=\sum_{\ell=0}^{\infty} \sum_{m=0}^{\# K} \alpha_{-m} \theta_k^m\,
\frac{(\lambda s)^{\ell+\frac{m}{N}}}{\Gamma\!\left(\vphantom{\frac aN}\right.\!\!
\ell+\frac{m+1}{N}\!\!\left.\vphantom{\frac aN}\right)}
=\sum_{m=0}^{\# K} \alpha_{-m} \,\theta_k^m (\lambda s)^{\frac{m}{N}} E_{1,\frac{m+1}{N}}(\lambda s)
$$
and~(\ref{a_sum}) becomes
$$
\sum_{j \in J} \frac{A_j}{\!\sqrt[N]{\lambda+\nu}-\frac{\theta_k}{\theta_j}\!\sqrt[N]{\lambda}}
=\int_0^{\infty}e^{-\nu s} \left(s^{\frac{1}{N}-1}e^{-\lambda s}
\sum_{m=0}^{\# K} \alpha_{-m} \,\theta_k^m (\lambda s)^{\frac{m}{N}}
E_{1,\frac{m+1}{N}}(\lambda s)\right) \mathrm{d}s.
$$
As a result, by introducing a convolution product, we obtain
\begin{align*}
\lqn{\int_0^{\infty} e^{-\lambda t} \left[\mathbb{E}\!\left(\!\!\vphantom{e^t}\right.\right.
e^{-\nu T(t)},\, X(t) \in \mathrm{d}x \!\!\left.\left.\vphantom{e^t}\right)\!/\mathrm{d}x\right]\mathrm{d}t}
&=
-\frac{1}{\lambda^{\frac{\# K -1}{N}}} \sum_{k \in K} B_k\theta_k e^{-\theta_k\!\!\sqrt[N]{\lambda}\,x}
\\
&\hphantom{=\,}
\times \int_0^{\infty} e^{-\nu s}\left(\int_0^s \frac{\sigma^{\frac{\# J-1}{N}-1}e^{-\lambda \sigma}}{\Gamma
\!\left(\vphantom{\frac aN}\right.\!\!\frac{\#J-1}{N}\!\!\left.\vphantom{\frac aN}\right)}
\times e^{-\lambda (s-\sigma)} \sum_{m=0}^{\# K} \alpha_{-m} \theta_k^m \lambda^{\frac m N}
(s-\sigma)^{\frac{m+1}{N}-1} E_{1,\frac{m+1}{N}} (\lambda(s-\sigma)) \,\mathrm{d}\sigma
\right)\mathrm{d}s.
\end{align*}
By removing the Laplace transforms with respect to the parameter $\nu$
of each member of the foregoing equality, we extract
\begin{align*}
\lqn{\int_0^{\infty} e^{-\lambda t} \left[\mathbb{P}\{T(t)\in \mathrm{d}s,
\, X(t) \in \mathrm{d}x\} /(\mathrm{d}s \,\mathrm{d}x)\right] \mathrm{d}t}
&=
-\frac{e^{-\lambda s}}{\lambda^{\frac{\# K -1}{N}}} \sum_{m=0}^{\#K} \alpha_{-m}\lambda^{\frac{m}{N}}
\left(\,\sum_{k \in K} B_k\theta_k^{m+1}  e^{-\theta_k\!\!\sqrt[N]{\lambda}\, x}\right)
\int_0^s \frac{\sigma^{\frac{\# J-1}{N}-1}}{\Gamma\!
\left(\vphantom{\frac aN}\right.\!\!\frac{\# J-1}{N}\!\!\left.\vphantom{\frac aN}\right)}
\,(s-\sigma)^{\frac{m+1}{N}-1} E_{1,\frac{m+1}{N}} (\lambda(s-\sigma))\,\mathrm{d}\sigma.
\end{align*}
The integral lying on the right-hand side of the previous equality can be evaluated as follows:
\begin{align*}
\int_0^s \frac{\sigma^{\frac{\# J-1}{N}-1}}{\Gamma\!
\left(\vphantom{\frac aN}\right.\!\!\frac{\# J-1}{N}\!\!\left.\vphantom{\frac aN}\right)}\,
(s-\sigma)^{\frac{m+1}{N}-1} E_{1,\frac{m+1}{N}} (\lambda(s-\sigma))\,\mathrm{d}\sigma
&=
\int_0^s \frac{\sigma^{\frac{\# J-1}{N}-1}}{\Gamma\!
\left(\vphantom{\frac aN}\right.\!\!\frac{\# J-1}{N}\!\!\left.\vphantom{\frac aN}\right)}
(s-\sigma)^{\frac{m+1}{N}-1} \sum_{\ell=0}^{\infty}
\frac{\lambda^{\ell} (s-\sigma)^{\ell}}{\Gamma\!
\left(\vphantom{\frac aN}\right.\!\!\ell+\frac{m+1}{N}\!\!\left.\vphantom{\frac aN}\right)}
\,\mathrm{d}\sigma
\\
&=
\sum_{\ell=0}^{\infty} \lambda^{\ell}
\int_0^s \frac{\sigma^{\frac{\# J-1}{N}-1}}{\Gamma\!
\left(\vphantom{\frac aN}\right.\!\! \frac{\# J-1}{N} \!\!\left.\vphantom{\frac aN}\right)}\,
\frac{(s-\sigma)^{\ell+\frac{m+1}{N}-1}}{\Gamma\!
\left(\vphantom{\frac aN}\right.\!\!\ell+\frac{m+1}{N}\!\!\left.\vphantom{\frac aN}\right)}
\,\mathrm{d}\sigma
\\
&=
\sum_{\ell=0}^{\infty} \frac{\lambda^{\ell} s^{\ell+\frac{m+\# J}{N}-1}}{\Gamma\!
\left(\vphantom{\frac aN}\right.\!\!\ell+\frac{m+\#J}{N}
\!\!\left.\vphantom{\frac aN}\right)}
=s^{\frac{m+\# J}{N}-1}E_{1,\frac{m+\#J}{N}} (\lambda s)
\end{align*}
from which we deduce~(\ref{wrt.nu}).
\EndDim

%
\begin{rem}
Let us integrate formula~(\ref{wrt.nu}) with respect to $x$ on $(-\infty,0]$. This gives
$$
\int_0^{\infty} e^{-\lambda t} \left[\mathbb{P} \{T(t) \in \mathrm{d}s,\,
X(t)\le 0\}/\mathrm{d}s\right] \mathrm{d}t
=\frac{e^{-\lambda s}}{(\lambda s)^{\frac{\# K}{N}}} \sum_{m=0}^{\# K} \alpha_{-m}\beta_m
(\lambda s)^{\frac{m}{N}} E_{1,\frac{m+\# J}{N}}(\lambda s).
$$
In view of~(\ref{set14}) and~(\ref{set15}), since $E_{1,1}(\lambda s)=e^{\lambda s}$,
\begin{align}
\int_0^{\infty} e^{-\lambda t}\left[\mathbb{P} \{T(t) \in \mathrm{d}s,\,
X(t) \le 0\}/\mathrm{d}s\right] \mathrm{d}t
&=
\frac{e^{-\lambda s}}{(\lambda s)^{\frac{\# K}{N}}} \left(\alpha_0\beta_0
E_{1,\frac{\# J}{N}}(\lambda s)+\alpha_{-\#K}\beta_{\#K}
(\lambda s)^{\frac{\# K}{N}} E_{1,1}(\lambda s)\right)
\nonumber\\
&=
\frac{e^{-\lambda s}}{(\lambda s)^{\frac{\# K}{N}}}  E_{1,\frac{\# J}{N}}(\lambda s)-1.
\label{star}
\end{align}
We can rewrite $E_{1,\frac{\# J}{N}}(\lambda s)$ as an integral by using
Lemma~\ref{8.2}. We obtain that
\begin{align*}
\lqn{\int_0^{\infty} e^{-\lambda t}\left[\mathbb{P} \{T(t)\in \mathrm{d}s,
\, X(t) \le 0\}/\mathrm{d}s\right] \mathrm{d}t}
&=
\frac{e^{-\lambda s}}{(\lambda s)^{\frac{\# K}{N}}}\, \frac{1}{(\lambda s)^{\frac{\# J}{N}-1}}
\left(e^{\lambda s}+\frac{\sin\!\left(\!\!\vphantom{\frac aN}\right.
\frac{\# J}{N}\pi\!\!\left.\vphantom{\frac aN}\right)}{\pi} \int_0^{\infty}
\frac{t^{1-\frac{\# J}{N}}}{t+1} \,e^{-\lambda s t}\,\mathrm{d}t\right)-1
\\
&=
\frac{\sin\!\left(\!\!\vphantom{\frac aN}\right. \frac{\# J}{N}\pi\!\!
\left.\vphantom{\frac aN}\right)}{\pi} e^{-\lambda s} \int_0^{\infty}
\frac{t^{\frac{\# K}{N}}}{t+1} \,e^{-\lambda s t}\,\mathrm{d}t
=
\frac{\sin\!\left(\!\!\vphantom{\frac aN}\right. \frac{\# J}{N}\pi\!\!
\left.\vphantom{\frac aN}\right)}{\pi s^{\frac{\# K}{N}}} e^{-\lambda s}
\int_0^{\infty} \frac{t^{\frac{\# K}{N}}}{t+s} \,e^{-\lambda t}\,\mathrm{d}t
\\
&=
\frac{\sin\!\left(\!\!\vphantom{\frac aN}\right. \frac{\# K}{N}\pi\!\!
\left.\vphantom{\frac aN}\right)}{\pi s^{\frac{\# K}{N}}}
\int_s^{\infty} \frac{(t-s)^{\frac{\# K}{N}}}{t} \,e^{-\lambda t}\,\mathrm{d}t.
\end{align*}
From this, we extract, for $0 < s <t$,
\begin{equation}
\mathbb{P}\{T(t) \in \mathrm{d}s,\, X(t) \le 0\}/\mathrm{d}s
=\frac{\sin\!\left(\!\!\vphantom{\frac aN}\right. \frac{\# K}{N}\pi\!\!
\left.\vphantom{\frac aN}\right)}{\pi t} \left(\frac{t-s}{s}\right)^{\!\frac{\# K}{N}}.
\label{f1}
\end{equation}
We retrieve Theorem 14 of \cite{2003}.
By integrating~(\ref{f1}) with respect to $s$, we get
$$
\mathbb{P} \{ X(t) \le 0\} =\frac{\sin\!\left(\!\!\vphantom{\frac aN}\right.
\frac{\# K}{N}\pi\!\!\left.\vphantom{\frac aN}\right)}{\pi t}
\int_0^t \left(\frac{t-s}{s}\right)^{\!\frac{\# K}{N}}\mathrm{d}s
=\frac{\sin\!\left(\!\!\vphantom{\frac aN}\right. \frac{\# K}{N}\pi\!\!
\left.\vphantom{\frac aN}\right)}{\pi}
B\!\left(\frac{\# K}{N}+1,1-\frac{\# K}{N}\right)
=\frac{\Gamma\!\left(\!\!\vphantom{\frac aN}\right. \frac{\# K}{N}+1\!\!
\left.\vphantom{\frac aN}\right)}{\Gamma\!\left(\!\!\vphantom{\frac aN}\right.
\frac{\# K}{N}\!\!\left.\vphantom{\frac aN}\right)}
$$
which simplifies to $\mathbb{P} \{ X(t) \le 0\} =\# K/N.$
We retrieve (11) of \cite{2003}.
\end{rem}
%
\begin{rem}
An alternative expression for formula~(\ref{wrt.nu}) is for $x\le 0$
\begin{align}
\lqn{\int_0^{\infty} e^{-\lambda t}\left[\mathbb{P}\{T(t)\in \mathrm{d}s, X(t) \in \mathrm{d}x\}
/(\mathrm{d}s \,\mathrm{d}x)\right] \mathrm{d}t}
=-\frac{e^{-\lambda s}}{\lambda^{\frac{\# K-1}{N}}s^{\frac{\# K}{N}}}
\sum_{\substack{j \in J \\k \in K}} A_j B_k\theta_k \,E_{\frac{1}{N},\frac{\# J}{N}}\!
\left(\frac{\theta_k}{\theta_j} \!\sqrt[N]{\lambda s}\right) e^{-\theta_k\!\!\sqrt[N]{\lambda}\,x}.
\label{wrt.nu2}
\end{align}
In effect, by~(\ref{wrt.nu}),
\begin{align*}
\lqn{\int_0^{\infty} e^{-\lambda t}\left[\mathbb{P}\{T(t)\in \mathrm{d}s,
X(t) \in \mathrm{d}x\}  /(\mathrm{d}s \,\mathrm{d}x)\right] \mathrm{d}t}
&=
-\frac{e^{-\lambda s}}{\lambda^{\frac{\# K-1}{N}}s^{\frac{\# K}{N}}} \sum_{\ell=0}^{\infty}
\sum_{m=0}^{\# K} \sum_{k\in K} \alpha_{-m} B_k\theta_k^{m+1}
\frac{(\lambda s)^{\ell+\frac{m}{N}}}{\Gamma\!
\left(\vphantom{\frac aN}\right.\!\!\ell+\frac{m+\# J}{N}\!\!\left.\vphantom{\frac aN}\right)}
\,e^{-\theta_k\!\!\sqrt[N]{\lambda}\, x}
\\
&=
-\frac{e^{-\lambda s}}{\lambda^{\frac{\# K-1}{N}}s^{\frac{\# K}{N}}}
\sum_{\ell=0}^{\infty} \sum_{m=0}^{N-1} \sum_{\substack{j \in J \\k \in K}} A_j B_k\theta_k
\left(\frac{\theta_k}{\theta_j}\right)^{\!m} \frac{(\lambda s)^{\ell+\frac{m}{N}}}{\Gamma\!
\left(\vphantom{\frac aN}\right.\!\!\ell+\frac{m+\# J}{N}\!\!\left.\vphantom{\frac aN}\right)}
\,e^{-\theta_k\!\!\sqrt[N]{\lambda}\, x}.
\end{align*}
In the last displayed equality, we have extended the sum with respect to $m$
to the range $0\le m\le N-1$ because, by~(\ref{set14}), the $\alpha_{-m}$,
$\# K+1 \le m \le N-1$, vanish. Let us introduce the index $r=\ell N+m$.
Since $\left(\frac{\theta_k}{\theta_j}\right)^{\!m}=\left(\frac{\theta_k}{\theta_j}\right)^{\!r}$,
we have
$$
\int_0^{\infty} e^{-\lambda t}\left[\mathbb{P}\{T(t)\in \mathrm{d}s,
X(t) \in \mathrm{d}x\} /(\mathrm{d}s \,\mathrm{d}x)\right] \mathrm{d}t
=-\frac{e^{-\lambda s}}{\lambda^{\frac{\# K-1}{N}}s^{\frac{\# K}{N}}}
\sum_{\substack{j \in J \\k \in K}} A_j B_k\theta_k \sum_{r=0}^{\infty}
\frac{\left(\frac{\theta_k}{\theta_j} \!\sqrt[N]{\lambda s}\right)^{\!r}}{\Gamma\!
\left(\vphantom{\frac aN}\right.\!\!\frac{r+\# J}{N}\!\!\left.\vphantom{\frac aN}\right)}
\,e^{-\theta_k\!\!\sqrt[N]{\lambda}\, x}
$$
which coincide with~(\ref{wrt.nu2}).
\end{rem}
%
\begin{ex}
Case $N=2$. Suppose $x\le 0$. The first expression~(\ref{wrt.nu}) reads
$$
\int_0^{\infty} e^{-\lambda t}\left[\mathbb{P}\{T(t) \in \mathrm{d}s, X(t)\in \mathrm{d}x\}
/ (\mathrm{d}s \,\mathrm{d}x)\right] \mathrm{d}t
=\frac{e^{-\lambda s}}{\sqrt s} \left(E_{1 ,\frac12} (\lambda s) -\sqrt{\lambda s}\,
E_{1,1}(\lambda s)\right) e^{\sqrt \lambda \, x}
$$
while the second expression~(\ref{wrt.nu2}) reads
$$
\int_0^{\infty} e^{-\lambda t} \left[\mathbb{P}\{T(t) \in \mathrm{d}s, X(t)\in \mathrm{d}x\} /
(\mathrm{d}s \,\mathrm{d}x)\right] \mathrm{d}t
=\frac{e^{-\lambda s}}{\sqrt s} \,
E_{\frac12,\frac12} \!\left(\vphantom{\sqrt a}\right.\!\!
-\sqrt{\lambda s}\!\left.\vphantom{\sqrt a}\right) e^{\sqrt \lambda \, x}.
$$
From Lemma~\ref{2half}, we have
\begin{equation}
E_{\frac12,\frac12} \!\left(\vphantom{\sqrt a}\right.\!\!
-\sqrt{\lambda s}\!\left.\vphantom{\sqrt a}\right)\!
=E_{1,\frac12}(\lambda s)-\sqrt{\lambda s} \,e^{\lambda s}
\label{kub}
\end{equation}
which proves the coincidence of both formulas. Moreover, from Lemma~\ref{2half},
for $x \le 0$,
$$
\int_0^{\infty} e^{-\lambda t}\left[\mathbb{P}\{T(t) \in \mathrm{d}s,
X(t)\in \mathrm{d}x\} /(\mathrm{d}s \,\mathrm{d}x)\right] \mathrm{d}t
=\left(\frac{e^{-\lambda s}}{\sqrt{\pi s}}-\sqrt{\lambda}
\,\mathrm{Erfc} (\sqrt{\lambda s})\right)
e^{\sqrt \lambda \,x}.
$$
We retrieve formula 1.4.6, page 129, of \cite{bs}.
\end{ex}
%
\begin{ex}
Case $N=3$. We have for $x\le 0$, when $\kappa_3=-1$:
$$
\int_0^{\infty} e^{-\lambda t} \left[\mathbb{P}\{ T(t)\in \mathrm{d}s,\, X(t) \in \mathrm{d}x\}
/(\mathrm{d}s \,\mathrm{d}x)\right] \mathrm{d}t=
e^{\sqrt[3]{\lambda} \,x} \left(\frac{e^{-\lambda s}}{\sqrt[3]{s}}\,
E_{1,\frac23}(\lambda s)-\sqrt[3] \lambda\,\right)
$$
and when $\kappa_3=1$:
\begin{align*}
\lqn{\int_0^{\infty} e^{-\lambda t} \left[\mathbb{P}\{ T(t)\in \mathrm{d}s,\, X(t) \in \mathrm{d}x\}
/(\mathrm{d}s \,\mathrm{d}x)\right]\mathrm{d}t}
&=
\frac{e^{-\lambda s+\frac{\sqrt[3]{\lambda}}{2}x}}{\sqrt 3\,\sqrt[3]{\lambda s^2}}
\left[\vphantom{\frac tt}\right.\!\!\sqrt 3 \,\cos\!\left(\vphantom{\frac tt}\right.\!\!
\frac{\sqrt 3 \,\sqrt[3] \lambda \,x}{2}\!\!\left.\vphantom{\frac tt}\right)\!\!
\left(\sqrt[3]{\lambda s}\,
E_{1,\frac23}(\lambda s)-(\lambda s)^{2/3} e^{\lambda s}\right)
\\
&\hphantom{=\,}
+\sin\!\left(\vphantom{\frac tt}\right.\!\!
\frac{\sqrt 3 \,\sqrt[3] \lambda \,x}{2}\!\!\left.\vphantom{\frac tt}\right)\!\!
\left(\sqrt[3]{\lambda s}\,
E_{1,\frac23}(\lambda s)+(\lambda s)^{2/3} e^{\lambda s}-2E_{1,\frac13} (\lambda s)\right)
\!\!\left.\vphantom{\frac tt}\right]\!.
\end{align*}
\end{ex}
%
\begin{ex}
Case $N=4$. We have, for $x\ge 0$,
\begin{align*}
\lqn{\int_0^{\infty} e^{-\lambda t} \left[\mathbb{P}\{ T(t)\in \mathrm{d}s,\, X(t) \in \mathrm{d}x\}
/(\mathrm{d}s \,\mathrm{d}x)\right] \mathrm{d}t}
&=
\frac{\sqrt 2\, e^{-\lambda s+\frac{\sqrt[4]\lambda}{\sqrt 2}x}}{\sqrt[4]\lambda \,\sqrt s}
\left[\vphantom{\frac tt}\right.\!\!\cos\!\left(\vphantom{\frac tt}\right.\!\!
\frac{\sqrt[4]\lambda \,x}{\sqrt 2}\!\!\left.\vphantom{\frac tt}\right)\!\!
\left(\sqrt[4]{\lambda s}\,E_{1,\frac34}(\lambda s) -\sqrt{\lambda s} \,e^{\lambda s}\right)
+\sin\!\left(\vphantom{\frac tt}\right.\!\!
\frac{\sqrt[4]\lambda \,x}{\sqrt 2}\!\!\left.\vphantom{\frac tt}\right)\!\!
\left(\sqrt[4]{\lambda s}\,E_{1,\frac34}(\lambda s) -E_{1,\frac12}(\lambda s)\right)
\!\!\left.\vphantom{\frac tt}\right]\!.
\end{align*}
\end{ex}
%

\subsection{The case $x\ge 0$}
%
\begin{teo}\label{theorem-wrt.nu+2}
The Laplace transform with respect to $t$ of the density of the couple
$(T(t),X(t))$ is given, when $x\ge 0$, by
\begin{align}
\lqn{\int_0^{\infty} e^{-\lambda t} [\mathbb{P}\{ T(t)\in \mathrm{d}s,\, X(t) \in \mathrm{d}x\}
/(\mathrm{d}s \,\mathrm{d}x)] \,\mathrm{d}t}
&=
-\frac{e^{-\lambda s}}{\lambda^{\frac{\# K-1}{N}}} \sum_{\substack{j \in J \\k \in K}}
A_j B_k\theta_k \int_0^s \sigma^{\frac 1N-1} E_{\frac 1N,\frac 1N}\!
\left(\frac{\theta_k}{\theta_j} \sqrt[N]{\lambda \sigma}\right)
I_{j,\# J-1}(s-\sigma;x)\,\mathrm{d}\sigma
\label{wrt.nu+1}
\end{align}
where the function $I_{j,\# J-1}$ is defined by~(\ref{set18}).
\end{teo}
%
\Dim
Recall~(\ref{wrt.mu}) in the case $x \ge 0$:
\begin{align*}
\lqn{\int_0^{\infty} e^{-\lambda t} \left[\mathbb{E}\!\left(\!\!\vphantom{e^t}\right.\right.
e^{-\nu T(t)},\, X(t) \in \mathrm{d}x \!\!\left.\left.\vphantom{e^t}\right)\!/\mathrm{d}x\right]\mathrm{d}t}
&=
\frac{1}{\lambda^{\frac{\# K-1}{N}}(\lambda+\nu)^{\frac{\# J-1}{N}}}
\sum_{j \in J} A_j \theta_j \left(\,\sum_{k \in K}
\frac{B_k\theta_k}{\theta_k\!\!\sqrt[N]{\lambda}-\theta_j\!\!\sqrt[N]{\lambda+\nu}}\right)
e^{-\theta_j\!\!\sqrt[N]{\lambda+\nu}\,x}.
\end{align*}
We have to invert the quantity
$\frac{e^{-\theta_j\!\!\sqrt[N]{\lambda+\nu}\,x}}{(\lambda+\nu)^{\frac{\# J-1}{N}}
\left(\vphantom{\sqrt t}\right.\!\!\!\sqrt[N]{\lambda+\nu}
-\frac{\theta_k}{\theta_j} \!\sqrt[N]{\lambda}\left.\vphantom{\sqrt t}\right)}$ with respect to $\nu$.
Recalling~(\ref{mittag}) and~(\ref{set24}),
\begin{align*}
\frac{1}{\!\sqrt[N]{\lambda +\nu}-\beta}
&=
\int_0^{\infty} e^{-\nu s}\left(s^{\frac{1}{N}-1}
e^{-\lambda s} E_{\frac 1N,\frac 1N} \!\left(\beta \sqrt[N]{s}\,\right)\right) \mathrm{d}s,
\\
\frac{e^{-\theta_j\!\!\sqrt[N]{\lambda+\nu}\,x}}{(\lambda+\nu)^{\frac{\# J-1}{N}}}
&=
\int_0^{\infty} e^{-\nu s}\left(e^{-\lambda s} I_{j,\# J-1} (s;x)\right)\mathrm{d}s,
\end{align*}
we get by convolution
\begin{align*}
\lqn{\frac{e^{-\theta_j\!\!\sqrt[N]{\lambda+\nu}\,x}}{(\lambda+\nu)^{\frac{\# J-1}{N}}
\left(\vphantom{\sqrt t}\right.\!\!\!\sqrt[N]{\lambda+\nu}-\frac{\theta_k}{\theta_j} \!\sqrt[N]{\lambda}
\left.\vphantom{\sqrt t}\right)}}
&=
\int_0^{\infty} e^{-\nu s} \left(\int_0^s \sigma^{\frac 1N-1} e^{-\lambda \sigma}
E_{\frac 1N,\frac1N} \!\left(\frac{\theta_k}{\theta_j} \!\sqrt[N]{\lambda \sigma}\right)
\times e^{-\lambda (s-\sigma)} I_{j,\# J-1} (s-\sigma;x)\,\mathrm{d}\sigma\right)\mathrm{d}s
\\
\lqn{}
&=
\int_0^{\infty} e^{-\nu s} \left(e^{-\lambda s} \int_0^s \sigma^{\frac 1N-1}
E_{\frac 1N,\frac1N} \!\left(\frac{\theta_k}{\theta_j} \!\sqrt[N]{\lambda \sigma}\right)
I_{j,\# J-1} (s-\sigma;x)\,\mathrm{d}\sigma\right)\mathrm{d}s.
\end{align*}
This immediately yields~(\ref{wrt.nu+1}).
\EndDim

%
\begin{rem} Noticing that
$$
E_{\frac 1N,\frac 1N} \!\left(\frac{\theta_k}{\theta_j} \!\sqrt[N]{\lambda \sigma}\right)
=\sum_{r=0}^{\infty} \frac{\theta_k^r}{\theta_j^r} \,
\frac{(\lambda \sigma)^{\frac{r}{N}}}{\Gamma\!\left(\frac{r+1}{N}\right)}
=\sum_{\ell=0}^{\infty} \sum_{m=0}^{N-1} \frac{\theta_k^m}{\theta_j^m}\,
\frac{(\lambda \sigma)^{\ell+\frac{m}{N}}}{\Gamma\!\left(\ell+\frac{m+1}{N}\right)}
=\sum_{m=0}^{N-1} \frac{\theta_k^m}{\theta_j^m} \,(\lambda \sigma)^{\frac{m}{N}}
E_{1,\frac{m+1}{N}}(\lambda \sigma)
$$
and reminding that, from~(\ref{set26}), the $\beta_m$, $1\le m \le \# K-1$, vanish,
we can rewrite~(\ref{wrt.nu+1}) in the following form. For $x\ge 0$,
\begin{align}
\lqn{\int_0^{\infty} e^{-\lambda t} \left[\mathbb{P}\{ T(t)\in \mathrm{d}s,\, X(t) \in \mathrm{d}x\}
/(\mathrm{d}s \,\mathrm{d}x)\right] \mathrm{d}t}
&=
-e^{-\lambda s} \sum_{m=\# K -1}^{N-1} \left(\,\sum_{k \in K} B_k\theta_k^{m+1}\right)
\lambda^{\frac{m-\# K +1}{N}} \int_0^s \sigma^{\frac{m+1}{N}-1}
E_{1,\frac{m+1}{N}} (\lambda \sigma) \left(\!\vphantom{\sum_{\in}}\right.\sum_{j \in J}
\frac{A_j}{\theta_j^m}I_{j,\#J-1}(s-\sigma;x) \left.\vphantom{\sum_{\in}}\!\right)\mathrm{d}\sigma
\nonumber \\
&=
-e^{-\lambda s} \sum_{m=\# K}^{N} \beta_m\, \lambda^{\frac{m-\# K}{N}} \int_0^s
\sigma^{\frac{m}{N}-1} E_{1,\frac{m}{N}} (\lambda \sigma) \,\Phi_m(s-\sigma;x) \,\mathrm{d}\sigma
\label{wrt.nu+2}
\end{align}
with $\Phi_m(\tau;x)= \sum_{j \in J} \frac{A_j}{\theta_j^{m-1}} \,I_{j,\#J-1}(\tau;x)$.
\end{rem}
%
\begin{rem}
Let us integrate~(\ref{wrt.nu+2}) with respect to $x$ on $[0,\infty)$. We first compute
\begin{align*}
\int_0^{\infty} \Phi_m(\tau;x) \,\mathrm{d}x
&=
\frac{N i}{2 \pi}\left(\!\vphantom{\sum_{\in}}\right.\sum_{j \in J}
\frac{A_j}{\theta_j^m}\left.\!\vphantom{\sum_{\in}}\right)\!\!
\left(e^{-i\frac{\#J}{N}\pi}-e^{i\frac{\#J}{N}\pi}\right)
\int_0^{\infty} \xi^{\# K-1} \,e^{-\tau \xi^N} \,\mathrm{d}\xi
\\
&=
\frac{\Gamma\!\left(\!\!\vphantom{\frac aN}\right.\frac{\# K}{N}\!\!
\left.\vphantom{\frac aN}\right)\sin\!\left(\!\!\vphantom{\frac aN}\right.
\frac{\# J}{N}\pi\!\!\left.\vphantom{\frac aN}\right)}{\pi \tau^{\frac{\# K}{N}}}\,\alpha_{-m}
=\frac{\alpha_{-m}}{\Gamma\!\left(\!\!\vphantom{\frac aN}\right. \frac{\# J}{N}\!\!\left.\vphantom{\frac aN}\right) \tau^{\frac{\# K}{N}}}.
\end{align*}
Then, with the aid of~(\ref{set26}) and~(\ref{set15}), we get
\begin{align*}
\lqn{\int_0^{\infty} e^{-\lambda t} \left[\mathbb{P}\{ T(t)\in \mathrm{d}s,\, X(t)\ge 0\}
/ \mathrm{d}s\right] \mathrm{d}t}
&=
-\frac{e^{-\lambda s}}{\Gamma\!\left(\!\!\vphantom{\frac aN}\right.\frac{\#J}{N}\!\!
\left.\vphantom{\frac aN}\right)}\left(\alpha_{-\# K}\,
\beta_{\#K} \int_0^s\frac{\sigma^{\frac{\#K}{N}-1}}{(s-\sigma)^{\frac{\# K}{N}}}
E_{1,\frac{\# K}{N}}(\lambda \sigma)\,\mathrm{d}\sigma + \alpha_{-N}\,\beta_{N} \,
\lambda^{\frac{\#J}{N}} \int_0^s\frac{E_{1,1}(\lambda \sigma)}{(s-\sigma)^{\frac{\# K}{N}}}
\,\mathrm{d}\sigma\right)
\\
&=
\frac{e^{-\lambda s}}{\Gamma\!\left(\!\!\vphantom{\frac aN}\right. \frac{\# J}{N}\!\!\left.\vphantom{\frac aN}\right)}\left(\int_0^s
\frac{E_{1,\frac{\# K}{N}}(\lambda \sigma)}{\sigma^{\frac{\#J}{N}}  (s-\sigma)^{\frac{\#K}{N}}}
\,\mathrm{d}\sigma-\lambda^{\frac{\#J}{N}}
\int_0^s\frac{e^{\lambda \sigma}}{(s-\sigma)^{\frac{\# K}{N}}}\,\mathrm{d}\sigma\right)
\\
&=
\frac{e^{-\lambda s}}{\Gamma\!\left(\!\!\vphantom{\frac aN}\right. \frac{\# J}{N}\!\!\left.\vphantom{\frac aN}\right)} \sum_{\ell=0}^{\infty}
\frac{\lambda^{\ell}}{\Gamma\!
\left(\vphantom{\frac aN}\right.\!\!\ell+\frac{\# K}{N}\!\!\left.\vphantom{\frac aN}\right)}
\int_0^s \frac{\sigma^{\ell-\frac{\#J}{N}}}{(s-\sigma)^{\frac{\#K}{N}}} \,\mathrm{d}\sigma
-\frac{\lambda^{\frac{\#J}{N}}}{\Gamma\!\left(\!\!\vphantom{\frac aN}\right. \frac{\# J}{N}\!\!\left.\vphantom{\frac aN}\right)}
\int_0^s \frac{e^{-\lambda( s-\sigma)}}{(s-\sigma)^{\frac{\# K}{N}}} \,\mathrm{d}\sigma
\\
&=
e^{-\lambda s} \sum_{\ell=0}^{\infty} \frac{(\lambda s)^{\ell}}{\ell!}
-\frac{1}{\Gamma\!\left(\!\!\vphantom{\frac aN}\right. \frac{\# J}{N}\!\!\left.\vphantom{\frac aN}\right)}
\int_0^{\lambda s} \frac{e^{-\sigma}}{\sigma^{\frac{\# K}{N}}} \,\mathrm{d}\sigma
=
1-\frac{1}{\Gamma\!\left(\!\!\vphantom{\frac aN}\right. \frac{\# J}{N}\!\!\left.\vphantom{\frac aN}\right)}
\int_0^{\lambda s} \sigma^{\frac{\# J}{N}-1} e^{-\sigma} \,\mathrm{d}\sigma.
\end{align*}
By Lemma~\ref{1alpha}, this simplifies into
\begin{equation}
\int_0^{\infty} e^{-\lambda t} \left[\mathbb{P}\{ T(t)\in \mathrm{d}s,\, X(t)\ge 0\}/
\mathrm{d}s\right] \mathrm{d}t
=1-(\lambda s)^{\frac{\# J}{N}} e^{-\lambda s} E_{1,\frac{\# J}{N}+1}(\lambda s).
\label{TLXpositive}
\end{equation}
Now, using  Lemma~\ref{8.2}, we derive another representation for
the foregoing Laplace transform:
\begin{align*}
\int_0^{\infty}e^{-\lambda t}\left[\mathbb{P}\{T(t)\in \mathrm{d}s, X(t) \ge 0\}/
\mathrm{d}s\right] \mathrm{d}t
&=
1-e^{-\lambda s} \left(e^{\lambda s}-\frac{\sin\!\left(\!\!\vphantom{\frac aN}\right. \frac{\# J}{N}\pi\!\!\left.\vphantom{\frac aN}\right)}{\pi}
\int_0^{\infty} \frac{t^{-\frac{\# J}{N}}}{t+1} \,e^{-\lambda s t}\,\mathrm{d}t\right)
\\
&=
\frac{\sin\!\left(\!\!\vphantom{\frac aN}\right. \frac{\# J}{N}\pi\!\!\left.\vphantom{\frac aN}\right)}{\pi} \,s^{\frac{\# J}{N}}
\int_0^{\infty} \frac{e^{-\lambda(t+s)}}{(t+s) \,t^{\frac{\# J}{N}}} \,\mathrm{d}t
\\
&=
\frac{\sin\!\left(\!\!\vphantom{\frac aN}\right. \frac{\# J}{N}\pi\!\!\left.\vphantom{\frac aN}\right)}{\pi} \,s^{\frac{\# J}{N}}
\int_{s}^{\infty} \frac{e^{-\lambda t}}{t (t-s)^{\frac{\# J}{N}}} \,\mathrm{d}t.
\end{align*}
As a result, we derive
\begin{equation}
\mathbb{P}\{T(t)\in \mathrm{d}s, X(t) \ge 0\}/ \mathrm{d}s
=\frac{\sin\!\left(\!\!\vphantom{\frac aN}\right. \frac{\# J}{N}\pi\!\!\left.\vphantom{\frac aN}\right)}{\pi t}
\left(\frac{s}{t-s}\right)^{\!\frac{\# J}{N}},\qquad 0<s<t.
\label{f2}
\end{equation}
This is formula~(11) of~\cite{2003}.
\end{rem}
%
\begin{rem}
If we add~(\ref{star}) and~(\ref{TLXpositive}), we find
$$
\int_0^{\infty}e^{-\lambda t}\left[\mathbb{P}\{T(t) \in \mathrm{d}s\} /\mathrm{d}s\right] \mathrm{d}t
=e^{-\lambda s} \left(\sum_{\ell=0}^{\infty}
\frac{(\lambda s)^{\ell+\frac{\# J}{N}-1}}{\Gamma\!
\left(\vphantom{\frac aN}\right.\!\!\ell+\frac{\# J}{N}\!\!\left.\vphantom{\frac aN}\right)}
-\sum_{\ell=0}^{\infty} \frac{(\lambda s)^{\ell+\frac{\# J}{N}}}{\Gamma\!
\left(\vphantom{\frac aN}\right.\!\! \ell+\frac{\# J}{N}+1 \!\!
\left.\vphantom{\frac aN}\right)}\right)
=\frac{e^{-\lambda s}}{\Gamma\!
\left(\vphantom{\frac aN}\right.\!\!\frac{\# J}{N}\!\!\left.\vphantom{\frac aN}\right) \!(\lambda s)^{\frac{\# K}{N}}}.
$$
It is easy to invert this Laplace transform. Indeed,
\begin{align*}
\frac{e^{-\lambda s}}{\lambda^{\frac{\#K}{N}}}=\frac{e^{-\lambda s}}{\Gamma\!
\left(\vphantom{\frac aN}\right.\!\!\frac{\# K}{N} \!\!\left.\vphantom{\frac aN}\right)}
\int_0^{\infty} t^{\frac{\#K}{N}-1} e^{-\lambda t}\,\mathrm{d}t
&=
\frac{1}{\Gamma\!\left(\vphantom{\frac aN}\right.\!\frac{\# K}{N}
\!\left.\vphantom{\frac aN}\right)} \int_{s}^{\infty} (t-s)^{\frac{\# K}{N}-1} e^{-\lambda t}\,\mathrm{d}t
\\
&=
\frac{1}{\pi}\,\Gamma\!\left(\frac{\# J}{N}\right) \,\sin\!\left(\frac{\# J}{N}\pi\right)
\int_{s}^{\infty} \frac{e^{-\lambda t}}{(t-s)^{\frac{\# J}{N}}} \,\mathrm{d}t.
\end{align*}
This implies that
$$
\mathbb{P}\{T(t) \in \mathrm{d}s\} /\mathrm{d}s
=\frac{\sin\!\left(\!\!\vphantom{\frac aN}\right. \frac{\# J}{N}\pi\!\!\left.\vphantom{\frac aN}\right)}{\pi}
\frac{\ind_{(0,t)}(s)}{ s^{\frac{\# K}{N}}(t-s)^{\frac{\# J}{N}}}
$$
which can be also obtained by adding directly~(\ref{f1}) and~(\ref{f2}).
Thus, we retrieve the famous counterpart to the Paul L\'evy's arc-sine
law stated in \cite{2003} (Corollary 9).
\end{rem}
%
\begin{rem}
For $x=0$, using formula~(\ref{wrt.nu}) which is valid for $x\le0$, we get,
by~(\ref{set26}), (\ref{set14}) and~(\ref{set15}),
\begin{align}
\lqn{\int_0^{\infty}e^{-\lambda t} \,\mathbb{P}\{T(t) \in \mathrm{d}s,\, X(t)\in \mathrm{d}x\}
/(\mathrm{d}s \,\mathrm{d}x) \Big|_{x=0}\,\mathrm{d}t}
&=
-\frac{e^{-\lambda s}}{\lambda^{\frac{\#K-1}{N}} s^{\frac{\# K}{N}}}
\sum_{m=0}^{\# K} \alpha_{-m}\,\beta_{m+1}\, (\lambda s)^{\frac mN} E_{1,\frac{m+\# J}{N}} (\lambda s)
\nonumber\\
&=
-\frac{e^{-\lambda s}}{\lambda^{\frac{\#K-1}{N}} s^{\frac{\# K}{N}}}
\left(\alpha_{1-\#K}\,\beta_{\#K}\, (\lambda s)^{\frac{\#K-1}{N}} E_{1,1-\frac1N}(\lambda s)
+\alpha_{-\#K}\,\beta_{\#K+1}\, (\lambda s)^{\frac{\#K}{N}} E_{1,1}(\lambda s)\right)
\nonumber
\\
&=
\frac{e^{-\lambda s}}{\lambda^{\frac{\#K-1}{N}} s^{\frac{\# K}{N}}}
\left(\!\vphantom{\sum_\in}\right.\sum_{j \in J} \theta_j (\lambda s)^{\frac{\# K-1}{N}}
E_{1,1-\frac1N}(\lambda s)+\sum_{k \in K} \theta_k (\lambda s)^{\frac{\#K}{N}}
e^{\lambda s}\left.\!\!\vphantom{\sum_\in}\right)
\nonumber\\
&=
\left(\!\vphantom{\sum_{\in}}\right.\sum_{j \in J} \theta_j
\left.\!\vphantom{\sum_{\in}}\right) \frac{e^{-\lambda s}}{\!\!\sqrt[N]{s}}
\left(E_{1,1-\frac 1N}(\lambda s)-\sqrt[N]{\lambda s}\,e^{\lambda s}\right)\!.
\label{LTX-zero1}
\end{align}
On the other hand, with formula~(\ref{wrt.nu+2}) which is valid for $x \ge 0$,
\begin{align}
\lqn{\int_0^{\infty}e^{-\lambda t} \,\mathbb{P}\{T(t) \in \mathrm{d}s,
\,X(t)\in \mathrm{d}x\}/(\mathrm{d}s\,\mathrm{d}x )\Big|_{x=0} \,\mathrm{d}t}
&=
-e^{-\lambda s} \sum_{m=\# K}^{N} \beta_m\,\lambda^{\frac{m-\# K}{N}}
\int_0^s \sigma^{\frac{m}{N}-1} E_{1,\frac mN} (\lambda \sigma) \,\Phi_m(s-\sigma;0) \,\mathrm{d}\sigma
\label{LTX-zero2}
\end{align}
with
\begin{align*}
\Phi_m(\tau;0)
&=
\frac{Ni}{2 \pi} \left(\!\vphantom{\sum_{\in}}\right.
\sum_{j \in J} \frac{A_j}{\theta_j^{m-1}} \left.\!\vphantom{\sum_{\in}}\right)\!\!
\left(e^{-i \frac{\# J-1}{N}\pi}-e^{i \frac{\#J-1}{N}\pi}\right)
\int_0^{\infty} \xi^{\# K} e^{-\tau \xi^N} \,\mathrm{d}\xi
\\
&=
\frac{\Gamma\!\left(\vphantom{\frac aN}\right.\!\!
\frac{\# K+1}{N} \!\!\left.\vphantom{\frac aN}\right)
\sin\!\left(\!\!\vphantom{\frac aN}\right. \frac{\# J-1}{N}\pi
\!\!\left.\vphantom{\frac aN}\right)}{\pi\,\tau^{\frac{\#K+1}{N}}}\,\alpha_{1-m}
= \frac{\alpha_{1-m}}{\Gamma\!\left(\vphantom{\frac aN}\right.\!\!\frac{\#J-1}{N}
\!\!\left.\vphantom{\frac aN}\right) \tau^{\frac{\#K+1}{N}}}.
\end{align*}
In view of~(\ref{set26}), (\ref{set14}) and~(\ref{set15}), we have
\begin{align*}
\lqn{\int_0^{\infty} e^{-\lambda t} \,\mathbb{P}\{T(t) \in \mathrm{d}s,
\, X(t)\in \mathrm{d}x\}/(\mathrm{d}s\,\mathrm{d}x ) \Big|_{x=0} \,\mathrm{d}t}
&=
\frac{e^{-\lambda s}}{\Gamma\!\left(\vphantom{\frac aN}\right.\!\!\frac{\# J-1}{N}
\!\!\left.\vphantom{\frac aN}\right)} \left[\left(\!\vphantom{\sum_{\in}}\right.
\sum_{j \in J} \theta_j \!\left.\vphantom{\sum_{\in}}\right)
\int_0^s \frac{\sigma^{\frac{\#K}{N}-1}}{(s-\sigma)^{\frac{\#K+1}{N}}}\,
E_{1,\frac{\#K}{N}} (\lambda \sigma)\,\mathrm{d}\sigma\right.
\\
&\hphantom{=\,}
\left.+ \left(\,\sum_{k \in K} \theta_k\right) \!\!\sqrt[N]{\lambda}
\int_0^s \frac{\sigma^{\frac{\#K+1}{N}-1}}{(s-\sigma)^{\frac{\#K+1}{N}}}\,
E_{1,\frac{\#K+1}{N}} (\lambda \sigma) \,\mathrm{d}\sigma\right]
\\
&=
\left(\!\vphantom{\sum_{\in}}\right.\sum_{j \in J} \theta_j\!\left.\vphantom{\sum_{\in}}\right)
\frac{e^{-\lambda s}}{\Gamma\!\left(\vphantom{\frac aN}\right.\!\!
\frac{\#J-1}{N} \!\!\left.\vphantom{\frac aN}\right)} \,e^{-\lambda s}
\left(\,\sum_{\ell=0}^{\infty} \frac{B\!\left(\!\!\vphantom{\frac aN}\right.
\ell+\frac{\# K}{N}, 1-\frac{\# K+1}{N}\!\!\left.\vphantom{\frac aN}\right)}{\Gamma\!\left(\vphantom{\frac aN}\right.\!\!
\ell+\frac{\# K}{N}\!\!\left.\vphantom{\frac aN}\right)}\,
\lambda^{\ell}s^{\ell-\frac 1N}\right.
\\
&\hphantom{=\,}
\left.-\!\sqrt[N]\lambda \,\sum_{\ell=0}^{\infty} \frac{B\!\left(\!\!\vphantom{\frac aN}\right.
\ell+\frac{\# K+1}{N}, 1-\frac{\# K+1}{N}\!\!\left.\vphantom{\frac aN}\right)}{\Gamma
\!\left(\vphantom{\frac aN}\right.\!\!\ell+\frac{\# K+1}{N}\!\!\left.\vphantom{\frac aN}\right)}\,
(\lambda s)^{\ell}\right)
\\
&=
\left(\!\vphantom{\sum_{\in}}\right.\sum_{j \in J} \theta_j
\left.\!\vphantom{\sum_{\in}}\right) \frac{e^{-\lambda s}}{\!\!\sqrt[N]s}
\left(E_{1, 1-\frac 1N}(\lambda s)-\sqrt[N]{\lambda s} \,e^{\lambda s}\right)\!.
\end{align*}
Thus, we have checked that the two different formulas~(\ref{LTX-zero1})
and~(\ref{LTX-zero2}) lead to the same result.
\end{rem}
%
\begin{ex}
Case $N=2$. Suppose $x\ge 0$. Formula~(\ref{wrt.nu+1}) reads, with the
numerical values of Example~\ref{example1},
$$
\int_0^{\infty} e^{-\lambda t} \left[\mathbb{P}\{T(t) \in \mathrm{d}s, X(t)\in \mathrm{d}x\}
/(\mathrm{d}s\,\mathrm{d}x)\right] \mathrm{d}t
=e^{-\lambda s} \int_0^s \frac{1}{\sqrt \sigma}\,
E_{\frac12,\frac12}\!\left(\vphantom{\sqrt a}\right.\!\!-\sqrt{\lambda \sigma}
\!\left.\vphantom{\sqrt a}\right) I_{1,0}(s-\sigma;x)\,\mathrm{d}\sigma
$$
while formula~(\ref{wrt.nu+2}) gives, because of $\Phi_1=\Phi_2=I_{1,0}$ and~(\ref{kub}),
\begin{align*}
\lqn{\int_0^{\infty} e^{-\lambda t}\left[\mathbb{P}\{T(t) \in \mathrm{d}s, \,X(t)\in \mathrm{d}x\}
/(\mathrm{d}s \,\mathrm{d}x)\right] \mathrm{d}t}
&=
e^{-\lambda s} \left(\int_0^s \frac{1}{\sqrt{\sigma}}\,
E_{1,\frac12} (\lambda \sigma) \,\Phi_1(s-\sigma;x) \,\mathrm{d}\sigma -\sqrt \lambda
\int_0^s E_{1,1}(\lambda \sigma) \,\Phi_2(s-\sigma;x) \,\mathrm{d}\sigma\right)
\\
&=
e^{-\lambda s} \int_0^s \frac{1}{\sqrt{\sigma}} \left(E_{1,\frac12}(\lambda \sigma)
-\sqrt{\lambda \sigma}\,e^{\lambda \sigma}\right) I_{1,0} (s-\sigma;x)\,\mathrm{d}\sigma
\\
&=
e^{-\lambda s} \int_0^s \frac{1}{\sqrt{\sigma}} \,
E_{\frac12,\frac12}\!\left(\vphantom{\sqrt a}\right.\!\!-\sqrt{\lambda \sigma}
\!\left.\vphantom{\sqrt a}\right) I_{1,0} (s-\sigma;x)\,\mathrm{d}\sigma.
\end{align*}
We have checked that the two different representations~(\ref{wrt.nu+1})
and~(\ref{wrt.nu+2}) lead to the same result.
Let us pursue the computations. In view of~(\ref{case2I}) and Lemma~\ref{2half}, we get
\begin{align*}
\lqn{\int_0^{\infty} e^{-\lambda t}\left[\mathbb{P}\{T(t) \in \mathrm{d}s,
\, X(t)\in \mathrm{d}x\} /(\mathrm{d}s \,\mathrm{d}x)\right] \mathrm{d}t}
&=
\frac{x e^{-\lambda s}}{2 \sqrt{\pi}} \int_0^s \frac{e^{-\frac{x^2}{4(s-\sigma)}}}{\sqrt \sigma\,
(s-\sigma)^{3/2}} \,E_{\frac12,\frac12}\!\left(\vphantom{\sqrt a}\right.
\!\!-\sqrt{\lambda \sigma}\!\left.\vphantom{\sqrt a}\right) \mathrm{d}\sigma
\\
&=
\frac{x e^{-\lambda s}}{2 \sqrt{\pi}} \int_0^s \frac{e^{-\frac{x^2}{4\sigma}}}{
\sigma^{3/2}\sqrt{s-\sigma}} \,E_{\frac12,\frac12}\!\left(\vphantom{\sqrt a}\right.
\!\!-\sqrt{\lambda(s-\sigma)}\!\left.\vphantom{\sqrt a}\right) \mathrm{d}\sigma
\\
&=
\frac{xe^{-\lambda s}}{2 \sqrt \pi} \int_0^s \frac{e^{-\frac{x^2}{4\sigma}}}{
\sigma^{3/2}\sqrt{s-\sigma}} \left(\frac{1}{\sqrt \pi} -\sqrt{\lambda (s-\sigma)}
\,e^{\lambda(s-\sigma)}\,
\mathrm{Erfc}\!\left(\vphantom{\sqrt a}\right.\!\!\sqrt{\lambda (s-\sigma)}
\!\left.\vphantom{\sqrt a}\right)\right)\mathrm{d}\sigma
\\
&=
\frac{xe^{-\lambda s}}{2 \sqrt \pi} \left[\frac{1}{\sqrt \pi}
\int_0^s \frac{e^{-\frac{x^2}{4\sigma}}}{\sigma^{3/2}\sqrt{s-\sigma}} \,\mathrm{d}\sigma
-\frac{2 \sqrt{\lambda}}{\sqrt \pi} \int_0^s \frac{e^{\lambda(s-\sigma)-\frac{x^2}{4\sigma}}}{\sigma^{3/2}}
\left(\int_{\sqrt{\lambda (s-\sigma)}}^{\infty} e^{-\xi^2}\,\mathrm{d}\xi\right) \mathrm{d}\sigma\right]\!.
\end{align*}
The first integral in the last displayed equality writes, with
the change of variable $\sigma=s^2/\tau$,
$$
\int_0^s \frac{e^{-\frac{x^2}{4\sigma}}}{ \sigma^{3/2} \sqrt{s-\sigma}} \,\mathrm{d}\sigma
=\frac{1}{s^{3/2}}\int_{s}^{\infty} \frac{e^{-\frac{x^2}{4 s^2}\tau}}{\sqrt{\tau-s}} \,\mathrm{d}\tau
=\frac{e^{-\frac{x^2}{4 s}}}{s^{3/2}} \int_0^{\infty}
\frac{e^{-\frac{x^2}{4 s^2}\tau}}{\sqrt \tau} \,\mathrm{d}\tau
=\frac{2 \sqrt \pi}{\sqrt s \,x} \,e^{-\frac{x^2}{4 s}},
$$
and then
\begin{align*}
\lqn{\int_0^{\infty} e^{-\lambda t}\left[\mathbb{P}\{T(t) \in \mathrm{d}s, X(t)\in \mathrm{d}x\}
/(\mathrm{d}s \,\mathrm{d}x)\right]\mathrm{d}t}
&=
\frac{1}{\sqrt{\pi s}} \, e^{-\lambda s-\frac{x^2}{4s}}
-\frac{x \,\sqrt \lambda}{\pi} \int_0^s \frac{e^{-\lambda\sigma-\frac{x^2}{4\sigma}}}{\sigma^{3/2}}
\left(\int_{\sqrt{\lambda(s-\sigma)}}^{\infty} e^{-\xi^2} \,\mathrm{d}\xi\right)\mathrm{d}\sigma.
\end{align*}
The computation of the integral lying on the right-hand side of the
foregoing equality being cumbersome is postponed to Lemma~\ref{secondint}
in the appendix. The final result is, for $x\ge 0$,
$$
\int_0^{\infty} e^{-\lambda t}\left[\mathbb{P}\{T(t) \in \mathrm{d}s, X(t)\in \mathrm{d}x\}
/(\mathrm{d}s \,\mathrm{d}x)\right] \mathrm{d}t=\frac{e^{-\lambda s-\frac{x^2}{4s}}}{\sqrt{\pi s}}
-\sqrt \lambda \,e^{\sqrt{\lambda} \,x}\, \mathrm{Erfc}\!\left(\frac{x}{2 \sqrt s} +\sqrt{\lambda s}\right)\!.
$$
This is formula 1.4.6, page 129, of~\cite{bs}.
\end{ex}
%
\begin{ex}
Case $N=3$. For $x \ge 0$, (\ref{wrt.nu+1}) supplies with the numerical values of
Example~\ref{example2}, when $\kappa_3=-1$,
\begin{align*}
\lqn{\int_0^{\infty} e^{-\lambda t} \left[\mathbb{P}\{ T(t)\in \mathrm{d}s,
\, X(t) \in \mathrm{d}x\}/(\mathrm{d}s \,\mathrm{d}x)\right] \mathrm{d}t}
&=
\frac{e^{-\lambda s}}{\sqrt 3} \left(e^{\frac{i \pi}{6}} \int_0^s \sigma^{-2/3}
E_{\frac13,\frac13}\!\left(\vphantom{\sqrt a}\right.\!\!
-e^{-i\frac{\pi}{3}} \sqrt[3]{\lambda\sigma}\!\left.\vphantom{\sqrt a}\right)
I_{1,1}(s-\sigma;x) \,\mathrm{d}\sigma\right.
\\
&\hphantom{=\,}
\left. +\, e^{-\frac{i \pi}{6}} \int_0^s \sigma^{-2/3}
E_{\frac13,\frac13}\!\left(\vphantom{\sqrt a}\right.\!\!
-e^{i\frac{\pi}{3}} \sqrt[3]{\lambda\sigma}\!\left.\vphantom{\sqrt a}\right)
I_{2,1}(s-\sigma;x) \,\mathrm{d}\sigma\right)
\end{align*}
and when $\kappa_3=1$,
\begin{align*}
\lqn{\int_0^{\infty} e^{-\lambda t} \left[\mathbb{P}\{ T(t)\in \mathrm{d}s,
\, X(t) \in \mathrm{d}x\}/(\mathrm{d}s \,\mathrm{d}x)\right] \mathrm{d}t}
&=
\frac{i\,e^{-\lambda s}}{\sqrt 3 \,\sqrt[3]\lambda} \left(\int_0^s \sigma^{-2/3}
E_{\frac13,\frac13}\!\left(\vphantom{\sqrt a}\right.\!\!
e^{-i\frac{2\pi}{3}} \sqrt[3]{\lambda\sigma}\!\left.\vphantom{\sqrt a}\right)
I_{1,0}(s-\sigma;x) \,\mathrm{d}\sigma\right.
\\
&\hphantom{=\,}
\left. - \int_0^s \sigma^{-2/3} E_{\frac13,\frac13}\!\left(\vphantom{\sqrt a}\right.\!\!
e^{i\frac{2\pi}{3}} \sqrt[3]{\lambda\sigma}\!\left.\vphantom{\sqrt a}\right)
I_{1,0}(s-\sigma;x) \,\mathrm{d}\sigma\right)\!.
\end{align*}
The functions $I_{1,0}$, $I_{1,1}$ and $I_{2,1}$ above are respectively given
by~(\ref{case3I10}), (\ref{case3I11}) and~(\ref{case3I21}).
\end{ex}
%
\begin{ex}
Case $N=4$. For $x \ge 0$, (\ref{wrt.nu+1}) supplies, with the numerical values of
Example~\ref{example3},
\begin{align*}
\lqn{\int_0^{\infty} e^{-\lambda t} \left[\mathbb{P}\{ T(t)\in \mathrm{d}s,\, X(t) \in \mathrm{d}x\}
/(\mathrm{d}s \,\mathrm{d}x)\right] \mathrm{d}t}
&=
-\frac{e^{-\lambda s}}{2 \sqrt[4] \lambda}\left(e^{i\frac{\pi}{4}}
\int_0^s \sigma^{-3/4} E_{\frac14,\frac14}\!\left(\vphantom{\sqrt a}\right.\!\!
-\!\sqrt[4]{\lambda\sigma}\!\left.\vphantom{\sqrt a}\right)
I_{1,1}(s-\sigma;x)\,\mathrm{d}\sigma\right.
\\
&\hphantom{=\,}
\left. +\,e^{-i\frac{3\pi}{4}} \int_0^s \sigma^{-3/4}
E_{\frac14,\frac14}\!\left(\vphantom{\sqrt a}\right.\!\!
-i\sqrt[4]{\lambda\sigma}\!\left.\vphantom{\sqrt a}\right)
I_{1,1}(s-\sigma;x)\,\mathrm{d}\sigma\right.
\\
&\hphantom{=\,}
\left. +\,e^{i\frac{3\pi}{4}} \int_0^s \sigma^{-3/4}
E_{\frac14,\frac14}\!\left(\vphantom{\sqrt a}\right.\!\!
i\sqrt[4]{\lambda\sigma}\!\left.\vphantom{\sqrt a}\right)
I_{2,1}(s-\sigma;x)\,\mathrm{d}\sigma\right.
\\
&\hphantom{=\,}
\left.+\,e^{-i\frac{\pi}{4}} \int_0^s \sigma^{-3/4}
E_{\frac14,\frac14}\!\left(\vphantom{\sqrt a}\right.\!\!
-\!\sqrt[4]{\lambda\sigma}\!\left.\vphantom{\sqrt a}\right)
I_{2,1}(s-\sigma;x)\,\mathrm{d}\sigma\right)\!.
\end{align*}
The functions $I_{1,1}$ and $I_{2,1}$ above are given by~(\ref{case4I}).
\end{ex}
%

\section{Inverting with respect to $\lambda$}\label{section-inverting-lambda}

In this section, we perform the last inversion in $F(\lambda,\mu,\nu)$ in order
to derive the distribution of the couple $(T(t),X(t))$. As in the previous section,
we treat separately the two cases $x\le 0$ and $x\ge 0$.

\subsection{The case $x\le 0$}

%
\begin{teo}
The distribution of the couple $(T(t),X(t))$ is given, for $x \le 0$, by
\begin{align}
\lqn{\mathbb{P}\{T(t) \in \mathrm{d}s, X(t) \in \mathrm{d}x\}/ \mathrm{d}s \,\mathrm{d}x}
&=
-\frac{Ni}{2 \pi}\sum_{m=0}^{\# K} \alpha_{-m} s^{\frac{m-\#K}{N}}
\int_0^{\infty} \xi^{m+\# J} e^{-(t-s)\xi^N}
\mathcal{K}_m(x \xi) \,E_{1,\frac{m+\# J}{N}}(-s \xi^N)\,\mathrm{d}\xi
\label{wrt.lambda}
\end{align}
where
$$
\mathcal{K}_m(z)=e^{-i \frac{\# K-m-1}{N}\pi} \sum_{k \in K} B_k\theta_k^{m+1}
e^{-\theta_k e^{i\frac{\pi}{N}}z}-e^{i \frac{\# K -m-1}{N}\pi} \sum_{k \in K} B_k\theta_k^{m+1}
e^{-\theta_k e^{-i\frac{\pi}{N}}z}.
$$
\end{teo}
%
\Dim
Assume $x \le 0$. Recalling~(\ref{wrt.nu}), we have
\begin{align}
\lqn{\int_0^{\infty} e^{-\lambda t} \,[\mathbb{P}\{ T(t)\in \mathrm{d}s,
\, X(t) \in \mathrm{d}x\}/(\mathrm{d}s \,\mathrm{d}x)] \mathrm{d}t}
&=
-\frac{e^{-\lambda s}}{\lambda^{\frac{\# K -1}{N}} s^{\frac{\# K}{N}}}
\sum_{m=0}^{\#K}\alpha_{-m}(\lambda s)^{\frac{m}{N}} E_{1,\frac{m+\# J}{N}}(\lambda s)
\sum_{k \in K} B_k\theta_k^{m+1} e^{-\theta_k\!\!\sqrt[N]{\lambda}\,x}
\nonumber \\
&=
-\sqrt[N]{\lambda}\,e^{-\lambda s} \sum_{m=0}^{\#K}\alpha_{-m} \sum_{\ell=0}^{\infty}
\frac{(\lambda s)^{\ell+\frac{m-\#K}{N}}}{\Gamma
\!\left(\vphantom{\frac aN}\right.\!\!\ell+\frac{m+\# J}{N}\!\!\left.\vphantom{\frac aN}\right)}
\sum_{k \in K} B_k\theta_k^{m+1} e^{-\theta_k\!\!\sqrt[N]{\lambda}\,x}
\nonumber \\
&=
-\sum_{\ell=0}^{\infty}  \sum_{m=0}^{\#K} \alpha_{-m}
\frac{s^{\ell+\frac{m-\#K}{N}}}{\Gamma\!\left(\vphantom{\frac aN}\right.\!\!
\ell+\frac{m+\# J}{N}\!\!\left.\vphantom{\frac aN}\right)}
\sum_{k \in K} B_k\theta_k^{m+1} \,\lambda^{\ell+\frac{m-\# K+1}{N}}
e^{-\lambda s-\theta_k\!\!\sqrt[N]{\lambda}\,x}.
\label{wrt.lambda.inter}
\end{align}
We need to invert the quantity $\lambda^{\ell+\frac{m-\# K+1}{N}}
e^{-\lambda s-\theta_k\!\!\sqrt[N]{\lambda}\,x}$
for $\ell \ge 0$ and $0 \le m \le \#K$ with respect to $\lambda$.
We intend to use~(\ref{set24}) which is valid for $0\le m\le N-1$.
Actually~(\ref{set24}) holds true also for $m \le 0$; the proof of this claim
is postponed to Lemma~\ref{brom} in the appendix.
As a byproduct, for any $\ell \ge 0$ and $0 \le m \le \# K$,
\begin{align}
\lambda^{\ell+\frac{m-\# K+1}{N}} e^{-\lambda s-\theta_k\!\!\sqrt[N]\lambda\,x}
&=
e^{-\lambda s} \int_0^{\infty} e^{-\lambda u}I_{k,\#K-\ell N-m-1}(u;x)\,\mathrm{d}u
\nonumber \\
&=
\int_s^{\infty} e^{-\lambda t} I_{k,\# K-\ell N-m-1}(t-s;x)\,\mathrm{d}t.
\label{wrt.lambda.inter2}
\end{align}
Then, by putting~(\ref{wrt.lambda.inter2}) into~(\ref{wrt.lambda.inter}) and
next by eliminating the Laplace transform with respect to $\lambda$, we extract
\begin{align*}
\lqn{\mathbb{P}\{T(s)\in \mathrm{d}s,\,X(t)\in \mathrm{d}x\}/(\mathrm{d}s \,\mathrm{d}x)}
&=
-\sum_{\ell=0}^{\infty} \sum_{m=0}^{\# K} \alpha_{-m}
\frac{s^{\ell+\frac{m-\# K}{N}}}{\Gamma\!\left(\vphantom{\frac aN}\right.\!\!
\ell+\frac{m+\# J}{N}\!\!\left.\vphantom{\frac aN}\right)}
\sum_{k \in K} B_k\theta_k^{m+1} \,I_{k,\#K-\ell N-m-1}(t-s;x)
\\
&=
-\frac{Ni}{2 \pi} \sum_{\ell=0}^{\infty} \sum_{m=0}^{\# K} \alpha_{-m}
\frac{s^{\ell+\frac{m-\# K}{N}}}{\Gamma\!\left(\vphantom{\frac aN}\right.\!\!\ell+\frac{m+\# J}{N}
\!\!\left.\vphantom{\frac aN}\right)}
\\
&\hphantom{=\,}
\times\sum_{k \in K} B_k\theta_k^{m+1} \left(e^{-i \frac{\# K -\ell N-m-1}{N}\pi}
\int_0^{\infty} \xi^{N-\#K+\ell N+m} \,e^{-(t-s) \xi^N
-\theta_k e^{i \frac{\pi}{N}} x \xi} \,\mathrm{d}\xi\right.
\\
&\hphantom{=\,}
\left.-\,e^{i \frac{\#K-\ell N-m-1}{N}\pi} \int_0^{\infty} \xi^{N-\#K+\ell N+m}
\,e^{-(t-s) \xi^N-\theta_k e^{-i \frac{\pi}{N}} x \xi} \,\mathrm{d}\xi\right)
\\
&=
-\frac{Ni}{2 \pi} \sum_{m=0}^{\# K} \alpha_{-m} s^{\frac{m-\# K}{N}}
\sum_{k \in K} B_k\theta_k^{m+1}
\\
&\hphantom{=\,}
\times\left(e^{-i \frac{\# K -m-1}{N}\pi} \int_0^{\infty} \left(\sum_{\ell=0}^{\infty}
\frac{\left(-s \xi^N\right)^{\!\ell}}{\Gamma\!\left(\vphantom{\frac aN}\right.\!\!
\ell+\frac{m+\#J}{N}\!\!\left.\vphantom{\frac aN}\right)} \right)
\xi^{m+\#J} \,e^{-(t-s) \xi^N-\theta_k
e^{i \frac{\pi}{N}} x \xi} \,\mathrm{d}\xi\right.
\\
&\hphantom{=\,}
\left.-\,e^{i \frac{\#K-m-1}{N}\pi} \int_0^{\infty} \left(\sum_{\ell=0}^{\infty}
\frac{\left(-s \xi^N\right)^{\!\ell}}{\Gamma\!\left(\vphantom{\frac aN}\right.\!\!
\ell+\frac{m+\#J}{N}\!\!\left.\vphantom{\frac aN}\right)}\right)
\xi^{m+\#J} \,e^{-(t-s) \xi^N-\theta_k
e^{-i \frac{\pi}{N}} x \xi} \,\mathrm{d}\xi\right)
\\
\lqn{}
&=
-\frac{Ni}{2 \pi}  \sum_{m=0}^{\# K} \alpha_{-m} s^{\frac{m-\#K}{N}}
\sum_{k \in K} B_k\theta_k^{m+1}
\\
&\hphantom{=\,}
\times\left(e^{-i \frac{\# K -m-1}{N}\pi}
\int_0^{\infty}  \xi^{m+\# J} \,e^{-(t-s) \,\xi^N-\theta_k e^{i \frac{\pi}{N}} x \xi} \,
E_{1,\frac{m+\# J}{N}} \!\left(-s \xi^N\right)\mathrm{d}\xi\right.
\\
&\hphantom{=\,}
\left.-\,e^{i \frac{\#K-m-1}{N}\pi} \int_0^{\infty}
\xi^{m+\#J}\,e^{-(t-s) \,\xi^N-\theta_k e^{-i \frac{\pi}{N}} x \xi} \,
E_{1,\frac{m+\# J}{N}}\!\left(-s \xi^N\right) \mathrm{d}\xi\right)\!.
\end{align*}
The proof of~(\ref{wrt.lambda}) is established.
\EndDim

%
\begin{rem}\label{remark-integ}
Let us integrate~(\ref{wrt.lambda}) with respect to $x$ on $(-\infty,0]$.
We first compute, by using~(\ref{set26}),
\begin{align*}
\int_{-\infty}^0 \mathcal{K}_m (x \xi) \,\mathrm{d}x
&=
-\frac{1}{\xi} \left(\,\sum_{k \in K} B_k\theta_k^m\right)\!\!
\left(e^{-i \frac{\# K-m}{N}\pi}-e^{i\frac{\#K-m}{N}\pi}\right)
\\
&=
\frac{2i}{\xi} \sin\!\left(\frac{\#K-m}{N}\pi\right) \beta_m
=\begin{cases}
0 & \mbox{if } 1 \le m \le \#K ,\\[1ex]
\displaystyle{\frac{2i}{\xi} \,\sin\!\left(\frac{\#K}{N}\pi\right)} & \mbox{if } m=0.
\end{cases}
\end{align*}
We then obtain
\begin{align*}
\mathbb{P}\{T(t)\in \mathrm{d}s, X(t)\le 0\}/\mathrm{d}s
&=
\frac{N\sin\!\left(\!\!\vphantom{\frac aN}\right. \frac{\# K}{N}\pi\!\!
\left.\vphantom{\frac aN}\right)}{\pi s^{\frac{\# K}{N}}}
\int_0^{\infty} \xi^{\#J-1}\, e^{-(t-s) \xi^N} \,
E_{1,\frac{\# J}{N}}\!\left(-s \xi^N\right)\mathrm{d}\xi
\\
&=
\frac{N\sin\!\left(\!\!\vphantom{\frac aN}\right. \frac{\# K}{N}\pi\!\!
\left.\vphantom{\frac aN}\right)}{\pi s^{\frac{\# K}{N}}}
\sum_{\ell=0}^{\infty} \frac{(-s)^{\ell}}{\Gamma\!\left(\vphantom{\frac aN}\right.\!\!
\ell+\frac{\#J}{N}\!\!\left.\vphantom{\frac aN}\right)}
\int_0^{\infty}\xi^{\ell N+\#J-1} e^{-(t-s) \xi^N} \,\mathrm{d}\xi
\\
&=
\frac{\sin\!\left(\!\!\vphantom{\frac aN}\right. \frac{\# K}{N}\pi\!\!
\left.\vphantom{\frac aN}\right)}{\pi s^{\frac{\# K}{N}} (t-s)^{\frac{\#J}{N}}}
\sum_{\ell=0}^{\infty} \left(-\frac{s}{t-s}\right)^{\!\ell}
=\frac{\sin\!\left(\!\!\vphantom{\frac aN}\right. \frac{\# K}{N}\pi\!\!
\left.\vphantom{\frac aN}\right)}{\pi t} \left(\frac{t-s}{s}\right)^{\!\frac{\#K}{N}}.
\end{align*}
We retrieve~(\ref{f1}).
\end{rem}
%
\begin{rem}
Let us evaluate $\mathbb{P}\{T(t)\in \mathrm{d}s, X(t)\in \mathrm{d}x\}/(\mathrm{d}s\,\mathrm{d}x)$
at $x=0$. For $0 \le m \le \# K$,
$$
\mathcal{K}_m(0)=e^{-i \frac{\# K-m-1}{N}\pi}
\sum_{k \in K} B_k\theta_k^{m+1} -e^{i \frac{\# K-m-1}{N}\pi} \sum_{k \in K} B_k\theta_k^{m+1}
=-2i \sin\!\left(\frac{\# K-m-1}{N}\,\pi\right) \beta_{m+1}.
$$
Observing that $\sin\!\left(\!\!\vphantom{\frac aN}\right.
\frac{\# K-m-1}{N}\pi\!\!\left.\vphantom{\frac aN}\right)=0$ if $m=\# K-1$, in view
of~(\ref{set26}), (\ref{set14}) and~(\ref{set15}), we get
\begin{align*}
\lqn{\mathbb{P}\{T(t) \in \mathrm{d}s,\, X(t) \in \mathrm{d}x\}/ \mathrm{d}s \Big|_{x=0}}
&=
\frac{N}{\pi} \sin\!\left(\frac \pi N\right)\alpha_{-\#K} \beta_{\#K+1}
\int_0^{\infty} \xi^N \,e^{-(t-s) \xi^N} E_{1,1}\!\left(-s \xi^N\right) \mathrm{d}\xi
\\
&=
\frac{N}{\pi} \sin\!\left(\frac \pi N\right) \!\!\left(\!\vphantom{\sum_\in}\right.
\sum_{j \in J} \theta_j\!\!\left.\vphantom{\sum_in}\right) \int_0^{\infty} \xi^N \,e^{-t \xi^N} \mathrm{d}\xi
=\frac{\sin\!\left(\frac \pi N\right) \!\Gamma\!\left(\vphantom{\frac aN}\right.\!\!
\frac{1}{N}\!\!\left.\vphantom{\frac aN}\right)}{N \pi \,t^{1+\frac 1N}}\sum_{j \in J}\theta_j.
\end{align*}
Thanks to~(\ref{set7}) and~(\ref{p-at-zero}), we see that
$$
\mathbb{P}\{T(t) \in \mathrm{d}s,\, X(t) \in \mathrm{d}x\}/ \mathrm{d}s \Big|_{x=0}=\frac{1}{t} \,p(t;0)
$$
and we deduce
$$
\mathbb{P}\{ T(t) \in \mathrm{d}s| X(t)=0\}/\mathrm{d}s=\frac{\inde_{(0,t)}(s)}{t},
$$
that is, $(T(t)| X(t)=0)$ has the uniform law on $(0,t)$.
This is Theorem 2.13 of~\cite{2003}.
\end{rem}
%

\subsection{The case $x\ge 0$}

The case $x \ge 0$ can be related to the case $x \le 0$ by using the duality.
Let us introduce the dual process $(X^*_t)_{t \ge 0}$ of $(X_t)_{t \ge 0}$
defined as $X^*_t=-X_t$ for any $t\ge 0$. It is known that (see~\cite{2003}):
\begin{itemize}
\item
If $N$ is even, the processes $X$ and $X^*$ are identical in distribution
(because of the symmetry of the heat kernel $p$): $X^*\stackrel{d}{=}X$;
\item
If $N$ is odd, we have the equalities in distribution
$(X^+)^*\stackrel{d}{=}X^-$ and $(X^-)^*\stackrel{d}{=}X^+$
where $X^+$ is the pseudo-process associated with $\kappa_{\!_{ N}}=+1$ and
$X^-$ the one associated with $\kappa_{\!_{ N}}=-1$.
\end{itemize}
When $N$ is even, we have $\{-\theta_j, j\in J\}=\{\theta_k, k \in K\}$. In this case,
for any $j \in J$, there exists a unique $k \in K$ such that $\theta_j=-\theta_k$ and then
$$
A_j=\prod_{i \in J \setminus \{j\}} \frac{\theta_i}{\theta_i-\theta_j}
=\prod_{i \in K \setminus \{k\}} \frac{-\theta_i}{-\theta_i+\theta_k}
=\prod_{i \in K \setminus \{k\}} \frac{\theta_i}{\theta_i-\theta_k}=B_k
$$
and
$$
\alpha_{m}=\sum_{j \in J}  A_j \theta_j^m=\sum_{k \in K} B_k (-\theta_k)^m=(-1)^m \beta_m.
$$
When $N$ is odd, we distinguish the roots of $\kappa_{\!_{ N}}$ in the cases
$\kappa_{\!_{ N}}=+1$ and $\kappa_{\!_{ N}}=-1$:
\begin{itemize}
\item
For $\kappa_{\!_{ N}}=+1$, let $\theta^+_i$, $1 \le i \le N$, denote the roots of $1$
and set $J^+=\{i\in\{1,\dots,N\}: \Re(\theta^+_i) >0\}$ and
$K^+=\{i\in\{1,\dots,N\}: \Re(\theta_i^+)<0\}$;
\item
For $\kappa_{\!_{ N}}=-1$, let $\theta^-_i$, $1 \le i \le N$, denote the roots of $-1$
and set $J^-=\{i\in\{1,\dots,N\}: \Re(\theta^-_i) >0\}$ and
$K^-=\{i\in\{1,\dots,N\}: \Re(\theta_i^-)<0\}$.
\end{itemize}
We have $\{\theta^-_j, i \in J^-\}=\{-\theta_k^+, k \in K^+\}$ and
$\{\theta_k^-, k\in K^-\}=\{-\theta^+_j, j \in J^+\}$.
In this case, for any $j \in J^-$, there exists a unique $k \in K^+$ such that
$\theta_j^-=-\theta_k^+$ and then
$$
A_j^-=\prod_{i \in J^-\setminus \{j\}} \frac{\theta^-_i}{\theta^-_i-\theta^-_j}
=\prod_{i \in K^+ \setminus \{k\}} \frac{-\theta^+_i}{-\theta^+_i-\theta^+_k}
=\prod_{i \in K^+ \setminus \{k\}} \frac{\theta^+_i}{\theta^+_i-\theta^+_k}=B^+_k
$$
and similarly $A_j^+=B_k^-$. Moreover, we have
$$
\alpha_{m}^-=\sum_{j \in J^-}A_j^-(\theta_j^-)^m=\sum_{k \in K^+} B_k^+(-\theta_k^+)^m
=(-1)^m \sum_{k \in K^+} B_k^+ (\theta_k^+)^m=(-1)^m \beta_m^+
$$
and similarly $\alpha_m^+=(-1)^m \beta_m^-$.\\

Now, concerning the connection between sojourn time and duality, we have
the following fact. Set
$$
\tilde{T}(t)=\int_0^t \ind_{(0,+\infty)}(X(u))\,\mathrm{d}u
\quad \mbox{and} \quad T^*(t)=\int_0^t \ind_{[0,+\infty)}(X^*(u))\,\mathrm{d}u.
$$
Since Spitzer's identity holds true interchanging the closed interval
$[0,+\infty)$ and the open interval $(0,+\infty)$, it is easy to see
that $T(t)$ and $\tilde{T}(t)$ have the same distribution.
On the other hand, we have
$$
\tilde{T}(t)=\int_0^t \ind_{(0,+\infty)}(X(u))\,\mathrm{d}u
=\int_0^t \ind_{(-\infty,0)}(X^*(u))\,\mathrm{d}u
=\int_0^t [1-\ind_{[0,+\infty)}(X^*(u))]\,\mathrm{d}u=t-T^*(t).
$$
We then deduce that $T(t)$ and $t-T^*(t)$ have the same distribution. Consequently,
we can state the lemma below.
%
\begin{lem}\label{identity-duality}
The following identity holds:
$$
\mathbb{P}\{T(t)\in \mathrm{d}s, X(t)\in \mathrm{d}x\}/ (\mathrm{d}s \,\mathrm{d}x)
=\mathbb{P}\{T^*(t) \in \mathrm{d}(t-s), X^*(t)\in \mathrm{d}(-x)\}/(\mathrm{d}s \,\mathrm{d}x).
$$
\end{lem}
%

As a result, the following result ensues.
%
\begin{teo}\label{theorem-wrt.lambda+1}
Assume $N$ is even. The distribution of $(T(t),X(t))$ is given, for $x \ge 0$, by
\begin{align}
\lqn{\mathbb{P}\{T(t)\in \mathrm{d}s, X(t)\in \mathrm{d}x\}/ (\mathrm{d}s \,\mathrm{d}x)}
&=
\frac{Ni}{2 \pi} \sum_{m=0}^{\# J} \beta_{-m}\, (t-s)^{\frac{m-\#J}{N}}
\int_0^{\infty} \xi^{m+\# K} \,e^{-s \xi^N} \,\mathcal{J}_m(x \xi)\,
E_{1,\frac{m+\# K}{N}} \!\left(-(t-s)\xi^N\right) \mathrm{d}\xi
\label{wrt.lambda+1}
\end{align}
where
$$
\mathcal{J}_m(z)=e^{-i \frac{\# J-m-1}{N}\pi} \sum_{j \in J} A_j \theta_j^{m+1}
e^{-\theta_j e^{i\frac{\pi}{N}}z}-e^{i \frac{\# J-m-1}{N}\pi}
\sum_{j \in J} A_j \theta_j^{m+1} e^{-\theta_j e^{-i\frac{\pi}{N}}z}.
$$
\end{teo}
%
\Dim
When $N$ is even, we know that $X^*$ is identical in distribution to $X$
and $(T^*(t),X^*(t))$ is then distributed like $(T(t),X(t))$. Thus,
by~(\ref{wrt.lambda}) and Lemma~\ref{identity-duality}, for $x\ge 0$,
\begin{align*}
\lqn{\mathbb{P}\{T(t)\in \mathrm{d}s,\, X(t)\in \mathrm{d}x\}/ (\mathrm{d}s \,\mathrm{d}x)}
&=
\mathbb{P}\{T(t) \in \mathrm{d} (t-s),\, X(t) \in \mathrm{d}(-x)\}/ (\mathrm{d}s \,\mathrm{d}x)
\\
&=
-\frac{Ni}{2 \pi} \sum_{m=0}^{\# K} \alpha_{-m} (t-s)^{\frac{m-\#K}{N}}
\int_0^{\infty} \xi^{m+\# J} e^{-s \xi^N}\,\mathcal{K}_m(-x \xi) \,
E_{1,\frac{m+\# J}{N}} (-(t-s) \xi^N)\,\mathrm{d}\xi.
\end{align*}
The discussion preceding Lemma~\ref{identity-duality} shows that
$$
\mathcal{K}_m(z)
=e^{-i \frac{\# J-m-1}{N}\pi} \sum_{j \in J} A_j (-\theta_j)^{m+1} e^{\theta_j e^{i\frac{\pi}{N}} z}
-e^{i \frac{\# J -m-1}{N}\pi} \sum_{j \in J} A_j (-\theta_j)^{m+1} e^{\theta_j e^{-i\frac{\pi}{N}}z}.
$$
We see that $\mathcal{K}_m(z)=(-1)^{m+1}\mathcal{J}_m(-z)$ where the function
$\mathcal{J}_m$ is written in Theorem~\ref{theorem-wrt.lambda+1}.
Finally, by replacing $\alpha_{-m}$ by $(-1)^m \beta_{-m}$ and $\#J$, $\# K$ by
$\# K$, $\# J$ respectively (which actually coincide since $N$ is even),
(\ref{wrt.lambda+1}) ensues.
\EndDim

If $N$ is odd, although the results are not justified, similar formulas can
be stated. We find it interesting to produce them here.
We set $T^{\pm}(t)=\int_0^t \ind_{[0,+\infty)}(X^{\pm} (u))\,\mathrm{d}u$.
%
\begin{teo}\label{theorem-wrt.lambda+2}
Suppose that $N$ is odd. The distribution of $(T^+(t),X^+(t))$ is given, for $x \ge 0$, by
\begin{align}
\lqn{\mathbb{P}\{T^+(t)\in \mathrm{d}s, \, X^+(t)\in \mathrm{d}x\}/ (\mathrm{d}s \,\mathrm{d}x)}
&=
\frac{Ni}{2 \pi} \sum_{m=0}^{\# J^+} \beta^+_{-m}(t-s)^{\frac{m-\#J^+}{N}}
\int_0^{\infty} \xi^{m+\# K^+} \,e^{-s \xi^N} \,\mathcal{J}^+_m(x \xi)
\,E_{1,\frac{m+\# K^+}{N}} \!\left(-(t-s)\xi^N\right) \mathrm{d}\xi
\label{wrt.lambda+2}
\end{align}
where
$$
\mathcal{J}^+_m(z)=e^{-i \frac{\# J^+-m-1}{N}\pi} \sum_{j \in J^+} A_j^+(\theta_j^+)^{m+1}
e^{-\theta_j^+ e^{i\frac{\pi}{N}}z}-e^{i \frac{\# J^+-m-1}{N}\pi}
\sum_{j \in J^+} A_j^+ (\theta_j^+)^{m+1} e^{-\theta_j^+ e^{-i\frac{\pi}{N}}z}.
$$
\end{teo}
%
\Dim
When $N$ is odd, we know that $(X^+)^*\stackrel{d}{=}X^-$ and then
$((T^+)^*(t),(X^+)^*(t))\stackrel{d}{=}(T^-(t),X^-(t))$. Thus, by~(\ref{wrt.lambda})
and Lemma~\ref{identity-duality}, for $x\ge0$,
\begin{align*}
\lqn{\mathbb{P}\{T^+(t)\in \mathrm{d}s, \,X^+(t)\in \mathrm{d}x\}/ (\mathrm{d}s \,\mathrm{d}x)}
&=
\mathbb{P}\{T^-(t) \in \mathrm{d} (t-s), \,X^-(t) \in \mathrm{d}(-x)\}/ (\mathrm{d}s \,\mathrm{d}x)
\\
&=
-\frac{Ni}{2 \pi} \sum_{m=0}^{\# K^-} \alpha^-_{-m} (t-s)^{\frac{m-\#K^-}{N}}
\int_0^{\infty} \xi^{m+\# J^-} \,e^{-s \xi^N} \,\mathcal{K}^-_m(-x \xi)
\,E_{1,\frac{m+\# J^-}{N}} \!\left(-(t-s)\xi^N\right)\,\mathrm{d}\xi
\end{align*}
where
$$
\mathcal{K}^-_m(z)=e^{-i \frac{\# K^--m-1}{N}\pi} \sum_{k \in K^-} B^-_k(\theta_k^-)^{m+1}
e^{-{\theta}_k^-e^{i\frac{\pi}{N}}z}-e^{i \frac{\# K^--m-1}{N}\pi}
\sum_{k \in K^-} B^-_k(\theta^-_k)^{m+1} e^{-\theta_k^-e^{-i\frac{\pi}{N}}z}.
$$
As in the proof of Theorem~\ref{theorem-wrt.lambda+1}, we can write
$\mathcal{K}^-_m(z)=(-1)^{m+1}\mathcal{J}^+_m(-z)$ where the function
$\mathcal{J}^+_m$ is defined in Theorem~\ref{theorem-wrt.lambda+2}.
Finally, by replacing $\alpha^-_m$ by $(-1)^m \beta^+_m$ and
$\#J^-$, $\# K^-$ by $\# K^+$, $\# J^+$ respectively, (\ref{wrt.lambda+2}) ensues.
\EndDim
\pagebreak

Formula~(\ref{wrt.lambda+2}) involves only quantities with associated `$+$' signs.
We have a similar formula for $X^-$ by changing all `$+$' into `$-$'.
So, we can remove these signs
in order to get a unified formula (this is~(\ref{wrt.lambda+1}))
which is valid for even $N$ and, at least formally, for odd $N$ without sign.

%
\begin{rem}
Let us integrate~(\ref{wrt.lambda+1}) with respect to $x$ on $[0,\infty)$.
We first calculate, recalling that $\mathcal{J}_m(z)=(-1)^{m+1}\mathcal{K}_m(-z)$
and referring to Remark~\ref{remark-integ},
$$
\int_0^{\infty} \mathcal{J}_m (x \xi) \,\mathrm{d}x
=(-1)^{m+1}\int_{-\infty}^0 \mathcal{K}_m (x \xi) \,\mathrm{d}x
=\begin{cases}
\displaystyle{0} &  \mbox{if } 1 \le m \le \#J ,
\\[1ex]
\displaystyle{-\frac{2i}{\xi} \sin\!\left(\frac{\#J}{N}\pi\right)} &\mbox{if } m=0.
\end{cases}
$$
Then,
\begin{align*}
\mathbb{P}\{T(t)\in \mathrm{d}s, X(t)\ge 0\}/\mathrm{d}s
&=
\frac{N\sin\!\left(\!\!\vphantom{\frac aN}\right. \frac{\# J}{N}\pi\!\!
\left.\vphantom{\frac aN}\right)}{\pi(t-s)^{\frac{\# J}{N}}}
\int_0^{\infty} \xi^{\#K-1} \, e^{-s \xi^N}
\,E_{1,\frac{\# K}{N}}\!\left(-(t-s) \xi^N\right) \mathrm{d}\xi
\\
&=
\frac{N\sin\!\left(\!\!\vphantom{\frac aN}\right. \frac{\# J}{N}\pi\!\!
\left.\vphantom{\frac aN}\right)}{\pi(t-s)^{\frac{\# J}{N}}}
\sum_{\ell=0}^{\infty} \frac{(-(t-s))^{\ell}}{\Gamma\!\left(\vphantom{\frac aN}\right.\!\!
\ell+\frac{\#K}{N}\!\!\left.\vphantom{\frac aN}\right)}
\int_0^{\infty}\xi^{N\ell+\#K-1} e^{-s \xi^N} \,\mathrm{d}\xi
\\
&=
\frac{\sin\!\left(\!\!\vphantom{\frac aN}\right. \frac{\# J}{N}\pi\!\!
\left.\vphantom{\frac aN}\right)}{\pi s^{\frac{\# K}{N}} (t-s)^{\frac{\#J}{N}}}
\sum_{\ell=0}^{\infty} \left(-\frac{t-s}{s}\right)^{\!\ell}
=\frac{\sin\!\left(\!\!\vphantom{\frac aN}\right. \frac{\# J}{N}\pi\!\!
\left.\vphantom{\frac aN}\right)}{\pi t} \left(\frac{s}{t-s}\right)^{\!\frac{\#J}{N}}.
\end{align*}
We retrieve~(\ref{f2}).
\end{rem}
%

\subsection{Examples}

In this part, we write out the distribution of the couple $(T(t),X(t))$ in
the cases $N=2$, $N=3$ and $N=4$.
%
\begin{ex}
Case $N=2$. Using the numerical values of Example~\ref{example1},
formula~(\ref{wrt.lambda}) yields, for $x \le 0$,
\begin{align*}
\mathbb{P}\{T(t) \in \mathrm{d}s, X(t) \in \mathrm{d}x\}/(\mathrm{d}s \,\mathrm{d}x)
&=
\frac{1}{\pi} \left(\frac{1}{\sqrt s} \int_0^{\infty}
\xi \,e^{-(t-s) \xi^2} \,\tilde{\mathcal{K}}_0(x \xi)
\,E_{1,\frac12}\!\left(-s \xi^2\right) \mathrm{d}\xi\right.
\\
&\hphantom{=\,}
\left. +\int_0^{\infty} \xi^2 \,e^{-(t-s)\xi^2} \,\tilde{\mathcal{K}}_1(x \xi)
\,E_{1,1}\!\left(-s \xi^2\right) \mathrm{d}\xi\right)
\end{align*}
with
$$
E_{1,1} \!\left(-s \xi^2\right)\!=e^{-s \xi^2},
\quad\tilde{\mathcal{K}}_0(z)=-i\,\mathcal{K}_0(z)=i(e^{iz}-e^{-iz}),
\quad\tilde{\mathcal{K}}_1(z)=-i\,\mathcal{K}_1(z)=e^{iz}+e^{-iz}.
$$
On the other hand, we learn from Lemma~\ref{1alpha} that, for $\xi\ge 0$,
$$
E_{1,\frac12}(-s \xi^2)=\frac{1}{\sqrt \pi} \left(\!\!\vphantom{\int}\right.
1-2 \xi e^{-s \xi^2} \int_0^{s \xi}
e^{\frac{u^2}{s}} \,\mathrm{d}u \!\!\left.\vphantom{\int}\right)\!.
$$
Therefore,
\begin{align*}
\mathbb{P}\{T(t) \in \mathrm{d}s, X(t) \in \mathrm{d}x\}/(\mathrm{d}s \,\mathrm{d}x)
&=
\frac{i}{ \pi^{3/2} \sqrt s} \int_0^{\infty} \xi \left(
e^{i x \xi -(t-s) \xi^2}-e^{-ix \xi-(t-s) \xi^2}\right) \mathrm{d}\xi
\\
&\hphantom{=\,}
-\frac{2i}{\pi^{3/2} \sqrt s} \int_0^{\infty} \xi^2 \left(e^{i x \xi -t \xi^2}
-e^{-i x \xi-t \xi^2}\right) \!\!\left(\!\!\vphantom{\int}\right.
\int_0^{s \xi} e^{\frac{u^2}{s}} \,\mathrm{d}u\!\!\left.\vphantom{\int}\right) \mathrm{d}\xi
\\
&\hphantom{=\,}
+\frac{1}{\pi} \int_0^{\infty} \xi^2 \left(e^{i x \xi -t \xi^2} +e^{-i x \xi-t \xi^2}\right) \mathrm{d}\xi
\\
&=
\frac{i}{ \pi^{3/2} \sqrt s} \int_{-\infty}^{\infty} \xi \,e^{ix\xi-(t-s)\xi^2} \mathrm{d}\xi
+\frac{1}{\pi} \int_{-\infty}^{\infty} \xi^2 e^{ix\xi-t\xi^2} \mathrm{d}\xi
\\
&\hphantom{=\,}
+\frac{4}{\pi^{3/2} \sqrt s} \int_0^{\infty} \xi^2 \sin(x \xi) \,e^{-t\xi^2}
\left(\!\!\vphantom{\int}\right.\int_0^{s \xi} e^{\frac{u^2}{s}}
\,\mathrm{d}u\!\!\left.\vphantom{\int}\right) \mathrm{d}\xi.
\end{align*}
Let us compute two intermediate integrals:
\begin{align*}
\int_{-\infty}^{\infty} \xi \, e^{ix\xi-(t-s) \xi^2} \,\mathrm{d}\xi
&=
e^{-\frac{x^2}{4(t-s)}} \int_{-\infty}^{\infty} \xi \,e^{-(t-s)
\left(\xi-\frac{ix}{2(t-s)}\right)^{\!2}}\mathrm{d}\xi=e^{-\frac{x^2}{4(t-s)}}
\int_{-\infty}^{\infty}\left(\xi+\frac{ix}{2(t-s)}\right)e^{-(t-s) \xi^2} \,\mathrm{d}\xi
\\
&=
\frac{i \sqrt \pi \,x}{2 (t-s)^{3/2}} e^{-\frac{x^2}{4(t-s)}},
\end{align*}
\begin{align*}
\int_{-\infty}^{\infty} \xi^2  e^{i x \xi -t \xi^2}  \,\mathrm{d}\xi
&=
e^{-\frac{x^2}{4t}} \int_{-\infty}^{\infty}
\left(\xi + \frac{i x}{2t}\right)^{\!2} e^{-t \xi^2} \,\mathrm{d}\xi
=e^{-\frac{x^2}{4t}} \left(\int_{-\infty}^{\infty} \xi^2 e^{-t \xi^2} \,\mathrm{d}\xi
-\frac{x^2}{4 t^2}  \int_{-\infty}^{\infty} e^{-t \xi^2} \,\mathrm{d}\xi\right)
\\
&=
\frac{\sqrt \pi}{2 t^{3/2}} \, e^{-\frac{x^2}{4t}} \left(1-\frac{x^2}{2t}\right)\!.
\end{align*}
We deduce the following representation: for $x\le 0$,
\begin{align*}
\mathbb{P}\{T(t) \in \mathrm{d}s, X(t) \in \mathrm{d}x\}/(\mathrm{d}s \,\mathrm{d}x)
&=
\frac{4\sqrt s}{\pi^{3/2}} \left[\!\vphantom{\int}\right.
\int_0^{\infty} \xi^2 \sin(x \xi) \, e^{-t \xi^2}
\left(\!\!\vphantom{\int}\right.\int_0^{\xi} e^{su^2} \,\mathrm{d}u\!\!
\left.\vphantom{\int}\right)\mathrm{d}\xi -\frac{\sqrt \pi}{8}
\frac{xe^{-\frac{x^2}{4(t-s)}}}{s(t-s)^{3/2}}\!\!\left.\vphantom{\int}\right]
\\
&\hphantom{=\,}
+\frac{1}{2 \sqrt \pi t^{3/2}} \left(1-\frac{x^2}{2t}\right) e^{-\frac{x^2}{4t}}.
\end{align*}

For $x \ge 0$, (\ref{wrt.lambda+1}) gives
\begin{align*}
\mathbb{P}\{T(t) \in \mathrm{d}s, X(t) \in \mathrm{d}x\}/(\mathrm{d}s \,\mathrm{d}x)
&=
\frac{1}{\pi} \left(\frac{1}{\sqrt{t-s}} \int_0^{\infty}
\xi \,e^{-s\xi^2} \,\tilde{\mathcal{J}}_0(x \xi)
\,E_{1,\frac12}\!\left(-(t-s)\xi^2\right) \mathrm{d}\xi\right.
\\
&\hphantom{=\,}
\left. +\int_0^{\infty} \xi^2 \,e^{-s\xi^2} \,\tilde{\mathcal{J}}_1(x \xi)
\,E_{1,1}\!\left(-(t-s)\xi^2\right) \mathrm{d}\xi\right)
\end{align*}
with
$$
\tilde{\mathcal{J}}_0(z)=i\,\mathcal{J}_0(z)=-i(e^{iz}-e^{-iz}),
\quad\tilde{\mathcal{J}}_1(z)=-i\,\mathcal{J}_1(z)=e^{iz}+e^{-iz}.
$$
As previously,
\begin{align*}
\mathbb{P}\{T(t) \in \mathrm{d}s, X(t) \in \mathrm{d}x\}/(\mathrm{d}s \,\mathrm{d}x)
&=
-\frac{i}{ \pi^{3/2} \sqrt{t-s}} \int_{-\infty}^{\infty} \xi \,e^{ix\xi-s\xi^2} \mathrm{d}\xi
+\frac{1}{\pi} \int_{-\infty}^{\infty} \xi^2 e^{ix\xi-t\xi^2} \mathrm{d}\xi
\\
&\hphantom{=\,}
-\frac{4}{\pi^{3/2} \sqrt{t-s}} \int_0^{\infty} \xi^2 \sin(x \xi) \,e^{-t\xi^2}
\left(\!\!\vphantom{\int}\right.\int_0^{(t-s) \xi} e^{\frac{u^2}{t-s}}
\,\mathrm{d}u\!\!\left.\vphantom{\int}\right) \mathrm{d}\xi
\\
&=
\frac{4\sqrt{t-s}}{\pi^{3/2}} \left[
\frac{\sqrt \pi}{8}\frac{xe^{-\frac{x^2}{4s}}}{s^{3/2}(t-s)}-\!\!\left.\vphantom{\int}
\!\vphantom{\int}\right.\int_0^{\infty} \xi^2 \sin(x \xi) \,e^{-t \xi^2}
\left(\!\!\vphantom{\int}\right.\int_0^{\xi} e^{su^2} \,\mathrm{d}u\!\!
\left.\vphantom{\int}\right)\mathrm{d}\xi
\right]
\\
&\hphantom{=\,}
+\frac{1}{2 \sqrt \pi t^{3/2}} \left(1-\frac{x^2}{2t}\right) e^{-\frac{x^2}{4t}}.
\end{align*}

Actually, the density of $(T(t),X(t))$ related to (rescaled) Brownian motion is well-known under
another form. For instance, in~\cite{bs} pages 129--131, we find that
$$
\mathbb{P}\{T(t) \in \mathrm{d}s, X(t) \in \mathrm{d}x\}/(\mathrm{d}s \,\mathrm{d}x)=\begin{cases}
\displaystyle{\int_0^{\infty} \frac{y(y-x)}{4 \pi s^{3/2}(t-s)^{3/2}}\,
e^{-\frac{y^2}{4 s}-\frac{(y-x)^2}{4(t-s)}} \,\mathrm{d}y} &  \mbox{if } x\le 0,
\\[2ex]
\displaystyle{\int_0^{\infty} \frac{y(y+x)}{4 \pi s^{3/2}(t-s)^{3/2}}\,
e^{-\frac{(y+x)^2}{4 s}-\frac{y^2}{4(t-s)}} \,\mathrm{d}y} & \mbox{if } x\ge 0.
\end{cases}
$$
The coincidence of our representation and that of~\cite{bs} can be checked by using
Lemma~\ref{s} and Lemma~\ref{ss} in the appendix.
\end{ex}
%
\begin{ex}
Case $N=3$.

\noindent $\bullet$ Suppose $\kappa_3=1$. Using $E_{1,1}\!\left(-s \xi^3\right)=e^{-s \xi^3}$
and the values of Example~\ref{example2}, (\ref{wrt.lambda}) writes, for $x \le 0$,
\pagebreak
\begin{align*}
\lqn{\mathbb{P}\{T(t) \in \mathrm{d}s, X(t) \in \mathrm{d}x\}/(\mathrm{d}s \,\mathrm{d}x)}
&=
\frac{\sqrt 3}{2\pi} \left(s^{-2/3} \int_0^{\infty} \xi \,e^{-(t-s) \xi^3}\,
\tilde{\mathcal{K}}_0(x \xi)\,E_{1,\frac13} \!\left(-s \xi^3\right) \mathrm{d}\xi\right.
\\
&\hphantom{=\,}
\left. +\, s^{-1/3} \int_0^{\infty} \xi^2 \,e^{-(t-s) \xi^3}\,\tilde{\mathcal{K}}_1(x \xi)\,
E_{1,\frac23}\!\left(-s \xi^3\right) \mathrm{d}\xi
+ \int_0^{\infty} \xi^3\, e^{-t\xi^3}\,\tilde{\mathcal{K}}_2(x \xi) \,\mathrm{d}\xi\right)
\end{align*}
where
\begin{align*}
\tilde{\mathcal{K}}_0(z)
&=
-i\sqrt 3\,\mathcal{K}_0(z)= e^z-e^{-z/2}
\left(\!\!\vphantom{\frac{\sqrt a}{a}}\right.\cos\frac{\sqrt 3\,z}{2}
+\sqrt 3 \,\sin\frac{\sqrt 3\,z}{2}\left.\!\!\vphantom{\frac{\sqrt a}{a}}\right)\!,
\\
\tilde{\mathcal{K}}_1(z)
&=
-i\sqrt 3\,\mathcal{K}_1(z)= -e^z+e^{-z/2}
\left(\!\!\vphantom{\frac{\sqrt a}{a}}\right.
\cos\frac{\sqrt 3\,z}{2}-\sqrt 3 \,\sin\frac{\sqrt 3\,z}{2}
\left.\!\!\vphantom{\frac{\sqrt a}{a}}\right)\!,
\\
\tilde{\mathcal{K}}_2(z)
&=
-i\sqrt 3\,\mathcal{K}_2(z)= e^z+2e^{-z/2} \cos\frac{\sqrt 3\,z}{2}.
\end{align*}
For $x \ge 0$, (\ref{wrt.lambda+2}) gives
\begin{align*}
\lqn{\mathbb{P}\{T(t) \in \mathrm{d}s, X(t) \in \mathrm{d}x\}/(\mathrm{d}s \,\mathrm{d}x)}
&=
\frac{3}{2\pi} \left(\frac{1}{\!\sqrt[3]{t-s}} \int_0^{\infty} \xi^2 \,e^{-s \xi^3}\,
\tilde{\mathcal{J}}_0(x\xi)\,E_{1,\frac23}\!\left(-(t-s) \xi^3\right) \mathrm{d}\xi
+\int_0^{\infty} \xi^3\,e^{-t \xi^3}\,\tilde{\mathcal{J}}_1(x\xi) \,\mathrm{d}\xi\right)
\end{align*}
where
\begin{align*}
\tilde{\mathcal{J}}_0(z)
&=
i\,\mathcal{J}_0(z)=2 \,e^{-z/2} \,\sin\frac{\sqrt 3\,z}{2},
\\
\tilde{\mathcal{J}}_1(z)
&=
-i\,\mathcal{J}_1(z)= \,e^{-z/2} \left(\!\!\vphantom{\frac{\sqrt a}{a}}\right.
\sqrt 3\,\cos\frac{\sqrt 3\,z}{2}-\sin \frac{\sqrt 3\,z}{2}\left.\!\!\vphantom{\frac{\sqrt a}{a}}\right)\!.
\end{align*}

\noindent $\bullet$ Suppose $\kappa_3=-1$. Likewise, for $x \le 0$,
\begin{align*}
\lqn{\mathbb{P}\{T(t)\in \mathrm{d}s, X(t) \in \mathrm{d}x\}/(\mathrm{d}s \,\mathrm{d}x)}
&=
\frac{3}{2\pi} \left(\frac{1}{\!\sqrt[3]{s}} \int_0^{\infty} \xi^2\,e^{-(t-s) \xi^3}\,
\tilde{\mathcal{K}}_0(x\xi)\,E_{1,\frac2 3}\!\left(-s \xi^3\right) \mathrm{d}\xi\right.
+ \left. \int_0^{\infty} \xi^3\,e^{-t\xi^3} \tilde{\mathcal{K}}_1(x\xi) \,\mathrm{d}\xi\right)
\end{align*}
where
\begin{align*}
\tilde{\mathcal{K}}_0(z)
&=
-i\,\mathcal{K}_0(z)=-2\,e^{z/2} \,\sin \frac{\sqrt 3\,z}{2},
\\
\tilde{\mathcal{K}}_1(z)
&=
-i\,\mathcal{K}_1(z)=e^{z/2}
\left(\!\!\vphantom{\frac{\sqrt a}{a}}\right.\sqrt 3 \,\cos\frac{\sqrt 3\,z}{2}
+\sin\frac{\sqrt 3\,z}{2}\left.\!\!\vphantom{\frac{\sqrt a}{a}}\right)\!.
\end{align*}
For $x \ge 0$,
\begin{align*}
\lqn{\mathbb{P}\{T(t) \in \mathrm{d}s, X(t) \in \mathrm{d}x\}/(\mathrm{d}s \,\mathrm{d}x)}
\\[-3ex]
&=
\frac{\sqrt3}{2\pi} \left((t-s)^{-2/3}
\int_0^{\infty} \xi \,e^{-s \xi^3}\,\tilde{\mathcal{J}}_0(x \xi)\,
E_{1,\frac13}\!\left(-(t-s) \xi^3\right) \mathrm{d}\xi\right.
\\
&\hphantom{=\,}
\left.+\,(t-s)^{-1/3} \int_0^{\infty} \xi^2 \,e^{-s \xi^3}\,
\tilde{\mathcal{J}}_1(x \xi)\,E_{1,\frac23}\!\left(-(t-s) \xi^3\right) \mathrm{d}\xi
+\int_0^{\infty} \xi^3\,e^{-t\xi^3}\,\tilde{\mathcal{J}}_2(x \xi) \,\mathrm{d}\xi\right)
\end{align*}
where
\begin{align*}
\tilde{\mathcal{J}}_0(z)
&=
i\sqrt 3\,\mathcal{J}_0(z)=e^{-z}-e^{z/2}
\left(\!\!\vphantom{\frac{\sqrt a}{a}}\right.\cos\frac{\sqrt 3 \,z}{2}
-\sqrt 3 \,\sin\frac{\sqrt 3 \,z}{2}\left.\!\!\vphantom{\frac{\sqrt a}{a}}\right)\!,
\\
\tilde{\mathcal{J}}_1(z)
&=
-i\sqrt 3\,\mathcal{J}_1(z)=-e^{-z}+e^{z/2} \left(\!\!\vphantom{\frac{\sqrt a}{a}}\right.
\cos \frac{\sqrt 3 \,z}{2}+\sqrt 3 \,\sin\frac{\sqrt 3 \,z}{2}\left.\!\!\vphantom{\frac{\sqrt a}{a}}\right)\!,
\\
\tilde{\mathcal{J}}_2(z)
&=
i\sqrt 3\,\mathcal{J}_2(z)=e^{-z}+2\,e^{z/2} \cos\frac{\sqrt 3 \,z}{2}.
\end{align*}
\end{ex}
%
\begin{ex}
Case $N=4$. Referring to Example~\ref{example3}, formula~(\ref{wrt.lambda})
writes, for $x \le 0$,
\begin{align*}
\lqn{\mathbb{P}\{T(t) \in \mathrm{d}s, X(t) \in \mathrm{d}x\}/(\mathrm{d}s \,\mathrm{d}x)}
&=
\frac{2}{ \pi} \left(\frac{1}{\sqrt{s}} \int_0^{\infty} \xi^2 \,e^{-(t-s) \xi^4}\,
\tilde{\mathcal{K}}_0(x \xi) \,E_{1,\frac12}\!\left(-s \xi^4\right) \mathrm{d}\xi\right.
\\
&\hphantom{=\,}
\left. +\,\frac{\sqrt 2}{\sqrt[4]{s}} \int_0^{\infty} \xi^3\,e^{-(t-s) \xi^4}\,
\tilde{\mathcal{K}}_1(x \xi) \,E_{1,\frac34}\!\left(-s \xi^4\right) \mathrm{d}\xi
+\int_0^{\infty} \xi^4 \,e^{-t\xi^4}\, \tilde{\mathcal{K}}_2(x \xi) \,\mathrm{d}\xi\right)
\end{align*}
where
\begin{align*}
\tilde{\mathcal{K}}_0(z)
&=
-i\,\mathcal{K}_0(z)= e^z-\cos z-\sin z,
\\
\tilde{\mathcal{K}}_1(z)
&=
-i\,\mathcal{K}_1(z)= -e^z+\cos z-\sin z,
\\
\tilde{\mathcal{K}}_2(z)
&=
-i\,\mathcal{K}_2(z)= e^z+\cos z+\sin z.
\end{align*}
For $x \ge 0$, (\ref{wrt.lambda+1}) reads
\begin{align*}
\lqn{\mathbb{P}\{T(t) \in \mathrm{d}s, X(t) \in \mathrm{d}x\}/(\mathrm{d}s \,\mathrm{d}x)}
&=
\frac{2}{\pi} \left(\frac{1}{\sqrt{t-s}} \int_0^{\infty} \xi^2 \,e^{-s \xi^4}\,
\tilde{\mathcal{J}}_0(x \xi)\,E_{1,\frac12}\!\left(-(t-s) \xi^4\right) \mathrm{d}\xi\right.
\\
&\hphantom{=\,}
\left.+\,\frac{\sqrt 2}{\!\sqrt[4]{t-s}} \int_0^{\infty} \xi^3 \,e^{-s \xi^4}\,
\tilde{\mathcal{J}}_1(x \xi)\,E_{1,\frac34}\!\left(-(t-s) \xi^4\right) \mathrm{d}\xi
+ \int_0^{\infty} \xi^4 \,e^{-t\xi^4}\,\tilde{\mathcal{J}}_2(x \xi) \,\mathrm{d}\xi\right)
\end{align*}
where
\begin{align*}
\tilde{\mathcal{J}}_0(z)
&=
i\,\mathcal{J}_0(z)=e^{-z}-\cos z+\sin z,
\\
\tilde{\mathcal{J}}_1(z)
&=
-i\,\mathcal{J}_1(z)=-e^{-z}+\cos z+\sin z,
\\
\tilde{\mathcal{J}}_2(z)
&=
i\,\mathcal{J}_2(z)=e^{-z}+\cos z-\sin z.
\end{align*}
\end{ex}
%

\section{Appendix}\label{section-appendix}

%
\begin{lem}[Spitzer]\label{lemma-spitzer}
Let $(\xi_k)_{k \ge 1}$ be a sequence of independent identically distributed
random variables and set $X_0=0$ and $T_0=0$ and, for any $k\ge 1$,
$$
X_k=\xi_1+\dots+\xi_k,\qquad T_k=\sum_{j=1}^{k}\ind_{[0,+\infty)}(X_k).
$$
Then, for $\mu\in\mathbb{R}$, $\nu>0$ and $|z|<1$,
\begin{align}
\sum_{k=0}^{\infty}\mathbb{E}\!\left[e^{i \mu X_k-\nu T_k}\right] z^k
&=
\exp\!\left(\,\sum_{k=1}^{\infty} \mathbb{E}\!
\left[e^{i \mu X_k-\nu k \inde_{[0,+\infty)}(X_k)}\right]\frac{z^k}{k}\right)\!,
\label{spitzer-identity}\\
\sum_{k=0}^{\infty}\mathbb{E}\!\left[e^{i \mu X_k-\nu T_k}\ind_{[0,+\infty)}(X_k)\right] z^k
&=
\frac{1}{e^{\nu}-1}\left[e^{\nu}\!\!-
\exp\!\left(\!-\sum_{k=1}^{\infty} \left(1-e^{-\nu k}\right) \mathbb{E}\!
\left[e^{i \mu X_k}\ind_{[0,+\infty)}(X_k)\right]\!\frac{z^k}{k}\right)\right]\!,
\label{spitzer-identitybis}
\\
\sum_{k=0}^{\infty}\mathbb{E}\!\left[e^{i \mu X_k-\nu T_k}\ind_{(-\infty,0)}(X_k)\right] z^k
&=
\frac{e^{\nu}}{e^{\nu}-1}\left[
\exp\!\left(\,\sum_{k=1}^{\infty} \left(1-e^{-\nu k}\right) \mathbb{E}\!
\left[e^{i \mu X_k}\ind_{(-\infty,0)}(X_k)\right]\frac{z^k}{k}\right)-1\right]\!.
\label{spitzer-identityter}
\end{align}
\end{lem}
%
\Dim
Formula~(\ref{spitzer-identity}) is stated in~\cite{spitzer}
without proof. So, we produce a proof below which is rather similar to one lying
in~\cite{spitzer} related to the maximum functional of the $X_k$'s.
\\

\noindent $\bullet$ \textbf{Step 1.}
Set, for any $(x_1,\dots, x_n)\in \mathbb{R}^n$ and $\sigma \in \mathfrak{S}_n$
($\mathfrak{S}_n$ being the set of the permutations of $1,2,\dots,n$),
$$
U(x_1,\dots, x_n)=\sum_{k=1}^{n}\ind_{[0,\infty)} \!\left(\vphantom{\sum_{\in}}\right.\!
\sum_{j=1}^{k}x_j \!\!\left.\vphantom{\sum_{\in}}\right)
$$
and
$$
V(\sigma; x_1,\dots, x_n)=\sum_{k=1}^{n_{\sigma}}\# c_k(\sigma)
\ind_{[0,\infty)}\!\left(\vphantom{\sum_{\in}}\right.\!
\sum_{j \in c_k(\sigma)}x_j \!\!\left.\vphantom{\sum_{\in}}\right)\!.
$$
In the definition of $V$ above, the permutation $\sigma$ is decomposed into
$n_{\sigma}$ cycles: $\sigma=(c_1(\sigma))(c_2(\sigma))\dots (c_{n_\sigma}(\sigma))$.

In view of Theorem 2.3 in \cite{spitzer}, we have the equality between
the two following sets:
$$
\{U(\sigma(x_1),\dots,\sigma(x_n)),\sigma \in \mathfrak{S}_n\}
=\{V(\sigma; x_1,\dots, x_n),\sigma \in \mathfrak{S}_n\}.
$$
We then deduce, for any bounded Borel functions $\phi$ and $F$,
\begin{align*}
\mathbb{E}\!\left[\phi(X_n) F(U(\xi_1,\dots,\xi_n))\right]
&=
\frac{1}{n!} \sum_{\sigma \in \mathfrak{S}_n} \mathbb{E}
\!\left[\vphantom{\sum_{\in}}\right.\!\phi\!\left(\vphantom{\sum_{\in}}\right.\!\sum_{j=1}^n
\xi_{\sigma(j)}\!\!\left.\vphantom{\sum_{\in}}\right)\! F(V(\sigma; \xi_1,\dots,\xi_n))
\!\left.\vphantom{\sum_{\in}}\right]\!.
\end{align*}
In particular, for $\phi(x)=e^{i\mu x}$ and $F(x)=e^{-\nu x}$ (where
$\mu\in\mathbb{R}$ and $\nu>0$ are fixed),
\begin{align*}
\mathbb{E}\!\left[e^{i \mu X_n-\nu \sum_{k=1}^{n}
\inde_{[0,+\infty)}(\sum_{j=1}^k \xi_j)}\right]
&=
\frac{1}{n!} \sum_{\sigma \in \mathfrak{S}_n} \mathbb{E}
\!\left[\vphantom{\sum_{\in}}\right.\!\exp
\!\left(\vphantom{\sum_{\in}}\right.\!
i \mu \sum_{k=1}^{n_\sigma}\sum_{j \in c_k(\sigma)}\xi_j -\nu \sum_{k=1}^{n_\sigma} \# c_k(\sigma)
\ind_{[0,\infty)} \!\left(\vphantom{\sum_{\in}}\right.\!
\sum_{j \in c_k(\sigma)}\xi_j \!\left.\vphantom{\sum_{\in}}\right)\!\!
\!\left.\vphantom{\sum_{\in}}\right)\!\!\left.\vphantom{\sum_{\in}}\right]\!
\\
&=
\frac{1}{n!} \sum_{\sigma \in \mathfrak{S}_n} \prod_{k=1}^{n_\sigma}
\mathbb{E} \!\left[\vphantom{\sum_{\in}}\right.\!
\exp\!\left(\vphantom{\sum_{\in}}\right.\!
i \mu \sum_{j \in c_k(\sigma)}\xi_j -\nu \,(\# c_k(\sigma))
\ind_{[0,\infty)}
\!\left(\vphantom{\sum_{\in}}\right.\!
\sum_{j \in c_k(\sigma)}\xi_j
\!\left.\vphantom{\sum_{\in}}\right)\!\!\!\left.\vphantom{\sum_{\in}}\right)
\!\!\left.\vphantom{\sum_{\in}}\right]\!
\\
&=
\frac{1}{n!} \sum_{\sigma \in \mathfrak{S}_n} \prod_{k=1}^{n_\sigma}
\mathbb{E} \!\left[\vphantom{\sum_{\in}}\right.\!
\exp\!\left(\vphantom{\sum_{\in}}\right.\! i \mu \sum_{j=1}^{\# c_k(\sigma)}\xi_j -\nu \,(\# c_k(\sigma))
\ind_{[0,\infty)}
\!\left(\vphantom{\sum_{\in}}\right.\!
\sum_{j =1}^{\# c_k(\sigma)}\xi_j
\!\left.\vphantom{\sum_{\in}}\right)\!\!\!\left.\vphantom{\sum_{\in}}\right)\!\!\left.\vphantom{\sum_{\in}}\right]\!.
\end{align*}
Denote by $r_\ell(\sigma)$ the number of cycles of length $\ell$ in $\sigma$ for any
$\ell\in\{1,\dots,n\}$. We have $r_1(\sigma)+2r_2(\sigma)+\dots+nr_n(\sigma)=n$. Then,
\begin{align*}
\mathbb{E}\!\left[e^{i \mu X_n-\nu T_n}\right]
&=
\frac{1}{n!} \sum_{\sigma \in \mathfrak{S}_n} \prod_{\ell=1}^{n} \left(\mathbb{E}\!\left[e^{i \mu X_l-\nu\ell\,
\inde_{[0,\infty)}(X_\ell)}\right]\right)^{\!r_\ell(\sigma)}
\\
&=
\frac{1}{n!} \sum_{\substack{k_1,\dots, k_n\ge0: \\k_1+2k_2+\dots+nk_n=n}}
N_{k_1,\dots,k_n} \prod_{\ell=1}^{n} \left(\mathbb{E}\!\left[e^{i \mu X_\ell-\nu \ell\,
\inde_{[0,\infty)}(X_\ell)}\right]\right)^{\!k_\ell}
\end{align*}
where $N_{k_1,\dots,k_n}$ is the number of the permutations $\sigma$ of
$n$ objects satisfying $r_1(\sigma)=k_1,\dots,r_n(\sigma)=k_n$; this number
is equal to
$$
N_{k_1,\dots,k_n}=\frac{n!}{(k_1! 1^{k_1}) (k_2! 2^{k_2})\dots (k_n! n^{k_n})}.
$$
Then,
\begin{align*}
\mathbb{E}\!\left[e^{i \mu X_n-\nu T_n}\right]
&=
\sum_{\substack{k_1,\dots, k_n\ge0: \\ k_1+2k_2+\dots+nk_n=n}}
\prod_{\ell=1}^{n}\frac{1}{k_\ell! \ell^{k_\ell}}
\left(\mathbb{E} \left[e^{i \mu X_\ell-\nu \ell \,
\inde_{[0,\infty)}(X_\ell)}\right]\right)^{\!k_\ell}.
\end{align*}
\\

\noindent $\bullet$ \textbf{Step 2.}
Therefore, the identity between the generating functions follows: for $|z|<1$,
\begin{align*}
\sum_{n=0}^{\infty} \mathbb{E}\!\left[e^{i \mu X_n-\nu T_n}\right] z^n
&=
\sum_{\substack{n\ge 0,k_1,\dots, k_n\ge0: \\k_1+2k_2+\dots+nk_n=n}}
\prod_{\ell=1}^{n} \frac{1}{k_\ell!}\left(\mathbb{E}\!\left[e^{i \mu X_\ell-\nu \ell
\,\inde_{[0,\infty)}(X_\ell)}\right] \frac{z^{\ell}}{\ell}\right)^{\!k_\ell}
\\
&=
\sum_{k_1,k_2,\dots\ge 0} \prod_{\ell=1}^{\infty} \frac{1}{k_\ell!} \left(\mathbb{E}
\left[e^{i \mu X_\ell-\nu \ell \,\inde_{[0,\infty)}(X_\ell)}\right]
\frac{z^{\ell}}{\ell}\right)^{\!k_\ell}
\\
&=
\prod_{\ell=1}^{\infty} \Bigg[\sum_{k=1}^{\infty} \frac{1}{k!}
\left(\mathbb{E} \left[e^{i \mu X_\ell-\nu \ell \,\inde_{[0,\infty)}(X_\ell)}\right]
\frac{z^{\ell}}{\ell}\right)^{\!k}\Bigg]
\\
&=
\prod_{\ell=1}^{\infty} \exp\!\left(\mathbb{E}
\!\left[e^{i \mu X_\ell-\nu \ell \,\inde_{[0,\infty)}(X_\ell)}\right]
\frac{z^{\ell}}{\ell}\right)
\\
\lqn{}
&=
\exp\!\left(\sum_{n=1}^{\infty}
\mathbb{E}\!\left[e^{i \mu X_n-\nu n \inde_{[0,+\infty)}(X_n)}\right]
\frac{z^n}{n}\right)\!.
\end{align*}
The proof of~(\ref{spitzer-identity}) is finished.
\\

\noindent $\bullet$ \textbf{Step 3.}

Using the elementary identity $e^{a \inde_A(x)}-1=(e^a-1)\ind_A(x)$ and
noticing that $T_k=T_{k-1}+\ind_{[0,+\infty)}(X_k)$, we get for any $k\ge 1$,
$$
\mathbb{E}\!\left[e^{i \mu X_k-\nu T_k}\ind_{[0,+\infty)}(X_k)\right]
=\mathbb{E}\!\left[e^{i \mu X_k-\nu T_k}\,\frac{e^{\nu\inde_{[0,+\infty)}(X_k)}-1}{e^{\nu}-1}\right]
=\frac{1}{e^{\nu}-1}\left[\mathbb{E}\!\left(e^{i \mu X_k-\nu T_{k-1}}\right)\!
-\mathbb{E}\!\left(e^{i \mu X_k-\nu T_k}\right)\right]\!.
$$
Now, since $X_k=X_{k-1}+\xi_k$ where $X_{k-1}$ and $\xi_k$ are independent
and $\xi_k$ have the same distribution as $\xi_1$, we have, for $k\ge 1$,
$$
\mathbb{E}\!\left(e^{i\mu X_k-\nu T_{k-1}}\right)=
\mathbb{E}\!\left(e^{i\mu\xi_1}\right)
\mathbb{E}\!\left(e^{i \mu X_{k-1}-\nu T_{k-1}}\right)\!.
$$
Therefore,
\begin{align}
\sum_{k=1}^{\infty} \mathbb{E}\!\left[e^{i \mu X_k-\nu T_k}\ind_{[0,+\infty)}(X_k)\right]z^k
&=
\frac{1}{e^{\nu}-1}\sum_{k=1}^{\infty} \left(
\mathbb{E}\!\left[e^{i \mu X_k-\nu T_{k-1}}\right]\!
-\mathbb{E}\!\left[e^{i \mu X_k-\nu T_k}\right]\right)z^k
\nonumber\\
&=
\frac{1}{e^{\nu}-1}\left(\mathbb{E}\!\left(e^{i\mu\xi_1}\right)
\sum_{k=1}^{\infty} \mathbb{E}\!\left[e^{i \mu X_{k-1}-\nu T_{k-1}}\right]z^k
-\sum_{k=1}^{\infty} \mathbb{E}\!\left[e^{i \mu X_k-\nu T_k}\right]z^k\right)
\nonumber\\
&=
\frac{1}{e^{\nu}-1}\left(\left(z\,\mathbb{E}\!\left(e^{i\mu\xi_1}\right)-1\right)
\sum_{k=0}^{\infty} \mathbb{E}\!\left[e^{i \mu X_k-\nu T_k}\right]z^k+1\right)\!.
\label{spitzer-inter}
\end{align}
By putting~(\ref{spitzer-identity}) into~(\ref{spitzer-inter}),
we extract
\begin{equation}
\sum_{k=0}^{\infty} \mathbb{E}\!\left[e^{i \mu X_k-\nu T_k}\ind_{[0,+\infty)}(X_k)\right]z^k
=\frac{1}{e^{\nu}-1}\left[e^{\nu}-\left(1-z\,\mathbb{E}\!\left(e^{i\mu\xi_1}\right)\right)
S(\mu,\nu,z)\right]
\label{spitzerbis-inter}
\end{equation}
where we set
$$
S(\mu,\nu,z)=\exp\!\left(\,\sum_{k=1}^{\infty} \mathbb{E}\!
\left[e^{i \mu X_k-\nu k \inde_{[0,+\infty)}(X_k)}\right]\frac{z^k}{k}\right)\!.
$$
Next, using the elementary identity $1-\zeta=\exp[\log(1-\zeta)]
=\exp\!\left[-\sum_{k=1}^{\infty} \zeta^k/k\right]$ valid for $|\zeta|<1$,
$$
1-z\,\mathbb{E}\!\left(e^{i\mu\xi_1}\right)
=\exp\!\left(-\sum_{k=1}^{\infty} \left[\mathbb{E}\!\left(e^{i\mu\xi_1}\right)\right]^k \frac{z^k}{k}\right)
=\exp\!\left(-\sum_{k=1}^{\infty} \mathbb{E}\!\left(e^{i\mu X_k}\right) \frac{z^k}{k}\right)
$$
and then
\begin{align}
\left(1-z\,\mathbb{E}\!\left(e^{i\mu\xi_1}\right)\right) S(\mu,\nu,z)
&=
\exp\!\left(\,\sum_{k=1}^{\infty} \mathbb{E}\!
\left[e^{i \mu X_k-\nu k \inde_{[0,+\infty)}(X_k)}-e^{i \mu X_k}\right]\frac{z^k}{k}\right)
\nonumber\\
&=
\exp\!\left(-\sum_{k=1}^{\infty} \left(1-e^{-\nu k}\right) \mathbb{E}\!
\left[e^{i \mu X_k}\ind_{[0,+\infty)}(X_k)\right]\frac{z^k}{k}\right)\!.
\label{spitzerbis-interbis}
\end{align}
Hence, by putting~(\ref{spitzerbis-interbis}) into~(\ref{spitzerbis-inter}),
formula~(\ref{spitzer-identitybis}) entails.

By subtracting~(\ref{spitzerbis-inter}) from~(\ref{spitzer-identity}), we obtain the intermediate
representation
$$
\sum_{k=0}^{\infty} \mathbb{E}\!\left[e^{i \mu X_k-\nu T_k}\ind_{(-\infty,0)}(X_k)\right]z^k
=\frac{1}{e^{\nu}-1}\left[\left(e^{\nu}-z\,\mathbb{E}\!
\left(e^{i\mu\xi_1}\right)\right)S(\mu,\nu,z)-e^{\nu}\right]\!.
$$
By writing, as previously,
$$
e^{\nu}-z\,\mathbb{E}\!\left(e^{i\mu\xi_1}\right)
=e^{\nu}\exp\!\left(-\sum_{k=1}^{\infty} \mathbb{E}\!\left(e^{i\mu X_k}\right) \frac{e^{-\nu k}z^k}{k}\right),
$$
we find
\begin{align*}
\left(e^{\nu}-z\,\mathbb{E}\!\left(e^{i\mu\xi_1}\right)\right) S(\mu,\nu,z)
&=
e^{\nu}\exp\!\left(\,\sum_{k=1}^{\infty} \mathbb{E}\!
\left[e^{i \mu X_k-\nu k \inde_{[0,+\infty)}(X_k)}-e^{i \mu X_k-\nu k}\right]\frac{z^k}{k}\right)
\\
&=
e^{\nu}\exp\!\left(\,\sum_{k=1}^{\infty} \left(1-e^{-\nu k}\right) \mathbb{E}\!
\left[e^{i \mu X_k}\ind_{(-\infty,0)}(X_k)\right]\frac{z^k}{k}\right)\!.
\end{align*}
Finally, (\ref{spitzer-identityter}) ensues.
\EndDim

%
\begin{lem}\label{lemma-vdm}
The following identities hold:
$$
\beta_{\#K}=(-1)^{\#K -1}\prod_{k \in K} \theta_k,\qquad
\beta_{\#K+1}=(-1)^{\#K -1}\left(\prod_{k \in K} \theta_k\right)\!\!\left(\sum_{k \in K} \theta_k\right)\!.
$$
\end{lem}
%
\Dim
We label the set $K$ as $\{1,2,3,\dots,\#K\}$.
By~(\ref{set13}), we know that the $B_k$'s solve a Vandermonde system. Then, by Cramer's formulas,
we can write them as fractions of some determinants: $B_k=V_k/V$ where
$$
V=\begin{vmatrix}
1                & \dots & 1                   \\
\theta_1         & \dots & \theta_{\#K}        \\
\theta_1^2       & \dots & \theta_{\#K}^2      \\
\vdots           &       & \vdots              \\[.5ex]
\theta_1^{\#K-1} & \dots & \theta_{\#K}^{\#K-1}
\end{vmatrix}
\quad\mbox{and}\quad
V_k=\begin{vmatrix}
1                & \dots & 1                    & 1      & 1                    & \dots & 1                   \\
\theta_1         & \dots & \theta_{k-1}         & 0      & \theta_{k+1}         & \dots & \theta_{\#K}        \\
\theta_1^2       & \dots & \theta_{k-1}^2       & 0      & \theta_{k+1}^2       & \dots & \theta_{\#K}^2      \\
\vdots           &       & \vdots               & \vdots & \vdots               &       & \vdots              \\[.5ex]
\theta_1^{\#K-1} & \dots & \theta_{k-1}^{\#K-1} & 0      & \theta_{k+1}^{\#K-1} & \dots & \theta_{\#K}^{\#K-1}
\end{vmatrix}.
$$
By expanding the determinant $V_k$ with respect to its $k^{\mathrm{th}}$ column and next
factorizing it suitably, we easily see that
\begin{align*}
V_k
&=
(-1)^{k+1}\begin{vmatrix}
\theta_1         & \dots & \theta_{k-1}         & \theta_{k+1}         & \dots & \theta_{\#K}        \\
\theta_1^2       & \dots & \theta_{k-1}^2       & \theta_{k+1}^2       & \dots & \theta_{\#K}^2      \\
\vdots           &       & \vdots               & \vdots               &       & \vdots              \\[.5ex]
\theta_1^{\#K-1} & \dots & \theta_{k-1}^{\#K-1} & \theta_{k+1}^{\#K-1} & \dots & \theta_{\#K}^{\#K-1}
\end{vmatrix}
\\
&=
(-1)^{k+1}\frac{\prod_{i\in K} \theta_i}{\theta_k}
\begin{vmatrix}
1                & \dots & 1                    & 1                    & \dots & 1                   \\
\theta_1         & \dots & \theta_{k-1}         & \theta_{k+1}         & \dots & \theta_{\#K}        \\
\theta_1^2       & \dots & \theta_{k-1}^2       & \theta_{k+1}^2       & \dots & \theta_{\#K}^2      \\
\vdots           &       & \vdots               & \vdots               &       & \vdots              \\[.5ex]
\theta_1^{\#K-2} & \dots & \theta_{k-1}^{\#K-2} & \theta_{k+1}^{\#K-2} & \dots & \theta_{\#K}^{\#K-2}
\end{vmatrix}.
\end{align*}
With this at hands, we have
$$
\beta_{\#K}=\sum_{k\in K} B_k\theta_k^{\#K}
=\frac{\prod_{k\in K} \theta_k}{V} \sum_{k\in K}(-1)^{k+1}\theta_k^{\#K-1}
\begin{vmatrix}
1                & \dots & 1                    & 1                    & \dots & 1                   \\
\theta_1         & \dots & \theta_{k-1}         & \theta_{k+1}         & \dots & \theta_{\#K}        \\
\theta_1^2       & \dots & \theta_{k-1}^2       & \theta_{k+1}^2       & \dots & \theta_{\#K}^2      \\
\vdots           &       & \vdots               & \vdots               &       & \vdots              \\[.5ex]
\theta_1^{\#K-2} & \dots & \theta_{k-1}^{\#K-2} & \theta_{k+1}^{\#K-2} & \dots & \theta_{\#K}^{\#K-2}
\end{vmatrix}.
$$
We can observe that the sum lying on the above right-hand side is nothing but
the expansion of the determinant $V$ with respect to its last row multiplied
by the sign $(-1)^{\#K -1}$. This immediately ensues that
$\beta_{\#K}=(-1)^{\#K -1}\prod_{k \in K} \theta_k$. Similarly,
$$
\beta_{\#K+1}=\sum_{k\in K} B_k\theta_k^{\#K+1}
=\frac{\prod_{k\in K} \theta_k}{V} \sum_{k\in K}(-1)^{k+1}\theta_k^{\#K}
\begin{vmatrix}
1                & \dots & 1                    & 1                    & \dots & 1                   \\
\theta_1         & \dots & \theta_{k-1}         & \theta_{k+1}         & \dots & \theta_{\#K}        \\
\theta_1^2       & \dots & \theta_{k-1}^2       & \theta_{k+1}^2       & \dots & \theta_{\#K}^2      \\
\vdots           &       & \vdots               & \vdots               &       & \vdots              \\[.5ex]
\theta_1^{\#K-2} & \dots & \theta_{k-1}^{\#K-2} & \theta_{k+1}^{\#K-2} & \dots & \theta_{\#K}^{\#K-2}
\end{vmatrix}.
$$
The above sum is the expansion with respect to its last row, multiplied by the sign $(-1)^{\#K -1}$,
of the determinant $V'$ defined as
$$
V'=\begin{vmatrix}
1                & \dots & 1                   \\
\theta_1         & \dots & \theta_{\#K}        \\
\theta_1^2       & \dots & \theta_{\#K}^2      \\
\vdots           &       & \vdots              \\[.5ex]
\theta_1^{\#K-2} & \dots & \theta_{\#K}^{\#K-2}\\
\theta_1^{\#K}   & \dots & \theta_{\#K}^{\#K}
\end{vmatrix}.
$$
Let $R_0,R_1,R_2,\dots,R_{\#K-2},R_{\#K-1}$ denote the rows of $V'$.
We perform the substitution $R_{\#K-1}\leftarrow R_{\#K-1}+\sum_{\ell=2}^{\#K}
(-1)^{\ell}\sigma_{\ell} R_{\#K-\ell}$ where the $\sigma_{\ell}$'s are defined
by~(\ref{set10}). This substitution does not affect the value of $V'$
and it transforms, e.g., the first term of the last row into
\begin{align*}
\theta_1^{\#K}+\sum_{\ell=2}^{\#K} (-1)^{\ell}\sigma_{\ell} \,\theta_1^{\#K-\ell}.
\end{align*}
Recall that $\sigma_{\ell}=\sum_{1\le k_1<\dots<k_\ell\le\#K}
\theta_{k_1} \dots \theta_{k_\ell}$.
We decompose $\sigma_{\ell}$, by isolating the terms involving $\theta_1$, into
$$
\theta_1 \sum_{2\le k_2<\dots<k_\ell\le \#K} \theta_{k_2} \dots \theta_{k_\ell}
+\sum_{2\le k_1<k_2<\dots<k_\ell\le \#K} \theta_{k_1}\theta_{k_2} \dots \theta_{k_\ell}
=\theta_1\,\sigma'_{\ell-1}+\sigma'_{\ell}
$$
where we set $\sigma'_{\#K}=0$ and $\sigma'_{\ell}=\sum_{2\le k_1<k_2<\dots<k_\ell\le \#K}
\theta_{k_1}\theta_{k_2} \dots \theta_{k_\ell}.$ Therefore, we have
\begin{align*}
\theta_1^{\#K}+\sum_{\ell=2}^{\#K} (-1)^{\ell}\sigma_{\ell} \,\theta_1^{\#K-\ell}
&=
\theta_1^{\#K}+\sum_{\ell=2}^{\#K} (-1)^{\ell}\sigma'_{\ell-1} \,\theta_1^{\#K-\ell+1}
+\sum_{\ell=2}^{\#K} (-1)^{\ell}\sigma'_{\ell} \,\theta_1^{\#K-\ell}
\\
&=
\theta_1^{\#K}+\sigma'_1\,\theta_1^{\#K-1}=\theta_1^{\#K-1}(\theta_1+\sigma'_1)
=\theta_1^{\#K-1}\left(\sum_{k \in K} \theta_k\right)\!.
\end{align*}
The foregoing manipulation works similarly for each term of the last row
of $V'$. So, we deduce that $V'=\left(\sum_{k \in K} \theta_k\right)V$
and finally $\beta_{\#K+1}=(-1)^{\#K -1}\left(\prod_{k \in K} \theta_k\right)
\!\!\left(\sum_{k \in K} \theta_k\right)$.
\EndDim

%
\begin{lem}\label{1alpha}
For $\alpha>0$, the Mittag-Leffler functions $E_{1,\alpha}$ and $E_{1,\alpha+1}$
admit the following integral representations:
\begin{align}
E_{1,\alpha}(x)
&=
\begin{cases}
\displaystyle{\frac{1}{\Gamma(\alpha)} \left(1+ x^{1-\alpha} e^{x}
\int_0^x u^{\alpha-1} e^{-u} \,\mathrm{d}u\right)} & \mbox{if } x>0,
\\[3ex]
\displaystyle{\frac{1}{\Gamma(\alpha)} \left(1-|x|^{1-\alpha} e^{x}
\int_0^{|x|} u^{\alpha-1} e^u \,\mathrm{d}u\right)} & \mbox{if } x<0,
\end{cases}
\label{mittag1}\\[1ex]
E_{1,\alpha+1}(x)
&=
\begin{cases}
\displaystyle \frac{e^x}{\Gamma(\alpha) \,x^{\alpha}}
\int_0^x u^{\alpha-1} e^{-u} \,\mathrm{d}u & \mbox{if } x>0,
\\[2ex]
\displaystyle \frac{e^x}{\Gamma(\alpha) \,|x|^{\alpha}}
\int_0^{|x|} u^{\alpha-1} e^u \,\mathrm{d}u & \mbox{if } x<0.
\end{cases}
\label{mittag2}
\end{align}
\end{lem}
%
\Dim
Using the series expansion of $E_{1,\alpha}$, we obtain
\begin{align*}
x E_{1,\alpha}'(x)
&=
\sum_{r=1}^{\infty} \frac{r x^r}{\Gamma(r+\alpha)}
=\sum_{r=1}^{\infty} \frac{(r +\alpha-1)x^r}{\Gamma(r+\alpha)}+
(1-\alpha)\sum_{r=1}^{\infty} \frac{x^r}{\Gamma(r+\alpha)}
\\
&=
x \sum_{r=1}^{\infty} \frac{x^{r-1}}{\Gamma(r+\alpha-1)}
+(1-\alpha) \left(E_{1,\alpha}(x)-\frac{1}{\Gamma(\alpha)}\right)
=\left(x+1-\alpha\right)E_{1,\alpha}(x)+\frac{\alpha-1}{\Gamma(\alpha)}.
\end{align*}
Hence, $E_{1,\alpha}$ solves the differential equation
$x E_{1,\alpha}'(x)=(x+1-\alpha)E_{1,\alpha}(x)+\frac{\alpha-1}{\Gamma(\alpha)}$.
In view of this equation, we know that $E_{1,\alpha}(x)$ is of the form
$$
E_{1,\alpha}(x)=\lambda(x)\,e^{\int \frac{x+1-\alpha}{x}\,
\mathrm{d}x}=\lambda(x) \, |x|^{1-\alpha} \, e^x
$$
where the unknown function $\lambda$ solves
$$
\lambda'(x)=\frac{\alpha-1}{\Gamma(\alpha)x} |x|^{\alpha-1} e^{-x}
=\begin{cases}
\displaystyle{\frac{\alpha-1}{\Gamma(\alpha)} \,x^{\alpha-2}} \, e^{-x} &\mbox{if } x>0,
\\[2ex]
\displaystyle{-\frac{\alpha-1}{\Gamma(\alpha)} \,|x|^{\alpha-2}}\, e^{-x} &\mbox{if } x<0.
\end{cases}
$$
This implies that, for a certain $x_0>0$ and $\lambda_0 \in \mathbb{R}$
($\lambda_0$ could be different for $x>0$ and $x <0$), we have,
for $x >0$,
$$
\lambda(x)=\frac{\alpha-1}{\Gamma(\alpha)}
\int_{x_0}^{x} u^{\alpha-2}\,e^{-u}\,\mathrm{d}u+\lambda_0
=\frac{1}{\Gamma(\alpha)}\left(x^{\alpha-1} \,e^{-x}-x_0^{\alpha-1} \,e^{-x_0}
+\int_{x_0}^{x} u^{\alpha-1}\,e^{-u} \,\mathrm{d}u\right)+\lambda_0,
$$
and, for $x<0$,
$$
\lambda(x)=\frac{\alpha-1}{\Gamma(\alpha)}
\int_{x_0}^{|x|} u^{\alpha-2}e^u\,\mathrm{d}u +\lambda_0
=\frac{1}{\Gamma(\alpha)}\left(|x|^{\alpha-1}\,e^{|x|}
-x_0^{\alpha-1}\,e^{x_0} - \int_{x_0}^{|x|}
u^{\alpha-1}\,e^{u} \,\mathrm{d}u\right)+\lambda_0.
$$
Then
\begin{equation}
E_{1,\alpha}(x)=\begin{cases}
\displaystyle{\frac{1}{\Gamma(\alpha)}+\left(\lambda_0
-\frac{e^{-x_0}x_0^{\alpha-1}}{\Gamma(\alpha)}\right)
x^{1-\alpha} e^{x}+\frac{x^{1-\alpha}e^{x}}{\Gamma(\alpha)}
\int_{x_0}^{x} u^{\alpha-1} e^{-u} \,\mathrm{d}u} & \mbox{if } x>0,
\\[3ex]
\displaystyle{\frac{1}{\Gamma(\alpha)}+\left(\lambda_0
-\frac{e^{x_0}x_0^{\alpha-1}}{\Gamma(\alpha)}\right)
|x|^{1-\alpha} e^{x}-\frac{|x|^{1-\alpha}e^{x}}{\Gamma(\alpha)}
\int_{x_0}^{|x|} u^{\alpha-1} e^u \,\mathrm{d}u} & \mbox{if } x<0.
\end{cases}
\label{mittag-inter}
\end{equation}
Because of the singularity of the differential equation at zero,
the initial value at zero does not determine $x_0$ and $\lambda_0$.
Nevertheless, we know that $E_{1,\alpha}$ is $C^1$ at zero. So, we need to
compute $E_{1,\alpha}'(x)$ for $x\neq 0$:
$$
E_{1,\alpha}'(x)=\begin{cases}
\displaystyle{\left(\lambda_0-\frac{x_0^{\alpha-1}e^{-x_0}}{\Gamma(\alpha)}
+ \frac{1}{\Gamma(\alpha)} \int_{x_0}^{x} u^{\alpha-1} e^{-u} \,\mathrm{d}u\right)\!
\left(x^{1-\alpha} e^{x}+(1-\alpha)x^{-\alpha} e^{x}\right)+\frac 1{\Gamma(\alpha)}}
&  \mbox{if } x>0,
\\[3ex]
\displaystyle{\left(\vphantom{\int_{x_0}^x}\right.\!\!
\lambda_0-\frac{x_0^{\alpha-1}e^{x_0}}{\Gamma(\alpha)}
-\frac{1}{\Gamma(\alpha)} \int_{x_0}^{|x|} u^{\alpha-1} e^u \,\mathrm{d}u
\!\!\left.\vphantom{\int_{x_0}^x}\right)\!
\left(|x|^{1-\alpha} e^{x}-(1-\alpha)|x|^{-\alpha} e^{x}\right)+\frac 1{\Gamma(\alpha)}}
& \mbox{if } x<0.
\end{cases}
$$
In order to have a $C^1$-function at $0$, we must have
\begin{equation}
\lambda_0-\frac{x_0^{\alpha-1}e^{-x_0}}{\Gamma(\alpha)}
= \frac{1}{\Gamma(\alpha)} \int_0^{x_0} u^{\alpha-1} e^{-u} \,\mathrm{d}u
\quad\mbox{ and }\quad\lambda_0-\frac{x_0^{\alpha-1}e^{x_0}}{\Gamma(\alpha)}
=-\frac{1}{\Gamma(\alpha)} \int_0^{x_0} u^{\alpha-1} e^u \,\mathrm{d}u.
\label{conditions}
\end{equation}
Putting~(\ref{conditions}) into~(\ref{mittag-inter}) yields~(\ref{mittag1}).
Next, formula~(\ref{mittag2}) can be deduced from~(\ref{mittag1}) by simply observing that,
e.g., for $x>0$,
$$
E_{1,\alpha+1}(x)=\frac{x^{-\alpha}e^x}{\Gamma(\alpha+1)} \left(x^{\alpha}e^{-x}+
\int_0^x u^{\alpha} e^{-u} \mathrm{d}u\right)=\frac{x^{-\alpha} e^x}{\Gamma(\alpha+1)}
\left(\alpha \int_0^x u^{\alpha-1} e^{-u} \,\mathrm{d}u\right)\!.
$$
\EndDim

%
\begin{lem} \label{8.2}
For $\alpha \in (0,2)$, the function $E_{1,\alpha}$ admits the following representation. For $x > 0$,
$$
E_{1,\alpha}(x)=x^{1-\alpha} \left(e^x +\frac{\sin (\alpha \pi)}{\pi}
\int_0^{\infty} \frac{\xi^{1-\alpha}}{\xi+1} \,e^{-x \xi} \,\mathrm{d}\xi\right)\!.
$$
\end{lem}
%
\Dim
Writing $\frac{e^{-\xi x}}{\xi+1}=e^x \int_x^{\infty} e^{-( \xi +1)u}\,\mathrm{d}u$, we obtain
\begin{align*}
\int_0^{\infty} \frac{\xi^{1-\alpha}}{\xi+1} \,e^{-x \xi} \,\mathrm{d}\xi
&=
e^{x} \int_0^{\infty} \xi^{1-\alpha} \,\mathrm{d}\xi \int_x^{\infty} e^{-(\xi+1)u} \,\mathrm{d}u
=\Gamma(2-\alpha) e^x \int_x^{\infty} u^{\alpha-2}e^{-u}\,\mathrm{d}u
\\
&=
\frac{\Gamma(2-\alpha)}{1-\alpha}\,e^x\!\left(x^{\alpha-1} e^{-x}-
\int_{x}^{\infty} u^{\alpha-1} e^{-u}\,\mathrm{d}u\right)
\\
&=
\Gamma(1-\alpha)\,e^x \!\left(x^{\alpha-1} e^{-x}-\Gamma(\alpha)
+\int_0^x u^{\alpha-1} e^{-u}\,\mathrm{d}u\right)\!.
\end{align*}
From this, it entails that
$$
\int_0^x u^{\alpha-1} e^{-u}\,\mathrm{d}u=\frac{e^{-x}}{\Gamma(1-\alpha)}
\int_0^{\infty} \frac{\xi^{1-\alpha}}{\xi+1} \,e^{-x \xi} \,\mathrm{d}\xi
+\Gamma(\alpha)-x^{\alpha-1} e^{-x}
$$
which, by~(\ref{mittag1}), proves Lemma~\ref{8.2}.
\EndDim

%
\begin{lem} \label{2half}
The Mittag-Leffler functions $E_{1,\frac12}$ and $E_{\frac12,\frac12}$
are related to the error function according to
\begin{align*}
E_{1,\frac12}(x)
&=
\begin{cases}
\displaystyle{\frac{1}{\sqrt \pi} +\sqrt x \,e^x \mathrm{Erf} (\sqrt x\,)} &\mbox{for } x\ge0,
\\[2ex]
\displaystyle{\frac{1}{\sqrt \pi} -\sqrt{|x|}\,e^x \mathrm{Erf} (\sqrt{|x|}\,)} &\mbox{for } x\le0,
\end{cases}
\\
E_{\frac12,\frac12} (x)
&=
\frac{1}{\sqrt \pi} + x e^{x^2}+|x| e^{x^2}
\mathrm{Erf} (|x|) \quad \mbox{ for } x \in \mathbb{R}.
\end{align*}
In particular,
$$
E_{\frac12 ,\frac12} (x)=\frac{1}{\sqrt \pi}+x e^{x^2}
\mathrm{Erfc}(|x|)\quad \mbox{ for } x \le 0.
$$
\end{lem}
%
\Dim
Recall that
$\mathrm{Erf} (x)=\frac{2}{\sqrt \pi} \int_0^x e^{-t^2}\,\mathrm{d}t=\frac{1}{\sqrt \pi}
\int_0^{x^2} \frac{e^{-u}}{\sqrt u} \,\mathrm{d}u$ and
$\mathrm{Erfc}(x)=1-\mathrm{Erf}(x)$. By~(\ref{mittag1}), we have,
e.g., for $x \ge 0$, that
$$
E_{1,\frac12}(x)= \frac{1}{\sqrt \pi} \left(1+\sqrt x \,e^x
\int_0^x \frac{e^{-u}}{\sqrt u} \,\mathrm{d}u\right)=\frac{1}{\sqrt \pi}
+\sqrt x \,e^x \mathrm{Erf} (\sqrt x).
$$
On the other hand, we have, for $x \le 0$,
$$
E_{\frac12,\frac12} (x)=\sum_{r=0}^{\infty} \frac{x^r}{\Gamma
\!\left(\vphantom{\frac aN}\right.\!\!\frac{r+1}{2}\!\!\left.\vphantom{\frac aN}\right)}
=\sum_{\ell=0}^{\infty} \frac{x^{2\ell}}{\Gamma\!\left(\vphantom{\frac aN}\right.\!\!
\ell+\frac{1}{2}\!\!\left.\vphantom{\frac aN}\right)}
+\sum_{\ell=0}^{\infty} \frac{x^{2\ell+1}}{\Gamma(\ell+1)}
=E_{1,\frac12}(x^2)+x e^{x^2}
$$
from which we immediately extract the aforementioned representation.
\EndDim

%
\begin{lem}\label{secondintlem}
For $\alpha\neq 0$, the following equality holds:
$$
\int_0^{\infty} \frac{e^{-w^2}}{w^2+\alpha^2}\,\mathrm{d}w
=\frac{\pi}{2 \alpha} \,e^{\alpha^2} \mathrm{Erfc}(\alpha).
$$
\end{lem}
%
\Dim
Put, for $\alpha\neq 0$,
$$
F(\alpha)=\alpha\int_0^{\infty} \frac{e^{-w^2}}{w^2+\alpha^2}\,\mathrm{d}w=
\int_0^{\infty} \frac{e^{-\alpha^2 w^2}}{w^2+1}\,\mathrm{d}w.
$$
We plainly have
\begin{align*}
F'(\alpha)
&=
-2\alpha \int_0^{\infty} \frac{w^2}{w^2+1}\,e^{-\alpha^2 w^2} \,\mathrm{d}w
=-2\alpha \int_0^{\infty} \left(1-\frac{1}{w^2+1}\right)e^{-\alpha^2 w^2} \,\mathrm{d}w
\\
&=
-2\alpha \int_0^{\infty}e^{-\alpha^2 w^2} \,\mathrm{d}w+2\alpha \,F(\alpha).
\end{align*}
So, the function $F$ solves the differential equation
$F'(\alpha)=2 \alpha F(\alpha)-\sqrt{\pi}$ with initial value
\mbox{$F(0)=\frac{\pi}{2}$}. We deduce that $F(\alpha)$ has the form
$F(\alpha)=G(\alpha)\,e^{\alpha^2}$ where the unknown function $G$
satisfies \mbox{$G'(\alpha)=-\sqrt{\pi} \,e^{-\alpha^2}$} and $G(0)=\frac{\pi}{2}$.
This implies that $G(\alpha)=\sqrt{\pi} \int_{\alpha}^{\infty}e^{-\xi^2}\,\mathrm{d}\xi
=\frac \pi 2 \,\mathrm{Erfc}(\alpha)$ and then
$F(\alpha)=\frac{\pi}{2} \,e^{\alpha^2} \mathrm{Erfc}(\alpha)$.
The proof of Lemma~\ref{secondintlem} is established.
\EndDim

%
\begin{lem} \label{secondint}
The following identity holds: for $\lambda>0$, $s>0$ and $x>0$,
$$
\int_0^s \frac{e^{-\lambda\sigma -\frac{x^2}{4\sigma}}}{\sigma^{3/2}}
\left(\int_{\sqrt{\lambda (s-\sigma)}}^{\infty} e^{-\xi^2} \,\mathrm{d}\xi\right)\mathrm{d}\sigma
=\frac{\pi}{x} \,e^{\sqrt \lambda \,x}\,
\mathrm{Erfc}\!\left(\frac{x}{2 \sqrt s}+\sqrt{\lambda s}\right)
$$
\end{lem}
%
\Dim
Put $I=\int_0^s \frac{e^{-\lambda\sigma-\frac{x^2}{4\sigma}}}{\sigma^{3/2}}
\left(\int_{\sqrt{\lambda (s-\sigma)}}^{\infty} e^{-\xi^2} \,\mathrm{d}\xi\right)\mathrm{d}\sigma$.
With the change of variables $\xi=\sqrt{\lambda s z}$ and $v=z+\sigma/s$, we have
$$
I=\frac{\sqrt{\lambda s}}{2}  \int_0^s \frac{e^{-\lambda\sigma
-\frac{x^2}{4\sigma}}}{\sigma^{3/2}} \left( \int_{(s-\sigma)/s}^{\infty}
\frac{e^{-\lambda s z}}{\sqrt{z}} \,\mathrm{d}z\right)\mathrm{d}\sigma
=\frac{\sqrt{\lambda s}}{2} \int_0^s \frac{e^{-\frac{x^2}{4\sigma}}}{\sigma^{3/2}}
\left(\int_{1}^{\infty} \frac{e^{-\lambda s v}}{\sqrt{v-\frac{\sigma}{s}}} \,\mathrm{d}v
\right)\mathrm{d}\sigma.
$$
With the new change of variable $u=s/\sigma$, it becomes
\begin{align*}
I
&=
\frac{\sqrt{\lambda s}}{2} \int_{1}^{\infty} \left(\frac{u}{s}\right)^{\!3/2}
\frac{s}{u^2} \,e^{-\frac{x^2}{4s}u} \left(\int_{1}^{\infty}
\frac{e^{-\lambda s v}}{\sqrt{v-\frac 1u}} \,\mathrm{d}v\right)\mathrm{d}u =\frac{\sqrt \lambda}{2}
\int_{1}^{\infty} \!\!\!\int_{1}^{\infty} \frac{e^{-\left(\frac{x^2}{4s}\right)u
-\lambda s v}}{\sqrt{uv-1}} \,\mathrm{d}u \,\mathrm{d}v
\\
&=
\frac{\sqrt \lambda}{2} \int_{1}^{\infty} \!\!\!\int_{1}^{\infty}
\frac{e^{-\frac a 2 u-\frac b 2 v}}{\sqrt{uv-1}} \,\mathrm{d}u \,\mathrm{d}v
\end{align*}
where we set $a=\frac{x^2}{2s}$ and $b=2 \lambda s$. With the change of
variable $(v,w)=(v,\frac a2 u+ \frac b2 v)$, this gives
\begin{align*}
I
&=
\frac{\sqrt{\lambda}}{\sqrt{ab}} \int_{w=\frac{a+b}{2}}^{w=\infty} \int_{v=1}^{v=\frac{2w-a}{b}}
\frac{e^{-w}}{\sqrt{\frac 2 b \,v w-v^2-\frac a b}} \,\mathrm{d}v \,\mathrm{d}w=\frac{1}{x}
\int_{\frac{a+b}{2}}^{\infty} e^{-w} \left(\int_{1}^{\frac{2w-a}{b}}
\frac{\mathrm{d}v}{\sqrt{\frac{w^2-ab}{b^2}-(v-\frac{w}{b})^2}}\right) \mathrm{d}w
\\
&=
\frac{1}{x}\int_{\frac{a+b}{2}}^{\infty} \left(\arcsin \frac{w-a}{\sqrt{w^2-ab}}
+ \arcsin \frac{w-b}{\sqrt{w^2-ab}}\right) e^{-w} \,\mathrm{d}w.
\end{align*}
Let us integrate by parts this last integral. First, we have
$$
\frac{\mathrm{d}}{\mathrm{d}w}\left(\arcsin \frac{w-a}{\sqrt{w^2-ab}}\right)
=\frac{\sqrt{a}\,(w-b)}{(w^2-ab)\sqrt{2w-a-b}},
$$
and then,
\begin{align*}
I
&=
\frac{1}{x} \left( \left[-\left(\arcsin \frac{w-a}{\sqrt{w^2-ab}}
+ \arcsin \frac{w-b}{\sqrt{w^2-ab}}\right)e^{-w}\right]_{\frac{a+b}{2}}^{\infty}
+\int_{\frac{a+b}{2}}^{\infty} \frac{\sqrt{a} \,(w-b)
+\sqrt b \,(w-a)}{(w^2-ab)\sqrt{2w-a-b}}\,e^{-w}\,\mathrm{d}w\right)
\\
&=
\frac{\sqrt{a}+\sqrt{b}}{x} \int_{\frac{a+b}{2}}^{\infty}
\frac{e^{-w}}{(w+\sqrt{ab})\sqrt{2w-(a+b)}}\,\mathrm{d}w=\frac{\sqrt{a}+\sqrt{b}}{x \sqrt 2}\,
e^{-\frac{a+b}{2}}\int_0^{\infty} \frac{e^{-w}}{(w+\frac{a+b}{2}+\sqrt{ab})\sqrt{w}}\,\mathrm{d}w
\\
&=
\sqrt 2 \,\frac{\sqrt{a}+\sqrt{b}}{x}\,e^{-\frac{a+b}{2}}
\int_0^{\infty} \frac{e^{-w^2}}{w^2+\left(\frac{\sqrt a+\sqrt b}{\sqrt 2}\right)^{\!2}} \,\mathrm{d}w.
\end{align*}
As a byproduct, in view of Lemma~\ref{secondintlem}, we have
$$
I=\sqrt 2 \,\frac{\sqrt{a}+\sqrt{b}}{x}
\,\frac{\pi e^{-\frac{a+b}{2}}}{\sqrt 2 (\sqrt a+\sqrt b)}\,
e^{\frac{a+b}{2}+\sqrt{ab}} \,\mathrm{Erfc}\!\left(\frac{\sqrt a+\sqrt b}{\sqrt 2}\right)
=\frac{\pi}{x}\,e^{\sqrt \lambda \,x} \,\mathrm{Erfc}\!\left(\frac{x}{2 \sqrt s}+\sqrt{\lambda s}\right)\!.
$$
Lemma~\ref{secondint} is proved.
\EndDim

%
\begin{lem} \label{brom}
For any integer $m\le N-1$ and any $x\ge 0$,
$$
\int_0^{\infty} e^{-\lambda u} I_{j,m}(u;x)\,\mathrm{d}u
=\lambda^{-\frac{m}{N}}e^{-\theta_j\!\!\sqrt[N]{\lambda}\,x}.
\eqno{(\ref{set24})}
$$
\end{lem}
%
\Dim
This formula is proved in~\cite{2007} for $0\le m\le N-1$. To prove that it holds
true also for negative $m$, we directly compute the Laplace transform
of $I_{j,m}(u;x)$. We have
$$
\int_0^{\infty} e^{-\lambda u} I_{j,m}(u;x) \,\mathrm{d}u
=\frac{N i}{2 \pi} \left(e^{-i \frac{m}{N}\pi } \int_0^{\infty} \frac{\xi^{N-m-1}}{\xi^N+\lambda}\,
e^{-\theta_j e^{i \frac{\pi}{N}} x \xi} \,\mathrm{d}\xi -e^{i \frac{m}{N}\pi } \int_0^{\infty}
\frac{\xi^{N-m-1}}{\xi^N+\lambda}\, e^{-\theta_j e^{-i \frac{\pi}{N}} x \xi} \,\mathrm{d}\xi\right)\!.
$$
Let us integrate the function $H:z \to \frac{z^{M-1}}{z^N+ \lambda} \,e^{-a z}$
for fixed $a$ and $M$ such that $\Re(a)>0$ and $M>0$ on the contour
$\Gamma_R=\{\rho e^{i \varphi}\in\mathbb{C}: \varphi=0, \rho \in [0,R]\}
\cup \{\rho e^{i \varphi}\in\mathbb{C}: \varphi \in (0,-\frac{2 \pi} N), \rho=R\}
\cup\{\rho e^{i \varphi}\in\mathbb{C}: \varphi=-\frac{2 \pi}N,$ $\rho\in(0,R]\}$.
We get, by residues theorem,
\begin{align*}
-\int_0^{\infty} \frac{z^{M-1}}{z^N+\lambda} \,e^{-a z} \,\mathrm{d}z
+e^{-2i \frac{M}{N}\pi} \int_0^{\infty} \frac{z^{M-1}}{z^N+\lambda}\,
e^{-a e^{-i\frac{2\pi}{N}}z} \,\mathrm{d}z
&=
2i \pi \,\mathrm{Residue}\!\left(H,\!\sqrt[N]{\lambda} \,e^{-i \frac{\pi}{N}}\right)
\\
&=
\frac{2i\pi}{N} \left(\!\!\!\vphantom{\sqrt a}\right.\sqrt[N] \lambda
\!\left.\vphantom{\sqrt a}\right)^{\!M-N} e^{-i \frac{M-N}{N}\pi}
e^{-a \!\sqrt[N] \lambda \,e^{-i\frac{\pi}{N}}}
\\
&=
\frac{2 \pi}{Ni} \,\lambda^{\frac{M}{N}-1} e^{-i\frac{M}{N}\pi}
e^{-a\!\sqrt[N] \lambda \,e^{-i\frac{\pi}{N}}}.
\end{align*}
For $M=N-m$ and $a=\theta_j\,e^{i \frac{\pi}{N}}x$, this yields
$$
\int_0^{\infty} e^{-\lambda u} I_{j;m} (u;x) \,\mathrm{d}u=
-e^{-i\frac{m}{N}\pi} \lambda^{-\frac{m}{N}} \times (-e^{i\frac{m}{N}\pi})
e^{-\theta_j\!\!\sqrt[N] \lambda \,x}=\lambda^{-\frac{m}{N}}e^{-\theta_j\!\!\sqrt[N] \lambda \,x}.
$$
Hence, (\ref{set24}) is valid for $m\le N-1$.
\EndDim

%
\begin{lem} \label{s}
The following identity holds: for any $x\in\mathbb{R}$ and $0<s<t$,
\begin{align}
\lqn{\int_0^{\infty} \xi^2 \sin(x \xi) e^{-t \xi^2} \!\left(\vphantom{\int}\right.\!
\int_0^{\xi} e^{s u^2} \,\mathrm{d}u\!\!\left.\vphantom{\int}\right) \mathrm{d}\xi}
&=
\frac{\sqrt \pi}{8} \frac{2t-s}{t^2 (t-s)^{3/2}} \,x e^{-\frac{x^2}{4(t-s)}}
-\frac{\pi}{8 \sqrt s \,t^{3/2}} \left(1-\frac{x^2}{2t}\right) e^{-\frac{x^2}{4t}}\,
\mathrm{Erf}\!\left(-\frac12 \sqrt{\frac{s}{t(t-s)}} \,x\right)\!.
\label{identity-s}
\end{align}
\end{lem}
%
\Dim
Set
\begin{align*}
A
&=
\int_0^{\infty} \xi^2 \sin(x \xi) e^{-t \xi^2} \!\left(\vphantom{\int}\right.\!
\int_0^{\xi} e^{s u^2} \,\mathrm{d}u\!\!\left.\vphantom{\int}\right) \mathrm{d}\xi,
\\
B
&=
\frac{\sqrt \pi}{8} \frac{2t-s}{t^2 (t-s)^{3/2}} \,x e^{-\frac{x^2}{4(t-s)}}
-\frac{\pi}{8\sqrt s \,t^{3/2}} \left(1-\frac{x^2}{2t}\right) e^{-\frac{x^2}{4t}}
\mathrm{Erf}\!\left(-\frac{1}{2} \sqrt{\frac{s}{t(t-s)}} \,x\right)\!.
\end{align*}
Using the expansion of the sine function, we get
$$
A=\sum_{n=0}^{\infty}(-1)^n \frac{x^{2n+1}}{(2n+1)!} \int_0^{\infty} \left(\xi^{2n+3}
e^{-t\xi^2}\right)\!\!\left(\vphantom{\int}\right.\!
\int_0^{\xi} e^{s u^2} \,\mathrm{d}u\!\!\left.\vphantom{\int}\right)\mathrm{d}\xi.
$$
In order to integrate by parts this last integral, we search for a primitive
of $\xi^{2n+3} e^{-t \xi^2}$. With the change of variable $\zeta=t\xi^2$, we have
$$
\int \xi^{2n+3} e^{-t \xi^2}\,\mathrm{d}\xi
=\frac{1}{2t^{n+2}} \int \zeta^{n+1}e^{-\zeta} \,\mathrm{d}\zeta
=-\frac{(n+1)!}{2t^{n+2}} \left(\,\sum_{k=0}^{n+1} \frac{\zeta^k}{k!}\right) e^{-\zeta}
=-\frac{(n+1)!}{2t^{n+2}} \sum_{k=0}^{n+1} \frac{t^k \xi^{2 k}}{k!} \, e^{-t\xi^2}.
$$
Then,
\begin{align}
A
&=
\frac12 \sum_{n=0}^{\infty} (-1)^n \frac{(n+1)!}{(2n+1)! \,t^{n+2}}
\left[\,\sum_{k=0}^{n+1} \frac{t^k}{k!}
\int_0^{\infty} \xi^{2 k} e^{-t \xi^2} \,\frac{\mathrm{d}}{\mathrm{d}\xi}
\left(\vphantom{\int}\right.\!\int_0^{\xi} e^{s u^2} \,\mathrm{d}u\!\!
\left.\vphantom{\int}\right)\mathrm{d}\xi\right] x^{2n+1}
\nonumber\\
&=
\frac{x}{4 t^2 \sqrt{t-s}} \sum_{n=0}^{\infty} (-1)^n \frac{(n+1)!}{(2n+1)! \,t^n}
\left[\,\sum_{k=0}^{n+1} \frac{\Gamma\!\left(\vphantom{\frac aN}\right.\!\! k+\frac12
\!\!\left.\vphantom{\frac aN}\right)}{k!} \left(\frac{t}{t-s}\right)^{\!k}\right]x^{2n}.
\label{expression1}
\end{align}
On the other hand, by expanding
$\mathrm{Erf}(\xi)=\frac{2}{\sqrt \pi} \int_0^{\xi} \sum_{n=0}^{\infty}
(-1)^n \frac{u^2}{n!} \,\mathrm{d}u=\frac{2}{\sqrt \pi}
\sum_{n=0}^{\infty} (-1)^n \frac{\xi^{2n+1}}{(2n+1)n!}$, we get
\begin{align}
B
&=
\frac{\sqrt \pi}{8} \frac{2t-s}{t^2 (t-s)^{3/2}} \,x e^{-\frac{x^2}{4(t-s)}}
\nonumber\\
&\hphantom{=\,}
-\frac{\pi}{8\sqrt s \,t^{3/2}} \left(\,\sum_{n=0}^{\infty} \frac{(-1)^n x^{2n}}{n! (4t)^n}
-2\sum_{n=0}^{\infty} \frac{(-1)^n x^{2n+2}}{n! (4t)^{n+1}}\right) \frac{2 x}{\sqrt \pi}
\sum_{n=0}^{\infty} \left(\frac12 \sqrt{\frac{s}{t(t-s)}}\,\right)^{\!2n+1}
\frac{(-1)^{n+1}x^{2n}}{(2n+1)n!}
\nonumber\\
\lqn{}
&=
\frac{\sqrt \pi}{8t^2} \frac{x}{\sqrt{t-s}} \left[\frac{2t-s}{t-s}
\sum_{n=0}^{\infty} \frac{(-1)^n x^{2n}}{n! (4(t-s))^n}\right.
\nonumber\\
&\hphantom{=\,}
\left.+\left(\,\sum_{p=0}^{\infty}(-1)^p \frac{2p+1}{p!} \frac{x^{2p}}{(4t)^p}\right) \!\!
\left(\,\sum_{q=0}^{\infty}\frac{(-1)^q}{(2q+1)q!}
\left(\frac{s}{4t(t-s)}\right)^{\!q} x^{2q}\right)\right]\!.
\label{somme1}
\end{align}
The computation of the above product of series can be carried out as follows:
\begin{align}
\lqn{\left(\,\sum_{p=0}^{\infty}(-1)^p \frac{2p+1}{p!} \frac{x^{2p}}{(4t)^p}\right) \!\!
\left(\,\sum_{q=0}^{\infty}\frac{(-1)^q}{(2q+1)q!} \left(\frac{s}{4t(t-s)}\right)^{\!q} x^{2q}\right)}
&=
\sum_{n=0}^{\infty} \left[\vphantom{\sum_p}\right.
\sum_{\substack{p,q\ge 0:\\ p+q=n}} \frac{2p+1}{(2q+1)p!q!}
\left(\frac{s}{t-s}\right)^{\!q} \!\left.\vphantom{\sum_p}\right]\frac{(-1)^nx^{2n}}{(4t)^n}
\nonumber\\
&=
\sum_{n=0}^{\infty} \left[\,\sum_{q=0}^n \frac{2n-2q+1}{2q+1} \binom{n}{q}\!\!
\left(\frac{s}{t-s}\right)^{\!q} \right]\frac{(-1)^nx^{2n}}{n!(4t)^n}
\nonumber\\
&=
\sum_{n=0}^{\infty} \left[(2n+2)\sum_{q=0}^n \frac{\binom{n}{q}}{2q+1}
\left(\frac{s}{t-s}\right)^{\!q} -\sum_{q=0}^n \binom{n}{q}\!\!
\left(\frac{s}{t-s}\right)^{\!q} \right]\frac{(-1)^nx^{2n}}{n!(4t)^n}
\nonumber\\
&=
\sum_{n=0}^{\infty} \left[(2n+2)\sum_{q=0}^n \frac{\binom{n}{q}}{2q+1}
\left(\frac{s}{t-s}\right)^{\!q} \right]\frac{(-1)^nx^{2n}}{n!(4t)^n}
-\sum_{n=0}^{\infty} \frac{(-1)^nx^{2n}}{n!(4(t-s))^n}.
\label{somme2}
\end{align}
By inserting~(\ref{somme2}) into~(\ref{somme1}), we derive
\begin{equation}
B=\frac{\sqrt \pi}{8 t^2} \frac{x}{\sqrt{t-s}} \sum_{n=0}^{\infty}
\left[\left(\frac{t}{t-s}\right)^{\!n+1}\! + \sum_{q=0}^n \frac{\binom{n}{q}}{2q+1}
\left(\frac{s}{t-s}\right)^{\!q}\right] \frac{(-1)^nx^{2n}}{n!(4t)^n}.
\label{somme3}
\end{equation}
Writing now $(\frac{s}{t-s})^q=(\frac{t}{t-s}-1)^q
=\sum_{k=0}^q (-1)^{q-k}\binom{q}{k} (\frac{t}{t-s})^k$, we have
\begin{equation}
\sum_{q=0}^n \frac{\binom{n}{q}}{2q+1}\left(\frac{s}{t-s}\right)^{\!q}
=\sum_{k=0}^n (-1)^k \left(\vphantom{\frac{\binom{q}{k}}{q}}\right.\!
\sum_{q=k}^n (-1)^q \left.\frac{\binom{q}{k} \binom{n}{q}}{2q+1}\right)
\!\!\left(\frac{t}{t-s}\right)^{\!k}.
\label{somme3bis}
\end{equation}
The sum with respect to $q$ can be calculated as follows:
\begin{align}
\sum_{q=k}^{n} (-1)^q \frac{\binom{q}{k} \binom{n}{q}}{2q+1}
&=
\binom{n}{k} \sum_{q=k}^{n} (-1)^q \frac{\binom{n-k}{q-k}}{2q+1}
=(-1)^k\binom{n}{k} \sum_{q=0}^{n-k} (-1)^q \frac{\binom{n-k}{q}}{2q+2k+1}
\nonumber\\
&=
(-1)^k\binom{n}{k} \int_0^1 \sum_{q=0}^{n-k} (-1)^q \binom{n-k}{q} \xi^{2q+2k}\,\mathrm{d}\xi
\nonumber\\
&=
(-1)^k\binom{n}{k} \int_0^1 (1-\xi^2)^{n-k}\xi^{2k}\,\mathrm{d}\xi
=(-1)^k\frac12\binom{n}{k} \int_0^1 (1-\xi)^{n-k}\xi^{k-1/2}\,\mathrm{d}\xi
\nonumber\\
&=
(-1)^k\frac{1}{2} \binom{n}{k} B\!\left(n-k+1,k+\frac12\right)
=(-1)^k\frac{n!\,\Gamma\!\left(\vphantom{\frac aN}\right.\!\!
k+\frac12 \!\!\left.\vphantom{\frac aN}\right)}{2k!\,\Gamma
\!\left(\vphantom{\frac aN}\right.\!\!n+\frac32 \!\!\left.\vphantom{\frac aN}\right)}.
\label{somme4}
\end{align}
Plugging~(\ref{somme4}) into~(\ref{somme3bis}), and next~(\ref{somme3bis})
into~(\ref{somme3}), we obtain
\begin{align}
B
&=
\frac{\sqrt \pi}{8 t^2} \frac{x}{\sqrt{t-s}} \sum_{n=0}^{\infty}
\frac{(-1)^n}{n! (4t)^n} \left[\left(\frac{t}{t-s}\right)^{\!n+1}+\frac{(n+1)!}{\Gamma
\!\left(\vphantom{\frac aN}\right.\!\!n+\frac32\!\!\left.\vphantom{\frac aN}\right)}
\sum_{k=0}^n \frac{\Gamma\!\left(\vphantom{\frac aN}\right.\!\!k+\frac12
\!\!\left.\vphantom{\frac aN}\right)}{k!} \left(\frac{t}{t-s}\right)^{\!k}\right] x^{2n}
\nonumber\\
&=
\frac{\sqrt \pi}{8 t^2} \frac{x}{\sqrt{t-s}} \sum_{n=0}^{\infty} \frac{(-1)^n}{n! (4t)^n}
\frac{(n+1)!}{\Gamma\!\left(\vphantom{\frac aN}\right.\!\!n+\frac32
\!\!\left.\vphantom{\frac aN}\right)} \left[\,\sum_{k=0}^{n+1} \frac{\Gamma
\!\left(\vphantom{\frac aN}\right.\!\!k+\frac12\!\!\left.\vphantom{\frac aN}\right)}{k!}
\left(\frac{t}{t-s}\right)^{\!k}\right] x^{2n}
\nonumber\\
&=
\frac{x}{4 t^2 \sqrt{t-s}} \sum_{n=0}^{\infty} (-1)^n \frac{(n+1)!}{(2n+1)! t^n}
\left[\,\sum_{k=0}^{n+1} \frac{\Gamma\!\left(\vphantom{\frac aN}\right.\!\!
k+\frac12 \!\!\left.\vphantom{\frac aN}\right)}{k!} \left(\frac{t}{t-s}\right)^{\!k}\right]x^{2n}
\label{expression2}
\end{align}
where we used in the last equality
$\Gamma\!\left(\vphantom{\frac aN}\right.\!\!n+\frac32\!\!\left.\vphantom{\frac aN}\right)\!
=\frac{\sqrt \pi}{2\cdot 4^n} \frac{(2n+1)!}{n!}$.
In view of~(\ref{expression1}) and~(\ref{expression2}), we see that
both members of~(\ref{identity-s}) are equal: $A=B$.
\EndDim

%
\begin{lem}\label{ss}
We have, for $x \le 0$,
\begin{align*}
\lqn{\int_0^{\infty} y(y-x) \,e^{-\frac{y^2}{4 s}-\frac{(y-x)^2}{4(t-s)}} \,\mathrm{d}y}
&=
-2\,\frac{s(t-s)^2}{t^2} \,x e^{-\frac{x^2}{4(t-s)}}+2\sqrt \pi\,\frac{s^{3/2} (t-s)^{3/2}}{t^{3/2}}
\left(1-\frac{x^2}{2t}\right) e^{-\frac{x^2}{4t}}\,
\mathrm{Erfc}\!\left(-\frac12 \,\sqrt{\frac{s}{t(t-s)}}\,x\right)\!,
\end{align*}
and, for $x \ge 0$,
\begin{align*}
\lqn{\int_0^{\infty} y(y+x) \,e^{-\frac{(y+x)^2}{4 s}-\frac{y^2}{4(t-s)}} \,\mathrm{d}y}
&=
2\,\frac{s^2(t-s)}{t^2} \,x e^{-\frac{x^2}{4s}}+2\sqrt \pi\,
\frac{s^{3/2} (t-s)^{3/2}}{t^{3/2}}\left(1-\frac{x^2}{2t}\right) e^{-\frac{x^2}{4t}}\,
\mathrm{Erfc}\!\left(\frac12 \,\sqrt{\frac{t-s}{s\,t}}\,x\right)\!.
\end{align*}
\end{lem}
%
\Dim
We only produce the proof for $x\le 0$ since the case $x\ge 0$ is quite similar.
We have
\begin{align*}
\lqn{\int_0^{\infty} y(y-x) \,e^{-\frac{y^2}{4 s}-\frac{(y-x)^2}{4(t-s)}} \,\mathrm{d}y}
&=
e^{-\frac{x^2}{4t}} \int_0^{\infty} y(y-x) \,e^{-\frac{t}{4s(t-s)} (y-\frac{sx}{t})^2} \,\mathrm{d}y
\\
&=
e^{-\frac{x^2}{4t}} \int_{-sx/t}^{\infty} \left(y+\frac{sx}{t}\right)\!\!
\left(y-\frac{(t-s)x}{t}\right) e^{-\frac{t}{4s(t-s)}y^2}\,\mathrm{d}y
\\
&=
e^{-\frac{x^2}{4t}} \int_{-sx/t}^{\infty} \left(y^2-\frac{s(t-s)}{t^2}\,x^2\right)
e^{-\frac{t}{4s(t-s)}y^2} \,\mathrm{d}y+\frac{2s-t}{t}\,x e^{-\frac{x^2}{4t}}
\int_{-sx/t}^{\infty} y e^{-\frac{t}{4s(t-s)}y^2} \,\mathrm{d}y
\\
&=
e^{-\frac{x^2}{4t}} \int_{-sx/t}^{\infty} y^2 e^{-\frac{t}{4s(t-s)}y^2} \,\mathrm{d}y
-\frac{s(t-s)}{t^2} \,x^2 e^{-\frac{x^2}{4t}}
\int_{-sx/t}^{\infty} e^{-\frac{t}{4s(t-s)}y^2}\,\mathrm{d}y
+ 2\,\frac{s(t-s)(2s-t)}{t^2} \, x e^{-\frac{x^2}{4(t-s)}}.
\end{align*}
Integrating by parts, we observe that
$$
\int_{-sx/t}^{\infty} y^2 e^{-\frac{t}{4s(t-s)}y^2} \,\mathrm{d}y
=2\,\frac{s^2(t-s)}{t^2}\,xe^{-\frac{sx^2}{4t(t-s)}}
+2\,\frac{s(t-s)}{t} \int_{-sx/t}^{\infty} e^{-\frac{t}{4s(t-s)}y^2}\,\mathrm{d}y,
$$
we finally get the result.
\EndDim


\begin{thebibliography}{3}

\bibitem{as}
Abramowitz, M. and Stegun, I.-A.: Handbook of mathematical functions
with formulas, graphs, and mathematical tables. Dover Publications, 1992.

\bibitem{bho}
Beghin, L., Hochberg, K.-J. and Orsingher, E.: Conditional maximal
distributions of processes related to higher-order heat-type equations.
Stochastic Process. Appl. 85 (2000), no.~2, 209--223.

\bibitem{bo}
Beghin, L. and Orsingher, E.: The distribution of the local time for
``pseudoprocesses'' and its connection with fractional diffusion equations.
Stochastic Process. Appl. 115 (2005), 1017--1040.

\bibitem{bor}
Beghin, L., Orsingher, O. and Ragozina, T.: Joint distributions of the
maximum and the process for higher-order diffusions.
Stochastic Process. Appl. 94 (2001), no.~1, 71--93.

\bibitem{bohr}
Bohr, H.: \"Uber die gleichm\"assige Konvergenz Dirichletscher Reihen.
J. f\"ur Math. 143 (1913), 203--211.

\bibitem{bs}
Borodin, A.-N. and Salminen, P.: Handbook of Brownian motion---facts and formulae.
First edition. Probability and its Applications. Birkh\"auser Verlag, 1996.

\bibitem{erdelyi}
Erd\'elyi, A., Ed.: Higher transcendental functions, vol. III,
Robert Krieger Publishing Company, Malabar, Florida, 1981.

\bibitem{hoch}
Hochberg, K.-J.: A signed measure on path space related to Wiener measure.
Ann. Probab. 6 (1978), no.~3, 433--458.

\bibitem{ho}
Hochberg, K.-J. and Orsingher, E.: The arc-sine law and its analogs
for processes governed by signed and complex measures.
Stochastic Process. Appl. 52 (1994), no.~2, 273--292.

\bibitem{kry}
Krylov, V. Yu.: Some properties of the distribution corresponding to the equation
$\frac{\partial u}{\partial t}=(-1)^{q+1}\frac{\partial^{2q} u}{\partial^{2q} x}$.
Soviet Math. Dokl. 1 (1960), 760--763.

\bibitem{2003}
Lachal, A.: Distribution of sojourn time, maximum and minimum for
pseudo-processes governed by higer-order heat-typer equations.
Electron. J. Probab. 8 (2003), paper no. 20, 1--53.

\bibitem{2006}
Lachal, A.: Joint law of the process and its maximum, first hitting time and place of
a half-line  for the pseudo-process driven by the equation
$\frac{\partial}{\partial t}=\pm\frac{\partial^N}{\partial x^N}$.
C. R. Acad. Sci. Paris, S\'er.~I 343 (2006), no.~8, 525--530.

\bibitem{2007}
Lachal, A.: First hitting time and place, monopoles and multipoles for
pseudo-procsses driven by the equation
$\frac{\partial}{\partial t}=\pm \frac{\partial^N}{\partial x^N}$.
Electron. J. Probab. 12 (2007), 300--353.

\bibitem{2008}
Lachal, A.: First hitting time and place for the pseudo-process driven
by the equation $\frac{\partial}{\partial t}=\pm\frac{\partial^n}{\partial x^n}$
subject to a linear drift. Stoch. Proc. Appl. 118 (2008), 1--27.

\bibitem{ns}
Nakajima, T. and Sato, S.: On the joint distribution of the first hitting time and the first
hitting place to the space-time wedge domain of a biharmonic pseudo process.
Tokyo J. Math. 22 (1999), no.~2, 399--413.

\bibitem{no}
Nikitin, Ya. Yu. and Orsingher, E.: On sojourn distributions of processes related to some higher-order
heat-type equations. J. Theoret. Probab. 13 (2000), no.~4, 997--1012.

\bibitem{nish1}
Nishioka, K.: Monopole and dipole of a biharmonic pseudo process.
Proc. Japan Acad. Ser. A 72 (1996), 47--50.

\bibitem{nish2}
Nishioka, K.: The first hitting time and place of a half-line by a biharmonic pseudo process.
Japan J. Math. 23 (1997), 235--280.

\bibitem{nish3}
Nishioka, K.: Boundary conditions for one-dimensional biharmonic pseudo process.
Electronic Journal of Probability 6 (2001), paper no.~13, 1--27.

\bibitem{ors}
Orsingher, E.: Processes governed by signed measures connected with third-order ``heat-type''
equations. Lithuanian Math. J. 31 (1991), no.~2, 220--231.

\bibitem{spitzer}
Spitzer, F.: A combinatorial lemma and its application to probability theory.
Trans. Amer. Math. Soc. 82 (1956), 323--339.




\end{thebibliography}
\end{document}